\pdfoutput=1
\RequirePackage{ifpdf}
\ifpdf 
\documentclass[pdftex]{sigma}
\else
\documentclass{sigma}
\fi

\numberwithin{equation}{section}

\newtheorem{Theorem}{Theorem}[section]

\newtheorem{Lemma}[Theorem]{Lemma}

 { \theoremstyle{definition}
\newtheorem{Definition}[Theorem]{Definition}

\newtheorem{Example}[Theorem]{Example}
\newtheorem{Remark}[Theorem]{Remark} }

\usepackage[mathscr]{eucal}
\usepackage{mathtools}
\usepackage{array}
\usepackage[thinlines]{easytable}
\usepackage{mathrsfs}
\usepackage[all]{xy}

\usepackage{longtable, booktabs}

\usepackage{todonotes}
\usepackage{tikz}
\usetikzlibrary{matrix}
\usetikzlibrary{arrows}

\usepackage{multicol}
\newcommand{\cQ}{\mathcal{Q}}

\newcommand{\cA}{\mathcal{A}}
\newcommand{\cB}{\mathcal{B}}

\newcommand{\CC}{\mathbb{C}}
\newcommand{\ZZ}{\mathbb{Z}}
\newcommand{\PP}{\mathbb{P}}
\newcommand{\QQ}{\mathbb{Q}}

\DeclareMathOperator{\SL}{SL}
\DeclareMathOperator{\GL}{GL}

\DeclareMathOperator{\age}{age}

\DeclareMathOperator{\Fix}{Fix}

\usepackage[labelsep=period]{caption}

\begin{document}

\allowdisplaybreaks

\newcommand{\arXivNumber}{1812.06200}

\renewcommand{\PaperNumber}{059}

\FirstPageHeading

\ShortArticleName{Mirror Symmetry for Nonabelian Landau--Ginzburg Models}

\ArticleName{Mirror Symmetry for Nonabelian Landau--Ginzburg\\ Models}

\Author{Nathan PRIDDIS~$^\dag$, Joseph WARD~$^\ddag$ and Matthew M.~WILLIAMS~$^\S$}

\AuthorNameForHeading{N.~Priddis, J.~Ward and M.M.~Williams}

\Address{$^\dag$~Brigham Young University, USA}
\EmailD{\href{mailto:priddis@math.byu.edu}{priddis@math.byu.edu}}

\Address{$^\ddag$~University of Utah, USA}
\EmailD{\href{mailto:josephwardjoseph@gmail.com}{josephwardjoseph@gmail.com}}

\Address{$^\S$~Colorado State University, USA}
\EmailD{\href{mailto:Matthew.Williams@colostate.edu}{Matthew.Williams@colostate.edu}}

\ArticleDates{Received September 24, 2019, in final form June 12, 2020; Published online June 27, 2020}

\Abstract{We consider Landau--Ginzburg models stemming from groups comprised of non-diagonal symmetries, and we describe a rule for the mirror LG model. In particular, we present the \emph{non-abelian dual group} $G^\star$, which serves as the appropriate choice of group for the mirror LG model. We also describe an explicit mirror map between the A-model and the B-model state spaces for two examples. Further, we prove that this mirror map is an isomorphism between the untwisted broad sectors and the narrow diagonal sectors for Fermat type polynomials.}

\Keywords{mirror symmetry; Landau--Ginzburg models; Calabi--Yau; nonabelian}

\Classification{14J32; 53D45; 14J81}

\section{Introduction}

Mirror symmetry is most easily explained for Calabi--Yau manifolds. The physics of string theory produces an A-model and a B-model for each Calabi--Yau manifold, and these come in dual pairs. Mirror symmetry essentially says that the A-model for a Calabi--Yau manifold is ``the same'' as the B-model on its mirror dual, meaning they produce the same physics.

The same physics can also be modeled with what is called a Landau--Ginzburg model, which is conjectured to be computationally more efficient. First constructed by physicists in \cite{IV}, Landau--Ginzburg models are built from an \emph{invertible polynomial} $W$ (also called a \emph{potential function}), and a group $G \leq G_{W}^{\max}$ of symmetries of $W$, both of which we describe later. One of the important structures of a Landau--Ginzburg model~-- both for the A-model and B-model~-- is that of a vector space called the \emph{state space}. This can also be given the structure of a~Frobenius algebra or a Frobenius manifold, and comes with invariants in every genus. The \emph{Landau--Ginzburg $($LG$)$ mirror symmetry conjecture} predicts that for an invertible polynomial $W$ with a~group~$G$ of \emph{admissible symmetries} of~$W$, there is a dual polynomial $W^T$ and dual group $G^T$ of symmetries of $W^T$ such that the Landau--Ginzburg A-model for the pair $(W, G)$ is isomorphic to the Landau--Ginzburg B-model for the pair $\big(W^T, G^T\big)$ (see~\cite{ber-hub} or~\cite{KrThesis}). In this article, we will focus on the LG mirror symmetry conjecture only at the level of state spaces.

In the past, mathematicians have primarily studied LG models of pairs $(W, G)$ where $G$ is an abelian group comprised of so-called \emph{diagonal symmetries} (see~\cite{FJRW}). Mirror theorems at various levels of structure have been proved in this setting. For example, in~\cite{KrThesis}, Krawitz showed that the state spaces for a pair $(W,G)$, and its mirror $\big(W^T,G^T\big)$ are isomorphic. In~\cite{FJJS}, the authors showed that the Frobenius algebras associated to $(W,G)$ and its mirror $\big(W^T,G^T\big)$ are isomorphic for a large class of polynomials and groups.

Chiodo and Ruan in~\cite{ChiRu:10} and Chiodo, Iritani, and Ruan in~\cite{CIR:11} proved an LG mirror theorem for genus zero invariants for Fermat polynomials $W$ with a certain group that we will later denote as~$J_W$. This theorem was extended in \cite{PLS:14,PrSh:13} to include all groups~$G$ satisfying a~certain condition called the Calabi--Yau condition. In \cite{Guere}, Gu\'er\'e proved an LG mirror theorem for genus zero invariants for so-called chain potentials, and in \cite{HLSW}, He--Li--Shen--Webb proved a~mirror theorem for all genus invariants for pairs $(W,G)$, with $G$ being the maximal group of diagonal symmetries (see Definition~\ref{def:Gdiag}).

Finally, let us briefly mention the articles \cite{Arte,BCP,CLPS,CP,FPS} wherein the authors show that in certain cases the Landau--Ginzburg mirror symmetry agrees with more geometric versions of mirror symmetry, such as mirror symmetry for K3 surfaces and Borcea--Voisin mirror symmetry.

All of the above mentioned results require $G$ to be a group of diagonal symmetries. There has been much interest in understanding the mirror symmetry for when $G$ is non-abelian, but until now there has not been a clear way to determine the mirror model, since the dual group was only defined when $G$ was a group of diagonal symmetries.

In this paper, we give a conjecture for the \emph{non-abelian dual group} $G^\star$, which extends the Landau--Ginzburg mirror symmetry conjecture to LG models built from non-diagonal symmetries. We describe the construction of the A- and B-model state spaces, and for two examples provide an explicit isomorphism (a mirror map) between them. More generally, we construct a canonical mirror map on certain natural subspaces of the A- and B-model state spaces for particular polynomials $W$ and groups $G$ containing non-diagonal symmetries, and we show this map is an isomorphism of bigraded vector spaces (see Theorem~\ref{thm:mirror}). This generalizes the mirror map for abelian LG models defined by Krawitz in~\cite{KrThesis}.

This construction of $G^\star$ was also discovered independently by Ebeling and Gusein-Zade, as described in \cite{Ebel1,Ebel2,Ebel3}. There, the authors describe the \emph{parity condition} (PC), which is a~condition on a group $G \leq G_W^{\max}$. They conjecture that if this condition is satisfied, then the Milnor fibers associated to $(W, G)$ and $\big(W^T, G^\star\big)$ have the same orbifold Euler characteristic. They prove this conjecture in several cases. The Milnor fiber is closely related to the A-model and B-model subspaces. In contrast to their work, in this article, we will describe the state spaces explicitly and attempt to find an explicit mirror map.

There are still many hurdles for considering structures beyond vector spaces, the foremost being the lack of definition of even a Frobenius product on the B-side. This is a possible direction for future work.

\section{Preliminary definitions}

In this section, we begin by introducing some definitions that will be vital to the construction of the LG A- and B-model state spaces. Recall, these are built from a pair $(W,G)$; we will first describe the potential function $W$ and then the group of symmetries $G$.

\subsection{Invertible polynomials}

We begin by describing necessary conditions on the potential function $W$.

\begin{Definition}
 A polynomial $W\colon \mathbb{C}^N \rightarrow \mathbb{C}$ is \emph{quasihomogeneous} if there exist positive rational numbers $q_1, \dots, q_N$ such that for every $c \in \mathbb{C}^*$, we have
 \[
 W\big(c^{q_1}x_1, \dots, c^{q_N}x_N\big) = cW(x_1,\dots, x_N).
 \]
 The numbers $q_1, \dots, q_N$ are called the \emph{weights} of the polynomial $W$.
\end{Definition}

\begin{Definition}A quasihomogeneous polynomial $W\colon \mathbb{C}^N \rightarrow \mathbb{C}$ is \emph{nondegenerate} if it has an isolated critical point at the origin, and it contains no monomials of the form $x_ix_j$ for $i \neq j$.
\end{Definition}

\begin{Definition} A quasihomogeneous, nondegenerate polynomial is \emph{invertible} if the polynomial has the same number of monomials as variables.
\end{Definition}

The definition of nondegeneracy implies that the weights $q_1, \dots, q_N$ of $W$ are uniquely determined and $q_i \in \big(0, \tfrac{1}{2}\big) \cap \mathbb{Q}$ for all $i$.

\begin{Example}\label{eg:invertible} Consider the polynomial $W\colon \mathbb{C}^4 \rightarrow \mathbb{C}$ defined by $W = x_{1}^4 + x_{2}^4 + x_{3}^4 + x_{4}^4$.
 One can easily check that $W$ is nondegenerate and that $W$ is quasihomogeneous with weights $\big(\tfrac{1}{4}, \tfrac{1}{4}, \tfrac{1}{4}, \tfrac{1}{4}\big)$.

Clearly this choice of weights is unique. Also, we can see that $W$ has four monomials and four variables, hence $W$ is invertible. We will continue to work with this particular polynomial throughout the first several sections of this article.
\end{Example}

\begin{Theorem}[Kreuzer--Skarke \cite{KrSk}]
Any invertible quasihomogeneous polynomial is a Thom--Sebastiani sum of polynomials $($meaning no two of the polynomials share a variable$)$ of one of the following three atomic types:
\begin{enumerate}\itemsep=0pt
 \item[$1)$] {Fermat type:} $x_1^{a_1}$,
 \item[$2)$] {chain type:} $x_1^{a_1}x_2 + x_2^{a_2}x_3 + \dots + x_N^{a_N}$, $N \geq 1$,
 \item[$3)$] {loop type:} $x_1^{a_1}x_2 + x_2^{a_2}x_3 + \dots + x_N^{a_N}x_1$, $N \geq 2$.
\end{enumerate}
In each case, we require $a_i\geq 2$ for all $i$.
\end{Theorem}

\begin{Example}\label{eg:chain}
The polynomial from Example~\ref{eg:invertible} defined by $W = x_1^4 + x_2^4 + x_3^4 + x_4^4$ with weights $\big(\tfrac{1}{4}, \tfrac{1}{4}, \tfrac{1}{4}, \tfrac{1}{4}\big)$ is a Fermat polynomial. An example of a chain polynomial is $ x_1^3x_2 + x_2^2x_3 + x_3^2$ with weights $\big(\tfrac{1}{4}, \tfrac{1}{4}, \tfrac{1}{2}\big)$ and an example of a loop polynomial is $ x_1^2x_2 + x_2^2x_3 + x_3^2x_1$, which has weights $\big(\tfrac{1}{3}, \tfrac{1}{3}, \tfrac{1}{3}\big)$.
\end{Example}

\subsection{Maximal symmetry group}

Now we describe the conditions on the groups $G$ that we will use to construct the Landau--Ginzburg models.

\begin{Definition}[Mukai \cite{Mukai}]\label{def:Gmax}
 Let $W\colon \mathbb{C}^N \rightarrow \mathbb{C}$ be an invertible polynomial with weights $(q_1,\dots,q_N)$. Then the \emph{maximal symmetry group of $W$}, denoted $G_{W}^{\max}$, is defined as follows:
 \begin{gather*}
 G_{W}^{\max} := \{g \in \GL_N(\mathbb{C}) \,|\, (g \cdot W)(x_1, \dots, x_N) = W(x_1, \dots, x_N)
 \text{ and } g_{ij} = 0 \text{ if } q_i \neq q_j\}.
 \end{gather*}
\end{Definition}
The condition $g_{ij} = 0$ if $q_i \neq q_j$ is equivalent to the condition that each $g\in G^{\max}_W$ commutes with the action of $\mathbb{C}^*$ (see~\cite{Mukai}) where $\CC^*$ acts on $(x_1, \dots, x_N)$ by
\[
c \cdot (x_1, \dots, x_N) = \big(c^{q_1}x_1, \dots, c^{q_N}x_N\big).
\]

\begin{Definition}\label{def:Gdiag}
 The \emph{diagonal symmetry group of W} is the group of diagonal linear transformations, defined by
 \[
 G_{W}^{\rm diag} := \big\{ (g_1, \dots, g_N) \in (\mathbb{C}^*)^N \,|\, W(g_1x_1, \dots, g_Nx_N) = W(x_1, \dots, x_N)\big\}.
 \]
\end{Definition}

This definition is the standard definition of diagonal symmetries (see, e.g.,~\cite{FJRW}). Note that~$G_{W}^{\rm diag}$ can be viewed as a subgroup of $G_{W}^{\max}$ via diagonal matrices. It is a standard fact that for $g = (g_1, \dots, g_N) \in G_W^{\rm diag}$ the entries $g_i$ as above are roots of unity (see, e.g.,~\cite{Arte}). For simplicity, we will typically represent these symmetries additively as $N$-tuples of rational numbers as follows:
\[
\big({\rm e}^{2\pi {\rm i} a_1},\dots,{\rm e}^{2 \pi {\rm i} a_N}\big)\leftrightarrow (a_1,\dots,a_N)\in(\QQ/\ZZ)^N.
\]

Furthermore, $G_W^{\rm diag}$ is generated by the entries of the inverse of the \emph{exponent matrix} $A_W$, which we define below (see \cite{Arte,KrThesis}). One can see that the exponential grading operator $j_W = (q_1, \dots, q_N)$ is an element of $G_W^{\rm diag}$, where $q_1, \dots, q_N$ are the weights of~$W$. We denote the group generated by~$j_W$ to be~$J_W$.

Two other important subgroups of $G_{W}^{\max}$ are $\SL_W$ and $\SL_{W}^{\rm diag}$. These are defined as
\begin{gather*}
\SL_W = \SL_N(\mathbb{C}) \cap G_{W}^{\max},\qquad
\SL_{W}^{\rm diag} = \SL_W\cap G^{\rm diag}_W.
\end{gather*}

\begin{Example}\label{eg:J_W dual}
 For $W = x_{1}^4 + x_{2}^4 + x_{3}^4 + x_{4}^4$, we have
 \[J_W = \big\langle \big(\tfrac{1}{4}, \tfrac{1}{4},\tfrac{1}{4}, \tfrac{1}{4}\big) \big\rangle \qquad \text{and} \qquad \SL_{W}^{\rm diag} = \big\langle \big(\tfrac{1}{4}, \tfrac{1}{4},\tfrac{1}{4}, \tfrac{1}{4}\big), \big(\tfrac{2}{4}, \tfrac{1}{4},\tfrac{1}{4}, 0\big), \big(\tfrac{1}{4}, \tfrac{2}{4},\tfrac{1}{4}, 0\big) \big\rangle.\]
\end{Example}

BHK mirror symmetry associates to an LG model $(W,G)$ another LG model $\big(W^T,G^T\big)$, which we work towards next.

\section{Dual polynomials and dual groups}

In this section, we will begin by reviewing the construction of $\big(W^T,G^T\big)$, known as BHK mirror symmetry. This is necessary to understanding the rule for mirror symmetry for nonabelian LG models. Then we will describe the rule for nonabelian LG models.

The following definition was first given by Berglund and H\"ubsch in \cite{ber-hub}.

\begin{Definition}\label{def:exp matrix} Let $W$ be an invertible polynomial. If we write $W = \sum\limits_{i=1}^{N} \prod\limits_{j=1}^{N} x_{j}^{a_{ij}},$ then the associated \emph{exponent matrix} is defined to be $A_W = (a_{ij})$. Notice we have suppressed coefficients of~$W$, as these can be scaled away. The \emph{dual polynomial}~$W^T$ is the invertible polynomial defined by the matrix~$A_{W}^{T}$.
\end{Definition}

\begin{Example}
 For $W = x_{1}^4 + x_{2}^4 + x_{3}^4 + x_{4}^4$, we have
 \[A_W =
 \begin{pmatrix}
 4 & 0 & 0 & 0 \\
 0 & 4 & 0 & 0 \\
 0 & 0 & 4 & 0 \\
 0 & 0 & 0 & 4
 \end{pmatrix} = A_{W}^{T}.\]
 Hence in this case, $W^T = W$.
\end{Example}

\begin{Example}
It is generally not the case that $W^T=W$. If, for example, $W$ is the chain polynomial $W = x_{1}^3x_2 + x_{2}^2x_3 + x_{3}^2$ from Example~\ref{eg:chain}, then
 \[A_W =
 \begin{pmatrix}
 3 & 1 & 0\\
 0 & 2 & 1\\
 0 & 0 & 2\\
 \end{pmatrix},\qquad \text{so} \qquad
 A_{W}^{T}=
 \begin{pmatrix}
 3 & 0 & 0\\
 1 & 2 & 0\\
 0 & 1 & 2\\
 \end{pmatrix}.\]
In this example, we see that $W^T = x_1^3 + x_1x_2^2 + x_2x_3^2$. Notice that $W^T$ is also invertible and that its weights are $\big(\tfrac{1}{3}, \tfrac{1}{3}, \tfrac{1}{3}\big)$.
\end{Example}

Note that the exponent matrix $A_W$ from Definition~\ref{def:exp matrix} is only defined up to a reordering of rows.

\begin{Definition}\label{def:Dual Group}
The \emph{dual group} of a subgroup $G \leq G_{W}^{\rm diag}$ is the set
\[
G^T = \big\{g \in G_{W^T}^{\rm diag} \,|\, gA_{W}h^T \in \mathbb{Z} \text{ for all } h \in G\big\},
\]
where we consider $g$ and $h$ in their additive form as row vectors.
\end{Definition}

The definition of the dual group was given initially by Berglund and Henningson in~\cite{ber-hen} and independently by Krawitz in \cite{KrThesis}. Ebeling and Takahashi proved in~\cite{ET} that the two definitions agree. It is an exercise to show that the definition given here is the same as the definition of Krawitz.

\begin{Example}
 Recall from Example~\ref{eg:J_W dual} the groups $J_W$ and $\SL_{W}^{\rm diag}$ for $W = x_1^4 + x_2^4 + x_3^4 + x_4^4 = W^T$. These groups are in fact the dual groups of each other. Observe
 \[(J_W)^T = \big\{g \in G_{W^T}^{\rm diag} \,|\, gA_{W}h^T \in \mathbb{Z} \text{ for all } h \in J_W\big\}.\]
 Let $g \in G_{W^T}^{\rm diag}$ and $h \in J_W$, then $g = \big(\tfrac{a_1}{4}, \tfrac{a_2}{4}, \tfrac{a_3}{4}, \tfrac{a_4}{4}\big)$ and $h = \big(\tfrac{b}{4}, \tfrac{b}{4}, \tfrac{b}{4}, \tfrac{b}{4}\big)$ where $a_1, a_2, a_3, a_4, b \in \{0, 1, 2, 3\}$. Then
\begin{gather*}
gA_{W}h^T = \big(\tfrac{a_1}{4}, \tfrac{a_2}{4}, \tfrac{a_3}{4}, \tfrac{a_4}{4}\big)
\begin{pmatrix}
4 & 0 & 0 & 0 \\
0 & 4 & 0 & 0 \\
0 & 0 & 4 & 0 \\
0 & 0 & 0 & 4
\end{pmatrix}
\big(\tfrac{b}{4}, \tfrac{b}{4}, \tfrac{b}{4}, \tfrac{b}{4}\big)^T
 = b\big(\tfrac{a_1}{4} + \tfrac{a_2}{4} + \tfrac{a_3}{4} + \tfrac{a_4}{4}\big).
\end{gather*}

This value is an integer for all $b \in \{0, 1, 2, 3\}$ if and only if $\big(\tfrac{a_1}{4} + \tfrac{a_2}{4} + \tfrac{a_3}{4} + \tfrac{a_4}{4}\big) \in \ZZ$, implying $g \in \SL_{W^T}^{\rm diag}$. Hence $(J_W)^T = \SL_{W^T}^{\rm diag}$. In fact, it is true that $(J_W)^T=\SL_{W^T}^{\rm diag}$ for any choice of invertible polynomial $W$ (see, e.g., \cite{Arte}).
\end{Example}

Let $W$ be an invertible polynomial and $G\leq G^{\rm diag}_W$. We can now define the \emph{BHK mirror} of a~pair $(W,G)$, as $\big(W^T, G^T\big)$.

As mentioned previously, most of the work done with Landau--Ginzburg models has been with subgroups of $G_W^{\rm diag}$. Next, we consider a group with a permutation as one of its generators, which is a~non-diagonal symmetry. In order to define the dual group we need to define the \emph{non-abelian dual group}, which we do in the next section.

\subsection{The non-abelian dual group}

We begin this section with an example to illustrate the sort of symmetry groups that we will encounter in the remainder of this article.

\begin{Example}\label{eg:nondiag}
With $W = x_1^4 + x_2^4 + x_3^4 + x_4^4$ as before, consider the subgroup
\[G = \langle j_W, (123) \rangle \leq G_{W}^{\max}, \qquad \text{where} \qquad (123) = \begin{pmatrix}
 0 & 1 & 0 & 0 \\
 0 & 0 & 1 & 0 \\
 1 & 0 & 0 & 0 \\
 0 & 0 & 0 & 1
\end{pmatrix}.\]
Here $(123)$ permutes the variables $x_1$, $x_2$, and $x_3$ under the action described in Definition \ref{def:Gmax}. Even though $G$ contains non-diagonal matrices, it is actually still abelian since the generators commute (because~$j_W$ lies in the center of~$\GL_N(\mathbb{C})$; see also the remark following Definition~\ref{def:Gmax}).
\end{Example}

Although $G$ is abelian in the above example, we cannot use the previously mentioned definition for $G^T$ since $G$ is not a subset of $G_{W}^{\rm diag}$, as required by Definition~\ref{def:Dual Group}. This brings us to the definition of the non-abelian dual group. First we need one more definition.

\begin{Definition}An element of $G_W^{\max}$ is called a \emph{pure permutation} if it acts on $\CC[x_1,\dots,x_N]$ by simply permuting the variables.
\end{Definition}

Notice that because of Definition~\ref{def:Gmax} a pure permutation can only permute variables that have the same weight with respect to~$W$. We are now ready to define the non-abelian dual group~$G^\star$.

\begin{Definition}\label{def:Gstar}
Let $G\leq G_{W}^{\max}$ be a group of the form
\[
G=H\cdot K,
\]
where $K \leq G$ is the subgroup of pure even permutations and $H\leq G \cap G_{W}^{\rm diag}$.
This product should be thought of as an interior product in $\GL_N(\mathbb{C})$. Since $K$ permutes only variables with the same weight, and is a symmetry of~$W$, we see that~$K$ can also be thought of as a~subgroup of $G^{\max}_{W^T}$. We define the \emph{non-abelian dual group} of~$G$ to be
 \[
 G^\star = H^T \cdot K \leq \GL_N(\CC).
 \]
\end{Definition}

\begin{Remark}
 Ebeling and Gusein--Zade use different notation, but one can check that the definition given here is equivalent. Notice $H$ is normal in $G$.
\end{Remark}

Also notice that since $K$ and $H^T$ are both subgroups of $G^{\max}_{W^T}$, we have $G^\star \leq G^{\max}_{W^T}$.

\begin{Example} \label{eg:Gstarex}
 If we consider $G = \langle j_W, (123) \rangle \leq G_{W}^{\max}$ from Example~\ref{eg:nondiag}, then
 \[
 G^\star = J_W^T\cdot\big\langle (123) \big\rangle = \SL_{W^T}^{\rm diag}\cdot\big\langle (123) \big\rangle.
 \]
 Explicitly, the elements of $G^\star$ are of the form $\big(\tfrac{a_1}{4}, \tfrac{a_2}{4}, \tfrac{a_3}{4}, \tfrac{a_4}{4}\big)(123)^k$, where $a_1 + a_2 + a_3 + a_4 \in 4\ZZ$ and $k \in \{0, 1, 2\}.$ In this example, we can see that $G^\star$ is non-abelian. For instance, consider the products of $\big(\tfrac{1}{2}, \tfrac{1}{4}, \tfrac{1}{4}, 0\big)(123)$ and $\big(\tfrac{1}{2}, \tfrac{1}{4}, \tfrac{1}{4}, 0\big)(132) \in G^\star$ in both ways. Observe
 \[
 \big(\tfrac{1}{2}, \tfrac{1}{4}, \tfrac{1}{4}, 0\big)(123) \cdot \big(\tfrac{1}{2}, \tfrac{1}{4}, \tfrac{1}{4}, 0\big)(132) = \big(\tfrac{3}{4}, \tfrac{1}{2}, \tfrac{3}{4}, 0\big),
 \]
 whereas
 \[
 \big(\tfrac{1}{2}, \tfrac{1}{4}, \tfrac{1}{4}, 0\big)(132) \cdot \big(\tfrac{1}{2}, \tfrac{1}{4}, \tfrac{1}{4}, 0\big)(123) = \big(\tfrac{3}{4}, \tfrac{3}{4}, \tfrac{1}{2}, 0\big).
 \]
\end{Example}

Now we have defined a rule relating two LG models $(W, G)$ and $\big(W^T,G^\star\big)$.

\section{Construction of the state space}\label{s:statespace}

We are now ready to construct the A- and B-model state spaces. We will begin with the A-model.

\subsection{A-model state space}

The A-model state space is defined using relative (orbifold) cohomology of the Milnor fiber (see~\cite{FJR:GLSM}). For simplicity we will give an alternate equivalent definition in terms of Milnor rings. First we need to define a few of the ingredients.

\begin{Definition}The \emph{Milnor ring} of a polynomial $W$ is defined to be
\[
\cQ_W = \frac{\mathbb{C}[x_1, \dots, x_N]}{\big( \frac{\partial W}{\partial x_1}, \dots, \frac{\partial W}{\partial x_N} \big)}.
\]
\end{Definition}

\begin{Definition}
 Let $W$ be a nondegenerate, quasihomogenenous polynomial with unique weights $(q_1, \dots, q_N)$, and let $G$ be a subgroup of $G_{W}^{\max}$. Then $G$ is \emph{A-admissible} if $G$ contains $j_W = (q_1, \dots, q_N)$.
\end{Definition}

\begin{Definition}
Given an element $g \in G_W^{\max}$, we let $\Fix(g)$ denote the subspace of $\mathbb{C}^N$ which is fixed by $g$, i.e.
\[
\Fix(g) = \{(a_1, \dots, a_N) \,|\, g\cdot(a_1, \dots, a_N) = (a_1, \dots, a_N)\}.
\]
\end{Definition}
To find $\Fix(g),$ we look for eigenvectors of $g$ with an eigenvalue of $1$, and these vectors will span $\Fix(g)$. We also write
\[
W_g = W|_{\Fix(g)}
\]
to denote the polynomial $W$ restricted to $\Fix(g)$. Let $N_g=\dim(\Fix(g)$).

\begin{Definition}\label{def:Astate}
Let $W$ be an invertible polynomial and $G$ be an A-admissible subgroup of~$G_{W}^{\max}$. The \emph{state space} for the A-model is defined as
\[
\mathcal{A}_{W, G} = \Bigg(\underset{g \in G}{\bigoplus}\cQ_{W_g} \cdot \omega_g\Bigg)^G,
\]
where $\omega_g$ is a volume form on the fixed locus of $g$.
\end{Definition}

As mentioned earlier, the original definition is given in terms of relative cohomologies of a~Milnor fiber of $W$. This is simply an equivalent definition for LG A-models that is easier to work with.

When $G$ is abelian, we can rewrite the state space definition as
\[
\mathcal{A}_{W, G} = \underset{g \in G}{\bigoplus}(\cQ_{W_g} \cdot \omega_g)^G,
\]
as the action of $G$ preserves each summand. This is the definition of the A-model more commonly seen for the A-model state space when using diagonal symmetries.
However, if~$G$ is non-abelian, then for $\gamma\in G$,
\begin{gather}\label{eq:Haction}
\gamma \cdot \big(\cQ_{W_g} \cdot \omega_g\big) \subseteq \cQ_{W_{\gamma g \gamma^{-1}}} \cdot \omega_{\gamma g \gamma^{-1}}.
\end{gather}
Here the action of $\gamma$ is induced by $\big(\gamma^T\big)^*$, where if $A$ is a matrix (thought of as a linear transformation), then $A^*$ is the action on the dual vector space. In other words, if $\sigma$ is a permutation, written in cycle notation, then $\sigma\cdot x_i=x_{\sigma^{-1}(i)}$.

We will use the notation $\lfloor P, g \rceil$ to denote an element of $\cQ_{W_g}\cdot\omega_g$, often suppressing the volume form where convenient. The volume form can be easily determined by $g$. We can form a basis of $\mathcal{A}_{W, G}$ using sums of the form
\[
\sum_{g_i \in [g]} \lfloor P, g_i \rceil,
\]
where $g_i$ are the group elements in the same conjugacy class $[g]$ of $G$, and $P \in \cQ_{W_{g_i}}\cdot \omega_g$.

The A-model can also be given a bigrading which in many cases is similar to the Hodge grading for Calabi--Yau manifolds. Since mirror symmetry for Calabi--Yau manifolds rotates the Hodge diamond, we expect a similar phenomenon for LG models. In fact, we will see that the B-model also has a bigrading and that the bigrading is preserved by mirror symmetry under certain conditions.

\begin{Definition}[Mukai \cite{Mukai}]\label{def:agedef}
 Let $G$ be a finite subgroup of the symmetry group of some nondegenerate quasihomogeneous polynomial in $\CC[x_1, \dots, x_N]$. We define the \emph{age} of $g \in G$ as
 \[
 \age g = \frac{1}{2 \pi {\rm i}} \sum_{j = 1}^N \log (\lambda_j),
 \]
 where $\lambda_1, \dots, \lambda_N$ are the eigenvalues of $g$ and the branch of the logarithmic function for $z \in \CC^*$ is chosen to satisfy $0 \leq \log (z) < 2 \pi {\rm i}$. Notice that all eigenvalues of $g\in G$ satisfy $|\lambda|=1$, since~$G$ has finite order.
\end{Definition}

\begin{Example}\label{eg:ageex}
 For $g \in G^{\rm diag}_W$, we can write $g=(a_1,\dots,a_N)$ additively.
 Then the age of $g$ is just $\sum\limits_{j=1}^{N} a_j$, where $a_j$ is chosen so that $0\leq a_j<1$.
\end{Example}

\begin{Definition}[Krawitz \cite{KrThesis}]\label{def:Agrading}
 The A-model \emph{bigrading} of $\lfloor P, g \rceil$ is defined to be the ordered pair
 \[
 (\deg P + \age g - \age j_W, N_g - \deg P + \age g - \age j_W),
 \]
 where $\deg P$ is the weighted degree of $P$. In this notation, note that the volume form $\omega_g$ contributes to $\deg P$. Recall $N_g$ is the dimension of~$\Fix(g)$.
\end{Definition}

\subsection{An extended example}

\begin{Example}\label{eg:AStateConstruction}
Let $W = x_{1}^4 + x_{2}^4 + x_{3}^4 + x_{4}^4$ and $G = \langle j_W, (123) \rangle$. We will determine a basis for~$\mathcal{A}_{W, G}$. Since in this case, $G$ is an abelian group, the conjugacy class for each $g \in G$ contains only $g$. Hence we can choose a basis of~$\mathcal{A}_{W, G}$ consisting of elements of the form $\lfloor P, g \rceil$ (i.e., single terms, instead of sums, although $P$ may have more than one summand). The elements of~$G$ can be expressed as $j_{W}^a(123)^b$ with $a \in \{0, 1, 2, 3\}$ and $b \in \{0, 1, 2\}.$ For each $g \in G$, we will need to find the basis elements of $(\cQ_{W_g} \cdot \omega_g)^G$. The choices of $g$ can be broken down into three different cases.

\textbf{Case 1: $g = (0, 0, 0, 0)$}.
When $g = (0, 0, 0, 0)$, then $W_g = W$, and the Milnor ring of $W_g$ is
\[
\cQ_{W_g} = \cQ_W = \frac{\mathbb{C}[x_1, x_2, x_3, x_4]}{\big( 4x_{1}^3, 4x_{2}^3, 4x_{3}^3, 4x_{4}^3 \big)}.
\]
The elements of $\cQ_W$ are sums of elements in the set $\big\{x_{1}^ax_{2}^bx_{3}^cx_{4}^d \,|\, 0 \leq a, b, c, d \leq 2\big\}.$ The volume form $\omega_g$ in this case is ${\rm d}x_1 \wedge {\rm d}x_2 \wedge {\rm d}x_3 \wedge {\rm d}x_4$. To find the elements of $(\cQ_W \cdot ({\rm d}x_1 \wedge {\rm d}x_2 \wedge {\rm d}x_3 \wedge {\rm d}x_4))^G$ we look for $p(x) \in \cQ_W$ such that $p(x) \cdot ({\rm d}x_1 \wedge {\rm d}x_2 \wedge {\rm d}x_3 \wedge {\rm d}x_4)$ is invariant under $j_W$ and $(123)$, the generators of~$G$. The volume form is invariant under~$j_W$ and under~$(123)$.
Thus, in this case we only need to be concerned with the actual polynomial~$p(x)$.

In order to be invariant under (123), the polynomial must be symmetric with respect to~$x_1$,~$x_2$, and $x_3$ and polynomials invariant under $j_W$ must have exponents in each term sum to a~multiple of~4; for example, the polynomial $x_1x_2x_3x_4 \in \cQ_W$ is invariant under both~$j_W$ and~(123).

One can check that the invariant elements of $\cQ_W\cdot\omega_g$ are spanned by the following polynomials:
\begin{gather*}
1,\\
x_1x_2x_3x_4,\\
x_1^2x_2^2x_3^2x_4^2,\\
x_1^2x_2^2 + x_{1}^2x_{3}^2 + x_{2}^2x_{3}^2,\\
x_1x_2x_3^2 + x_1^2x_2x_3 + x_1x_2^2x_3,\\
x_1x_2x_4^2 + x_1x_3x_4^2 + x_2x_3x_4^2,\\
x_1x_2^2x_4 + x_2x_3^2x_4 + x_3x_1^2x_4,\\
x_1^2x_2x_4 + x_2^2x_3x_4 + x_3^2x_1x_4,\\
x_1^2x_4^2 + x_2^2x_4^2 + x_3^2x_4^2.
\end{gather*}
The $9$-dimensional vector space generated by these elements is called the \emph{untwisted broad sector} of $\mathcal{A}_{W, G}$.

We now turn our attention to the bigrading. Since $\age j_W= 1$ as mentioned in Example~\ref{eg:ageex}, the bidegree for each element reduces to
\[
(\deg P + \age g - 1, N_g - \deg P + \age g - 1),
\]

Furthermore, when $g$ is the identity, we get $\age g = 0$ and $N_g = 4$, so the bidegree simplifies to
\[
(\deg P - 1, 3 - \deg P).
\]
Recall that there were $9$ polynomials in our basis for this choice of $g$.
We will now list the bidegree for all of the basis elements in this sector:
\begin{alignat*}{3}
&\text{{basis element}} \qquad &&\text{{bidegree}}& \\
&\lfloor 1, (0, 0, 0, 0) \rceil \qquad &&(0, 2)&\\
&\lfloor x_1x_2x_3x_4, (0, 0, 0, 0) \rceil \qquad &&(1, 1)&\\
&\big\lfloor x_1^2x_2^2x_3^2x_4^2, (0, 0, 0, 0) \big\rceil \qquad &&(2, 0)&\\
&\big\lfloor x_1^2x_2^2 + x_{1}^2x_{3}^2 + x_{2}^2x_{3}^2, (0, 0, 0, 0) \big\rceil \qquad &&(1, 1)&\\
&\big\lfloor x_1x_2x_3^2 + x_1^2x_2x_3 + x_1x_2^2x_3, (0, 0, 0, 0) \big\rceil \qquad &&(1, 1)&\\
&\big\lfloor x_1x_2x_4^2 + x_1x_3x_4^2 + x_2x_3x_4^2, (0, 0, 0, 0) \big\rceil \qquad &&(1, 1)&\\
&\big\lfloor x_1x_2^2x_4 + x_2x_3^2x_4 + x_3x_1^2x_4, (0, 0, 0, 0) \big\rceil \qquad &&(1, 1)&\\
&\big\lfloor x_1^2x_2x_4 + x_2^2x_3x_4 + x_3^2x_1x_4, (0, 0, 0, 0) \big\rceil \qquad &&(1, 1)&\\
&\big\lfloor x_1^2x_4^2 + x_2^2x_4^2 + x_3^2x_4^2, (0, 0, 0, 0) \big\rceil \qquad &&(1, 1)&
\end{alignat*}

This completes the construction of the untwisted sector. Next we consider the pure permutations.

\textbf{Case 2: $g = (123)$ or $g = (132)$.} Let
\[
g=(123) = \begin{pmatrix}
0&1&0&0\\
0&0&1&0\\
1&0&0&0\\
0&0&0&1
\end{pmatrix}.
\]
This representation of $(123)$ as a matrix matches the convention that $\sigma\cdot x_i=x_{\sigma^{-1}(i)}$. To find $\Fix(123)$, we look for eigenvectors of $(123)$ with an eigenvalue of~1. Diagonalizing $(123)$ gives $(123)=QDQ^{-1}$, where
\[
D = \begin{pmatrix}
1&0&0&0\\
0&1&0&0\\
0&0&{\rm e}^\frac{4 \pi {\rm i}}{3}&0\\
0&0&0&{\rm e}^\frac{2 \pi {\rm i}}{3}
\end{pmatrix} \qquad \text{and}\qquad
Q = \begin{pmatrix}
0&1&{\rm e}^\frac{4 \pi {\rm i}}{3}&{\rm e}^\frac{2 \pi {\rm i}}{3}\\
0&1&{\rm e}^\frac{2 \pi {\rm i}}{3}&{\rm e}^\frac{4 \pi {\rm i}}{3}\\
0&1&1&1\\
1&0&0&0
\end{pmatrix}.
\]

Thus the eigenvectors with eigenvalue 1 are $(1, 1, 1, 0)$ and $(0, 0, 0, 1)$, and the span of these two vectors is $\Fix(123)$. If we call the coordinates of these two vectors $y_1$ and $y_4$, then we have $W_g = c_1y_1^4+y_4^4$ for some constant $c_1$ whose value does not matter for our purposes.
The volume form here is ${\rm d}y_1 \wedge {\rm d}y_4=({\rm d}x_1+{\rm d}x_2+{\rm d}x_3)\wedge {\rm d}x_4$. This is invariant under $(123)$, which acts trivially when considering only $y_1$ and $y_4$. However, this volume form is not invariant under $j_W$ since
\[
j_W\cdot ({\rm d}y_1 \wedge {\rm d}y_4) = -({\rm d}y_1 \wedge {\rm d}y_4) \neq {\rm d}y_1 \wedge {\rm d}y_4.
\]
To balance this, in order for an element of $\cQ_{W_g} \cdot \omega_g$ to be invariant under $j_W$, the polyno\-mial~$p(x)$ must have degree equal to $2$ (mod~$4$). This gives us three anti-invariant polynomials:
\begin{gather*}
y_1^2 = (x_1 + x_2 + x_3)^2,\\
y_1y_4 = (x_1+x_2+x_3)x_4,\\
y_4^2 = x_4^2.
\end{gather*}
Each one of these, together with the volume form, is another element in the basis of $\mathcal{A}_{W, G}$.

The case of $g=(132)$ is almost identical, so we exclude the work here. The subspaces produced from non-identity group elements $g$ with $\Fix g \neq \{0\}$ are known as \emph{twisted broad sectors}.

We now consider the bigrading for these elements.
The values of $\age (123)$ and $\age(132)$ are both $1$.
There were $3$ elements in the $G$-invariant subspace of the Milnor ring, and their degrees are all equal to~1.

Thus the bidegree of each of these elements is
\[
(\deg P + \age g - 1, N_g - \deg P + \age g - 1)
= (1, 1).
\]

\textbf{Case 3: Other values of $g$.} The eigenvalues of $j_W$ are all ${\rm e}^\frac{2 \pi {\rm i}}{4}$, so $g = j_W$ has trivial fixed locus. Thus $W|_{\Fix(j_W)} = 0$. This implies that for $g = j_W$, we get $\cQ_{W_g} \cdot \omega_g \cong \CC$. This sector only produces a single basis element of $\mathcal{A}_{W, G}$, being $\lfloor 1, j_W \rceil$. Sectors with $\Fix(g)= \{0\}$, as is the case here, are called \emph{narrow sectors}. The action of $G$ on these narrow sectors is trivial, so each contributes to the basis for $\cA_{W,G}$. Similarly, the group elements $(j_W)^2$, $(j_W)^3$, $j_W(123)$, $(j_W)^2(123)$, $(j_W)^3(123)$, $(j_W)(132)$, $(j_W)^2(132)$, and $(j_W)^3(132)$ are also narrow sectors. In total, there are $9$ narrow sectors in $\mathcal{A}_{W, G}$.

To compute the bidegrees of these elements, we first notice that
$N_g = 0$, and $\deg P=0$ in the formula for bidegree. Thus formula for bidegree thus reduces to
\[
(\age g - \age j_W, \age g - \age j_W) =
(\age g - 1, \age g - 1).
\]
Hence in this case, the only thing we need to actually compute is $\age g$. When $g$ is a multiple of $(j_W)$, we have
\[
\age j_W = 1, \qquad \age j_W^2 = 2, \qquad \text{and}\qquad \age j_W^3 = 3.
\]

The rest of the elements are non-diagonal, so we must find the eigenvalues as in the previous case. The resulting age is the same for all of them, which is~$2$. Thus the bidegree for the rest of the narrow sectors is $(1, 1)$. To conclude this example, we have found that there are 9 narrow sectors, the untwisted broad sector has dimension~9, and the two twisted broad sectors from $(123)$ and $(132)$ each contribute dimension~3 to the state space. Hence $\mathcal{A}_{W,G}$ has dimension~24.
In following table, Table~\ref{tab:Amodelbasis}, we can see all basis elements and their bidegree.

\begin{table}[t!]\centering
 \begin{longtable}{ l c} \toprule
 \multicolumn{2}{c}{A-model basis elements} \\
{A-model basis} &{Bidegree} \\
 \midrule
 \endfirsthead
 \toprule
{A-model basis} &{Bidegree} \\
 \midrule
 \endhead
 $\lfloor 1, (0, 0, 0, 0) \rceil$ &$(0, 2)$\\
 $\lfloor x_1x_2x_3x_4, (0, 0, 0, 0) \rceil$ & $(1, 1)$\tsep{1pt}\\
 $\big\lfloor x_1^2x_2^2x_3^2x_4^2, (0, 0, 0, 0) \big\rceil$ & $(2, 0)$\tsep{1pt}\\
 $\big\lfloor x_1^2x_2^2 + x_{1}^2x_{3}^2 + x_{2}^2x_{3}^2, (0, 0, 0, 0) \big\rceil$ &$(1, 1)$\tsep{1pt}\\
 $\big\lfloor x_1x_2x_3^2 + x_1^2x_2x_3 + x_1x_2^2x_3, (0, 0, 0, 0) \big\rceil$ &$(1, 1)$\tsep{1pt}\\
 $\big\lfloor x_1x_2x_4^2 + x_1x_3x_4^2 + x_2x_3x_4^2, (0, 0, 0, 0) \big\rceil$ &$(1, 1)$\tsep{1pt}\\
 $\big\lfloor x_1x_2^2x_4 + x_2x_3^2x_4 + x_1^2x_3x_4, (0, 0, 0, 0) \big\rceil$ &$(1, 1)$\tsep{1pt}\\
 $\big\lfloor x_1^2x_2x_4 + x_2^2x_3x_4 + x_1x_3^2x_4, (0, 0, 0, 0) \big\rceil$ &$(1, 1)$\tsep{1pt}\\
 $\big\lfloor x_1^2x_4^2 + x_2^2x_4^2 + x_3^2x_4^2, (0, 0, 0, 0) \big\rceil$ &$(1, 1)$\tsep{1pt}\\
 $\big\lfloor (x_1 + x_2 + x_3)^2, (123) \big\rceil$ &$(1, 1)$\tsep{1pt}\\
 $\lfloor (x_1 + x_2 + x_3)x_4, (123) \rceil$ &$(1, 1)$\tsep{1pt}\\
 $\big\lfloor x_4^2, (123) \big\rceil$ &$(1, 1)$\tsep{1pt}\\
 $\big\lfloor (x_1 + x_2 + x_3)^2, (132) \big\rceil$ &$(1, 1)$\tsep{1pt}\\
 $\lfloor (x_1 + x_2 + x_3)x_4, (132) \rceil$ &$(1, 1)$\tsep{1pt}\\
 $\big\lfloor x_4^2, (132) \big\rceil$ &$(1, 1)$\tsep{1pt}\\
 $\lfloor 1, j_W \rceil$ &$(0, 0)$\tsep{1pt}\\
 $\lfloor 1, (j_W)^2 \rceil$ &$(1, 1)$\tsep{1pt}\\
 $\big\lfloor 1, (j_W)^3 \big\rceil$ &$(2, 2)$\tsep{1pt}\\
 $\lfloor 1, j_W(123) \rceil$ &$(1, 1)$\tsep{1pt}\\
 $\big\lfloor 1, (j_W)^2(123) \big\rceil$ &$(1, 1)$\tsep{1pt}\\
 $\big\lfloor 1, (j_W)^3(123) \big\rceil$ &$(1, 1)$\tsep{1pt}\\
 $\lfloor 1, j_W(132) \rceil$ &$(1, 1)$\tsep{1pt}\\
 $\big\lfloor 1, (j_W)^2(132) \big\rceil$ &$(1, 1)$\tsep{1pt}\\
 $\big\lfloor 1, (j_W)^3(132) \big\rceil$ &$(1, 1)$\tsep{1pt}\\
 \bottomrule
\caption{A-model basis elements.}\label{tab:Amodelbasis}
\end{longtable}
\end{table}

If we arrange these as a Hodge diamond, we have
\[
\begin{array}{@{}ccc}
 & 1 &\\
 1 & 20 & 1 \\
 & 1 &
\end{array}
\]
The reader may notice that this is the Hodge diamond of a K3 surface. Indeed if we consider the K3 surface $X_W$ defined by the polynomial $x_1^4 + x_2^4+x_3^4 + x_4^4$ in $\PP^3$, then one can quotient by the group action of~$(123)$ on~$X_W$. The minimal resolution of $X_W/\langle (123)\rangle$ is a K3 surface, corresponding to the LG model we have just considered (the Fermat quartic is a K3 surface and~$(123)$ acts symplectically on it; see, e.g., \cite{CLPS}).
\end{Example}

We end this section with a lemma that will be useful for the proof of Theorem~\ref{thm:mirror}.

\begin{Lemma}\label{lem:Abidegree}
Suppose $G$ is of the form $G=H\cdot K$, where $K$ is a group of pure permutations and $H \leq G^{\rm diag}_W$ is an A-admissible group. If $\gamma \in G$ and $\lfloor P, g \rceil \in \cQ_{W_g}\cdot \omega_g$, then $\gamma \cdot \lfloor P, g \rceil$ has the same bidegree as $\lfloor P, g \rceil$.
\end{Lemma}

\begin{proof}
First, recall the A-model bigrading from Definition~\ref{def:Agrading}:
\[
(\deg P + \age g - \age j_W, N_g - \deg P + \age g - \age j_W).
\]
We also recall \eqref{eq:Haction} from right after Definition~\ref{def:Astate}, which says that
\[
\gamma \cdot \lfloor P, g \rceil = \big\lfloor \gamma \cdot P, \gamma g \gamma^{-1} \big\rceil.
\]

Note that $\age j_W$ will clearly be unaffected by the action of $\gamma$ on $\lfloor P, g \rceil$. We aim to show that $\age \big(\gamma g \gamma^{-1}\big) = \age g$, that $N_{\gamma g \gamma^{-1}} = N_g$, and that $\deg (\gamma \cdot P) = \deg P$. Recall from Definition~\ref{def:agedef} that the age of $g$ is dependent only on the eigenvalues of $g$. Since $g$ and $\gamma g \gamma^{-1}$ are similar matrices, they must have the same eigenvalues, so $\age \big(\gamma g \gamma^{-1}\big) = \age g$. This also gives us $N_{(\gamma g \gamma^{-1})} = N_g$, since $N_g$ is the number of eigenvalues of $g$ which are equal to one. To show that $\deg (\gamma \cdot P) = \deg P$, we will consider two cases for~$\gamma$: either $\gamma$ is a pure permutation or $\gamma$ is a~diagonal symmetry. A third case would be when $\gamma$ is a product of a~pure permutation and a~diagonal symmetry, but this follows from the previous two cases.

\textbf{Case 1:} Suppose $\gamma=\sigma$ is a pure permutation. Recall from Definition~\ref{def:Gmax} that the elements of $G_W^{\max}$ only permute variables with the same weight. The degree of the volume form is also unaffected for the same reason. Then $\sigma \cdot P$ simply renames the indexes of the variables which will not change the degree at all. Thus $\deg (\sigma \cdot P) = \deg P$.

\textbf{Case 2:} Suppose $\gamma=h$ is a diagonal symmetry, and put $h=(a_1, a_2, \dots, a_N)$, written additively.
Then $h \cdot P = cP$ for some $c \in \mathbb{C}^*$,
so $\deg (h \cdot P) = \deg P$.

In any case, $\gamma \cdot \lfloor P, g \rceil$ has the same bidegree as $\lfloor P, g \rceil$ in $\mathcal{A}_{W, G}$.
\end{proof}

\subsection{The B-model state space}
Having constructed the A-model as a bigraded vector space, we can begin our construction of the B-model. We expect the B-model for $\big(W^T, G^\star\big)$ to be isomorphic to the A-model for $(W, G)$.

\begin{Definition}
 Let $W$ be a nondegenerate quasihomogenenous polynomial with unique \linebreak weights $(q_1, \dots, q_N)$, and let $G$ be a subgroup of $G_{W}^{\max}$. Then $G$ is \emph{B-admissible} if $G\subset\SL_{W}$.
\end{Definition}

\begin{Definition}
Let $W$ be an invertible polynomial and $G$ a B-admissible group. The \emph{state space} for the B-model is defined as
\[
\mathcal{B}_{W, G} = \Bigg(\underset{g \in G}{\bigoplus}\cQ_{W_g} \cdot \omega_g\Bigg)^G,
\]
where $\omega_g$ is a volume form on the fixed locus of $g$.
\end{Definition}

This is exactly analogous to Definition~\ref{def:Astate}, except that the associated group $G$ has different requirements than the group used for the A-model. For the B-model, this is the original definition, whereas Definition~\ref{def:Astate} for the A-model was only an equivalent definition. While the state spaces have similar definitions, the grading and product structures are very different~-- for example, the multiplicative identity on the B-side lies in the untwisted sector, whereas the multiplicative identity on the A-side lies in the sector indexed by $j_W$. We won't discuss the pro\-duct structure any further here, but the interested reader can find the definition of the B-model multiplication in~\cite{BTW} or~\cite{HLL}, and the definition of the A-model product in~\cite{FJR} or~\cite{FJRW}, when~$G$ and~$G^*$ are groups of diagonal symmetries. When~$G$ and~$G^*$ are not diagonal, the B-model product has not yet been defined, as far as we know, whereas the A-model structure comes from invariants of the associated Gauged Linear Sigma Model (see, e.g.,~\cite{FJR:GLSM}, though we are unaware of the invariants having been computed for any examples). For this article, the most interesting aspect of this state space is the relationship with the A-model under mirror symmetry.

If we use $(W, G)$ to construct the A-model with $G \leq G_W^{\rm diag}$, then Krawitz~\cite{KrThesis} showed that{\samepage
\[
\cA_{W,G}\cong \cB_{W^T,G^T}
\]
as bigraded vector spaces.}

For groups of non-diagonal matrices, in order for mirror symmetry to hold we replace $G^T$ by the non-abelian dual group $G^\star$, defined in Definition~\ref{def:Gstar}.

Just like with the A-model, there is also a bigrading on the B-model state space.

\begin{Definition}[Krawitz \cite{KrThesis}]\label{def:Bgrading}
 The \emph{B-model bigrading} of $\lfloor P, g \rceil$ is defined to be the ordered pair
 \[
 \big(\deg P + \age g - \age j_W, \deg P + \age g^{-1} - \age j_W\big).
 \]
As with the A-model bigrading from Definition~\ref{def:Agrading}, the volume form $\omega_g$ contributes to $\deg P$.
\end{Definition}

\subsection{An extended example}

\begin{Example}\label{eg:BStateConstruction}
 Let $W = x_1^4 + x_2^4 + x_3^4 + x_4^4$, with $G = \langle j_W, (123) \rangle$. From
 Example~\ref{eg:Gstarex}, we saw that $W^T=W$ and that $G^\star = \SL_{W^T}^{\rm diag}\cdot\langle (123) \rangle $. The elements of $\SL_{W^T}^{\rm diag}$ are of the form $\big(\tfrac{a_1}{4},\tfrac{a_2}{4},\tfrac{a_3}{4},\tfrac{a_4}{4}\big)\cdot(123)^k$; again, the notation $\big(\tfrac{a_1}{4}, \tfrac{a_2}{4}, \tfrac{a_3}{4}, \tfrac{a_4}{4}\big)$ refers to a 4x4 diagonal matrix with diagonal entries on the complex unit circle. The entries also satisfy $4|(a_1+a_2+a_3+a_4)$~-- the requirement to be in $\SL_4(\mathbb{C})$. Furthermore, $(123)\in\SL_W$, since it is an even permutation. The group $G^\star$ is generated by $(123)$, $j_W$, $K$, and $L$, where $j_W=\big(\frac{1}{4},\frac{1}{4},\frac{1}{4},\frac{1}{4}\big)$, $K=\big(\frac{1}{2},\frac{1}{4},\frac{1}{4},0\big)$, and $L=\big(\frac{1}{4},\frac{1}{2},\frac{1}{4},0\big)$.

 The A-model had $24$ basis elements, with $20$ of them having a bidegree of $(1,1)$ and $1$ of each of the following: $(0, 0)$, $(2, 0)$, $(0, 2)$, and $(2, 2)$. For the A- and B-models to be isomorphic as bigraded vector spaces, we should see the same breakdown of elements for the B-model as well.

 As we begin to construct $\mathcal{B}_{W^T, G^\star}$, we need to pay attention to centralizers and conjugacy classes. As on the A side (see~\eqref{eq:Haction}), we also have the property
 \begin{gather}\label{eq:BsideHaction}
 \gamma \cdot (\mathcal{B}_{W^T_g} \cdot \omega_g) \subseteq \mathcal{B}_{W^T_{\gamma g \gamma^{-1}}} \cdot \omega_{\gamma g \gamma^{-1}}.
 \end{gather}
 On the A-side, $j_W$ commuted with (123), so the centralizer of every element was $G$ and the conjugacy class of every element was itself. That is not the case for $G^\star$.

 \textbf{Case 1: $g = (0, 0, 0, 0)$.} Given that $W^T = W$, the Milnor ring here will be exactly the same as in case 1 of Example~\ref{eg:AStateConstruction}. However, the list of polynomials invariant under $G^\star$ will not be the same as that for $G$ since $G^\star$ has different generators. Since $(123)$, $j_W \in G^\star$, this list of polynomials will be a subset of the~$9$ from earlier, but we also need to check if those~$9$ polynomials are invariant under $K=\big(\frac{1}{2},\frac{1}{4},\frac{1}{4},0\big)$, and $L=\big(\frac{1}{4},\frac{1}{2},\frac{1}{4},0\big)$ as well. The only polynomials that will work are those where each monomial has the same exponent for $x_1$, $x_2,$ and~$x_3$.
 This $G^\star$-invariant subspace has dimension $3$ spanned by $1$, $x_1x_2x_3x_4$, and $x_1^2x_2^2x_3^2x_4^2$ (again suppressing the volume form).

 We now consider the bigrading for these. As with the A-model, we know that $\age j_W = 1$. Hence the bidegree for all of the elements in the B-model can be reduced to
 \[
 \big(\deg P + \age g - 1, \deg P + \age g^{-1} - 1\big).
 \]
 In this case, we also have $\age g = 0$, and thus $\age g^{-1} = 0$. Hence the bidegree for the elements in this sector is
 \[
 (\deg P - 1, \deg P - 1).
 \]
Thus we have three basis elements with bidegree $(0,0)$, $(1,1)$ and $(2,2)$, respectively.

\textbf{Case 2: $g$ is conjugate to $(123)$ or $(132)$.} The conjugacy class of $(123)$ contains all elements of the form $L^iK^j(123)$, for $0\leq i,j\leq 3$. Thus there are 16 elements of $G$ in the conjugacy class.

 Since $j_W$ lies in the center of $G$, we know that the polynomials in either of these sectors must be generated by those found in the A-model. Recall there were the same three polynomials for both $(123)$ and $(132)$. However, recalling~\eqref{eq:BsideHaction}, we will act on each of these three monomials in~$\cB_{(123)}$ to construct an invariant element for each of these two conjugacy classes. For example, if we consider the element $\big\lfloor (x_1+x_2+x_3)^2, (123) \big\rceil$, the following table lists each element in the conjugacy class of $(123)$ together with a monomial in $\mathcal{B}_{W^T_{\gamma (123) \gamma^{-1}}} \cdot \omega_{\gamma (123) \gamma^{-1}}$. If we add up each of these monomials, we get an invariant element of $\cB_{W^T,G^*}$. The middle column lists an eigenvector of the corresponding matrix with eigenvalue 1, so the reader can easily find $\Fix g$.
 \begin{table}[h!]\centering \footnotesize
 \begin{tabular}{|c |c| c|}
 \hline
 $g\in [(123)]$ & eigenvector & monomial\tsep{1pt}\bsep{1pt}\\
 \hline\hline
 $(123)$ & $(1, 1, 1, 0)$ & $\big\lfloor (x_1 + x_2 + x_3)^2({\rm d}x_1 + {\rm d}x_2 + {\rm d}x_3){\rm d}x_4, (123) \big\rceil$\tsep{2pt}\bsep{2pt}\\
 \hline
 $L(123)L^{-1}$ & $(1, i, 1, 0)$ & $-i\big\lfloor (x_1 + ix_2 + x_3)^2({\rm d}x_1 + i{\rm d}x_2 + {\rm d}x_3){\rm d}x_4, L(123)L^{-1} \big\rceil$\tsep{2pt}\bsep{2pt}\\
 \hline
 $L^2(123)L^2$ & $(1, -1, 1, 0)$ & $-\big\lfloor (x_1 - x_2 + x_3)^2({\rm d}x_1 - {\rm d}x_2 + {\rm d}x_3){\rm d}x_4, L^2(123)L^2 \big\rceil$\tsep{2pt}\bsep{2pt}\\
 \hline
 $L^3(123)L$ & $(1, -i, 1, 0)$ & $i\big\lfloor (x_1 - ix_2 + x_3)^2({\rm d}x_1 - i{\rm d}x_2 + {\rm d}x_3){\rm d}x_4, L^3(123)L \big\rceil$\tsep{2pt}\bsep{2pt}\\
 \hline
 $K(123)K^3$ & $(1, -i, -i, 0)$ & $-\big\lfloor (x_1 -ix_2 -ix_3)^2({\rm d}x_1 -i{\rm d}x_2 -i{\rm d}x_3){\rm d}x_4, K(123)K^3 \big\rceil$\tsep{2pt}\bsep{2pt}\\
 \hline
 $K^2(123)K^2$ & $(1, -1, -1, 0)$ & $\big\lfloor (x_1 - x_2 - x_3)^2({\rm d}x_1 - {\rm d}x_2 - {\rm d}x_3){\rm d}x_4, K^2(123)K^2 \big\rceil$\tsep{2pt}\bsep{2pt}\\
 \hline
 $K^3(123)K$ & $(1, i, i, 0)$ & $-\big\lfloor (x_1 + ix_2 + ix_3)^2({\rm d}x_1 + i{\rm d}x_2 + i{\rm d}x_3){\rm d}x_4, (KL)(123)(KL)^{-1} \big\rceil$\tsep{2pt}\bsep{2pt}\\
 \hline
 $(KL)(123)(KL)^{-1}$ & $(1, 1, -i, 0)$ & $i\big\lfloor (x_1 + x_2 - ix_3)^2({\rm d}x_1 + {\rm d}x_2 - i{\rm d}x_3){\rm d}x_4, (KL)(123)(KL)^{-1} \big\rceil$\tsep{2pt}\bsep{2pt}\\
 \hline
 $(KL)^2(123)(KL)^2$ & $(1, 1, -1, 0)$ & $-\big\lfloor (x_1 + x_2 - x_3)^2({\rm d}x_1 + {\rm d}x_2 - {\rm d}x_3){\rm d}x_4, (KL)^2(123)(KL)^2 \big\rceil$\tsep{2pt}\bsep{2pt}\\
 \hline
 $(KL)^3(123)(KL)$ & $(1, 1, i, 0)$ & $-i\big\lfloor (x_1 + x_2 + ix_3)^2({\rm d}x_1 + {\rm d}x_2 + i{\rm d}x_3){\rm d}x_4, (KL)^3(123)(KL) \big\rceil$\tsep{2pt}\bsep{2pt}\\
 \hline
 $\big(K^2L\big)(123)\big(K^2L\big)^{-1}$ & $(1, -i, -1, 0)$ & $-i\big\lfloor (x_1 - ix_2 - x_3)^2({\rm d}x_1 - i{\rm d}x_2 - {\rm d}x_3){\rm d}x_4, \big(K^2L\big)(123)\big(K^2L\big)^{-1} \big\rceil$\tsep{2pt}\bsep{2pt}\\
 \hline
 $\big(K^3L\big)(123)\big(K^3L\big)^{-1}$ & $(1, -1, i, 0)$ & $i\big\lfloor (x_1 - x_2 + ix_3)^2({\rm d}x_1 - {\rm d}x_2 +i{\rm d}x_3){\rm d}x_4, \big(K^3L\big)(123)\big(K^3L\big)^{-1} \big\rceil$\tsep{2pt}\bsep{2pt}\\
 \hline
 $\big(KL^2\big)(123)\big(KL^2\big)^{-1}$ & $(1, i, -i, 0)$ & $\big\lfloor (x_1 + ix_2 -ix_3)^2({\rm d}x_1 + i{\rm d}x_2 -i{\rm d}x_3){\rm d}x_4, \big(KL^2\big)(123)\big(KL^2\big)^{-1} \big\rceil$\tsep{2pt}\bsep{2pt}\\
 \hline
 $\big(KL^3\big)(123)\big(KL^3\big)^{-1}$ & $(1, -1, -i, 0)$ & $-i\big\lfloor (x_1 - x_2 - ix_3)^2({\rm d}x_1 - {\rm d}x_2 - i{\rm d}x_3){\rm d}x_4, \big(KL^3\big)(123)\big(KL^3\big)^{-1} \big\rceil$\tsep{2pt}\bsep{2pt}\\
 \hline
 $\big(K^2L^3\big)(123)\big(K^2L^3\big)^{-1}$ & $(1, i, -1, 0)$ & $i\big\lfloor (x_1 + ix_2 - x_3)^2({\rm d}x_1 + i{\rm d}x_2 - {\rm d}x_3){\rm d}x_4, \big(K^2L^3\big)(123)\big(K^2L^3\big)^{-1} \big\rceil$\tsep{2pt}\bsep{2pt}\\
 \hline
 $\big(K^3L^2\big)(123)\big(K^3L^2\big)^{-1}$ & $(1, -i, i, 0)$ & $\big\lfloor (x_1 - ix_2 + ix_3)^2({\rm d}x_1 - i{\rm d}x_2 + i{\rm d}x_3){\rm d}x_4, \big(K^3L^2\big)(123)\big(K^3L^2\big)^{-1} \big\rceil$\tsep{2pt}\bsep{2pt}\\
 \hline
 \end{tabular}
 \end{table}

 A similar computation can be made for the two elements $\lfloor (x_1+x_2+x_3)x_4, (123) \rceil$ and $\big\lfloor x_4^2, (123) \big\rceil$, and again for the conjugacy class of $(132)$.

 Recall from the A-model that $\age\, (123) = 1$ and $\age\, (132) = 1$. Also, the polynomials in this case are exactly the same as those from the A-model, where we found $\deg P = 1$ for all such polynomials. One can check that the bidegree for the basis elements in these sectors is $(1, 1)$.

 \textbf{Case 3: $g \in G_W^{\rm diag}$ and has a trivial fixed locus.}
 This case means $g = \big(\tfrac{a_1}{4}, \tfrac{a_2}{4}, \tfrac{a_3}{4}, \tfrac{a_4}{4}\big)$, for $1 \leq a_i \leq 3$. Any sector where $a_1$, $a_2$, $a_3$, and $a_4$ are all nonzero will be narrow. Furthermore, one can easily check that for any $\gamma\in G^\star$, we have $x \in \Fix(g)$ if and only if $\gamma^{-1}\cdot x\in \Fix\big(\gamma g \gamma^{-1}\big)$, so the conjugates of narrow group elements remain narrow.

 Since the sum of $a_1$, $a_2$, $a_3$, and $a_4$ must be a multiple of four, then $(a_1, a_2, a_3, a_4)$ will need to be an ordering of one the following:
 \begin{align*}
(1, 1, 1, 1),\qquad
(2, 2, 2, 2),\qquad
(3, 3, 3, 3),\qquad
(3, 3, 1, 1),\qquad
(3, 2, 2, 1).
 \end{align*}
 The first three lines the above where the components are all equal are powers of~$j_W$. In any of those~$3$ cases, the conjugacy class is trivial since each of them lies in the center of~$G^{\max}_W$.

 There are 12 different orderings of (1, 2, 2, 3). Conjugation by $j_W, K,$ or $L$ does nothing, but conjugation by $(123)$ creates a conjugacy class of size $3$, implying there will be $4$ conjugacy classes of this type. There are 6 orderings of $(1, 1, 3, 3)$, so this choice gives $2$ additional conjugacy classes. The powers of $j_W$ give three more classes. Thus in this case we found a~total of~$9$ conjugacy classes. The sums of the elements in each conjugacy class form a~basis vector for a~narrow sector. A~few examples of these are the following:
 \begin{gather*}
 \lfloor 1, j_W \rceil,\\
 \big\lfloor 1, \big(\tfrac{3}{4}, \tfrac{3}{4}, \tfrac{1}{4}, \tfrac{1}{4}\big) \big\rceil + \big\lfloor 1, \big(\tfrac{3}{4}, \tfrac{1}{4}, \tfrac{3}{4}, \tfrac{1}{4}\big) \big\rceil + \big\lfloor 1, \big(\tfrac{1}{4}, \tfrac{3}{4}, \tfrac{3}{4}, \tfrac{1}{4}\big) \big\rceil,\\
 \big\lfloor 1, \big(\tfrac{2}{4}, \tfrac{2}{4}, \tfrac{3}{4}, \tfrac{1}{4}\big) \rceil + \lfloor 1, \big(\tfrac{3}{4}, \tfrac{2}{4}, \tfrac{2}{4}, \tfrac{1}{4}\big) \big\rceil + \big\lfloor 1, \big(\tfrac{2}{4}, \tfrac{3}{4}, \tfrac{2}{4}, \tfrac{1}{4}\big) \big\rceil.
 \end{gather*}
 The rest are listed in a table at the end of this section. In Section~\ref{s:mirror}, we will show that these $9$ sums of narrow sectors correspond to the $9$ untwisted broad sectors from the A-model.

 Now we turn to the bigrading of these elements.
 All of the polynomials in these sectors will have degree 0, so the bidegree depends just on $\age g$ and $\age g^{-1}$.

 When $g = (j_W)^i$ where $i \in \{1, 2, 3\}$, then $\age (j_W)^i = i$ and $\age (j_W)^{-i} = 4 - i$. Hence the bidegree for these $3$ elements is
 \[
 (i - 1, 3-i).
 \]
 Again, we have three basis elements of bidegree $(0,2)$, $(1,1)$ and $(2,0)$, respectively.

 For the remaining six (sums of) narrow sectors we found above, the associated $g$ and $g^{-1}$ both
 have an age of $2$. Thus the bidegree for all of them will be $(2 - 1, 2 - 1) = (1, 1)$.

 \textbf{Case 4: $g$ is nondiagonal and narrow.} Next we consider the narrow sectors coming from non-diagonal elements of $G^\star$. In particular, $j_W(123)$, $(j_W)^2(123)$, $(j_W)^3(123)$, $j_W(132)$, $(j_W)^2(132)$, and $(j_W)^3(132)$ all have trivial fixed locus as we have seen in Case~3 of Example~\ref{eg:AStateConstruction}. However, on the B-side, these elements have nontrivial conjugacy classes.

 For example, one can check that the conjugacy class of $j_W(123)$ contains $K^iL^jj_W(123)$ for any $i,j\in \{0,1,2,3\}$. The same is true for the other five classes mentioned above. Thus we obtain 16 elements in each conjugacy class. There are~6 such conjugacy classes, and these comprise all of the remaining narrow group elements. Thus we get a contribution of 6 more narrow sectors to the state space.

 We now consider the bidegree of each of these six narrow sectors. Since all of these sectors are narrow, we have $\deg P = 0$, so once again the bidegree depends solely on $\age g$ and $\age g^{-1}$. All of these elements appeared in the A-model as well, where we found that they all have an age of $2$. The inverse of $(j_W)^i(123)^j$ is $(j_W)^{4-i}(123)^{3-j}$, which also has an age of $2$. One can check that each conjugate will also have age~2. Thus the bidegree for all of the elements in this case is $(2 - 1, 2 - 1) = (1, 1)$.

\textbf{Case 5: $g$ has non-trivial fixed locus.} We are left with $42$ elements of $G^\star$.
 One can check that none of these remaining elements contribute to the state space.

In conclusion, the B-model state space contains three basis elements from the untwisted broad sector, six basis elements from the two twisted broad sectors coming from $(123)$ and $(132)$, nine narrow sectors from Case~3, and six more narrow sectors from Case~4, for a total of $24$ basis elements. Recall that there were $24$ basis elements in the A-model as well. If we consider the bigrading, we also see that the dimensions of each graded piece match the A-model.

As with the A-model, we now present of the basis elements in the B-model with their bi\-gra\-ding, seen in Table~\ref{tab:Bmodelbasis}.

\begin{table}[th]\centering
 \begin{longtable}{ l c}\toprule
\multicolumn{2}{c}{B-model basis elements} \\
{B-model basis} &{bidegree} \\
\midrule
\endfirsthead
\toprule
{B-model basis} &{bidegree} \\
\midrule
\endhead
 $\lfloor 1, (0, 0, 0, 0) \rceil$ &$(0, 0)$\tsep{1pt}\\
 $\lfloor x_1x_2x_3x_4, (0, 0, 0, 0) \rceil$ &$(1, 1)$\tsep{1pt}\\
 $\big\lfloor x_1^2x_2^2x_3^2x_4^2, (0, 0, 0, 0) \big\rceil$ &$(2, 2)$\tsep{1pt}\\
 $\big\lfloor (x_1 + x_2 + x_3)^2, (123) \big\rceil+ (15 \text{ elements in conj. class})$ &$(1, 1)$\tsep{1pt}\\
 $\lfloor (x_1 + x_2 + x_3)x_4, (123) \rceil+ (15 \text{ elements in conj. class})$ &$(1, 1)$\tsep{1pt}\\
 $\big\lfloor x_4^2, (123) \big\rceil+ (15 \text{ elements in conj. class})$ &$(1, 1)$\tsep{1pt}\\
 $\big\lfloor (x_1 + x_2 + x_3)^2, (132) \big\rceil+ (15 \text{ elements in conj. class})$ &$(1, 1)$\tsep{1pt}\\
 $\lfloor (x_1 + x_2 + x_3)x_4, (132) \rceil+ (15 \text{ elements in conj. class})$ &$(1, 1)$\tsep{1pt}\\
 $\big\lfloor x_4^2, (132) \big\rceil+ (15 \text{ elements in conj. class})$ &$(1, 1)$\tsep{1pt}\\
 $\lfloor 1, j_W \rceil$ &$(0, 2)$\tsep{1pt}\\
 $\big\lfloor 1, (j_W)^2 \big\rceil$ &$(1, 1)$\tsep{1pt}\\
 $\big\lfloor 1, (j_W)^3 \big\rceil$ &$(2, 0)$\tsep{1pt}\\
 $\big\lfloor 1, \big(\tfrac{3}{4}, \tfrac{3}{4}, \tfrac{1}{4}, \tfrac{1}{4}\big) \big\rceil + \big\lfloor 1, \big(\tfrac{3}{4}, \tfrac{1}{4}, \tfrac{3}{4}, \tfrac{1}{4}\big) \big\rceil + \big\lfloor 1, \big(\tfrac{1}{4}, \tfrac{3}{4}, \tfrac{3}{4}, \tfrac{1}{4}\big) \big\rceil$ &$(1, 1)$\tsep{1pt}\\
 $\big\lfloor 1, \big(\tfrac{2}{4}, \tfrac{2}{4}, \tfrac{3}{4}, \tfrac{1}{4}\big) \big\rceil + \big\lfloor 1, \big(\tfrac{3}{4}, \tfrac{2}{4}, \tfrac{2}{4}, \tfrac{1}{4}\big) \big\rceil + \big\lfloor 1, \big(\tfrac{2}{4}, \tfrac{3}{4}, \tfrac{2}{4}, \tfrac{1}{4}\big) \big\rceil$ &$(1, 1)$\tsep{1pt}\\
 $\big\lfloor 1, \big(\tfrac{2}{4}, \tfrac{2}{4}, \tfrac{1}{4}, \tfrac{3}{4}\big) \big\rceil + \big\lfloor 1, \big(\tfrac{2}{4}, \tfrac{1}{4}, \tfrac{2}{4}, \tfrac{3}{4}\big) \big\rceil + \big\lfloor 1, \big(\tfrac{1}{4}, \tfrac{2}{4}, \tfrac{2}{4}, \tfrac{3}{4}\big) \big\rceil$ &$(1, 1)$\tsep{1pt}\\
 $\big\lfloor 1, \big(\tfrac{2}{4}, \tfrac{3}{4}, \tfrac{1}{4}, \tfrac{2}{4}\big) \big\rceil + \big\lfloor 1, \big(\tfrac{1}{4}, \tfrac{2}{4}, \tfrac{3}{4}, \tfrac{2}{4}\big) \big\rceil + \big\lfloor 1, \big(\tfrac{3}{4}, \tfrac{1}{4}, \tfrac{2}{4}, \tfrac{3}{4}\big) \big\rceil$ &$(1, 1)$\tsep{1pt}\\
 $\big\lfloor 1, \big(\tfrac{3}{4}, \tfrac{2}{4}, \tfrac{1}{4}, \tfrac{2}{4}\big) \big\rceil + \big\lfloor 1, \big(\tfrac{1}{4}, \tfrac{3}{4}, \tfrac{2}{4}, \tfrac{2}{4}\big) \rceil + \big\lfloor 1, \big(\tfrac{2}{4}, \tfrac{3}{4}, \tfrac{1}{4}, \tfrac{3}{4}\big) \big\rceil$ &$(1, 1)$\tsep{1pt}\\
 $\big\lfloor 1, \big(\tfrac{3}{4}, \tfrac{1}{4}, \tfrac{1}{4}, \tfrac{3}{4}\big) \big\rceil + \big\lfloor 1, \big(\tfrac{1}{4}, \tfrac{3}{4}, \tfrac{1}{4}, \tfrac{3}{4}\big) \big\rceil + \big\lfloor 1, \big(\tfrac{1}{4}, \tfrac{1}{4}, \tfrac{3}{4}, \tfrac{3}{4}\big) \big\rceil$ &$(1, 1)$\tsep{1pt}\\
 $\lfloor 1, j_W(123) \rceil + (15 \text{ elements in conj. class})$ &$(1, 1)$\tsep{1pt}\\
 $\big\lfloor 1, (j_W)^2(123) \big\rceil + (15 \text{ elements in conj. class})$ &$(1, 1)$\tsep{1pt}\\
 $\big\lfloor 1, (j_W)^3(123) \big\rceil + (15 \text{ elements in conj. class})$ &$(1, 1)$\tsep{1pt}\\
 $\lfloor 1, j_W(132) \rceil + (15 \text{ elements in conj. class})$ &$(1, 1)$\tsep{1pt}\\
 $\big\lfloor 1, (j_W)^2(132) \big\rceil + (15 \text{ elements in conj. class})$ &$(1, 1)$\tsep{1pt}\\
 $\big\lfloor 1, (j_W)^3(132) \big\rceil + (15 \text{ elements in conj. class})$ &$(1, 1)$\tsep{1pt}\\
 \bottomrule
 \caption{B-model basis elements.}\label{tab:Bmodelbasis}
 \end{longtable}
\end{table}

If we arrange these as a Hodge diamond, we have
\[
\begin{array}{@{}ccc}
 &1&\\
 1 & 20 & 1\\
 &1&
\end{array}
\]
Notice this is the same diamond as with our A-model example. This is enough to prove that the given A- and B-models are isomorphic as bigraded vector spaces, however, we would like to have a canonical isomorphism. This is the goal of the next section.
\end{Example}

As with Lemma~\ref{lem:Abidegree}, we again want to know that the bigrading of an element is unchanged when acted upon by a symmetry in $G^\star,$ so we prove the same fact for B-models.

\begin{Lemma}\label{lem:bidegree}
Given $\gamma \in G_{W^T}^{\max},$ and $\lfloor P, g \rceil \in \mathcal{Q}_{W^T_g}\cdot\omega_{g}$, the element $\gamma \cdot \lfloor P, g \rceil$ has the same bidegree as $\lfloor P, g \rceil$.
\end{Lemma}
\begin{proof}
Recall the B-model bigrading from Definition~\ref{def:Bgrading}:
\[
\big(\deg P + \age g - \age j_{W^T}, \deg P + \age g^{-1} - \age j_{W^T}\big).
\]
This proof follows the same as the proof of Lemma~\ref{lem:Abidegree}. We already proved that $\deg(\gamma \cdot P) = \deg P$ and $\age(\gamma g \gamma^{-1}) = \age g$ in Lemma~\ref{lem:Abidegree}. The work to show that $\age(\gamma g^{-1}\gamma^{-1}) = \age g^{-1}$ is the same, since $\gamma g^{-1}\gamma^{-1}$ and $g^{-1}$ are similar matrices, implying that they too have the same eigenvalues.
Thus $\gamma \cdot \lfloor P, g \rceil$ has the same bidegree as $\lfloor P, g \rceil$.
\end{proof}

\section{The mirror map}\label{s:mirror}

Thus far in our example from Section~\ref{s:statespace}, we have shown that the specified A- and B-models have $24$ basis elements with the same number of elements for each bidegree. While this in itself would be sufficient for claiming that they are isomorphic, we aim to create a canonical map which will better demonstrate which elements on one side correspond to elements on the other. In particular, we expect that this map will exchange narrow and broad sectors. This follows the map given by Krawitz~\cite{KrThesis} for A- and B-models built from abelian groups. This isomorphism between A- and B-models is known as the \emph{mirror map}.
Before we describe the mirror map in general, we will illustrate with our example from Section~\ref{s:statespace}.

\begin{Example}\label{eg:mirrormap}
We will continue with the same A- and B-models as in Examples~\ref{eg:AStateConstruction} and~\ref{eg:BStateConstruction}. To begin constructing the mirror map, we will first look at the part of the map that is already laid out for us by matching the $4$ elements on either side with unique bidegree.
\begin{table}[h!]\centering
 \begin{longtable}{ l ll}\toprule
 \multicolumn{3}{c}{Mirror map: first elements} \\
{bidegree} & {A-model basis} &{B-model basis} \\
 \midrule
 \endfirsthead
 \toprule
{bidegree} & {A-model basis} &{B-model basis} \\
 \midrule
 \endhead
$(0, 0)$ & $\lfloor 1, j_W \rceil $ & $\lfloor 1, (0, 0, 0, 0) \rceil$\tsep{1pt}\\
$(2, 2)$ & $\big\lfloor 1, (j_W)^3 \big\rceil$ & $\big\lfloor x_1^2x_2^2x_3^2x_4^2, (0, 0, 0, 0) \big\rceil$\tsep{1pt}\\
$(0, 2)$ & $\lfloor 1, (0, 0, 0, 0) \rceil $ & $\lfloor 1, j_W \rceil$\tsep{1pt}\\
$(2, 0)$ & $\big\lfloor x_1^2x_2^2x_3^2x_4^2, (0, 0, 0, 0) \big\rceil $ & $\big\lfloor 1, (j_W)^3 \big\rceil$\tsep{1pt}\\
\bottomrule
\end{longtable}
\end{table}

\noindent
This illuminates $2$ more corresponding elements:
\begin{table}[h!]\centering
\begin{longtable}{ l ll}\toprule
\multicolumn{3}{c}{Mirror map: first elements} \\
{bidegree} & {A-model basis} &{B-model basis} \\
\midrule
\endfirsthead
\toprule
{bidegree} & {A-model basis} &{B-model basis} \\
\midrule
\endhead
$(1, 1)$ & $\big\lfloor 1, (j_W)^2 \big\rceil$ & $\lfloor x_1x_2x_3x_4, (0, 0, 0, 0) \rceil$\tsep{1pt}\\
$(1, 1)$& $\lfloor x_1x_2x_3x_4, (0, 0, 0, 0) \rceil$ & $\big\lfloor 1, (j_W)^2 \big\rceil$\tsep{1pt}\\
\bottomrule
\end{longtable}
\end{table}

A nice generalization of the six element maps above can be seen by
\[
\big\lfloor 1, (j_W)^i \big\rceil \leftrightarrow \big\lfloor x_1^{i-1}x_2^{i-1}x_3^{i-1}x_4^{i-1}, (0, 0, 0, 0) \big\rceil.
\]

Recall that the dimension of the untwisted broad sector in the A-model was 9, and since three of those are seen above, there are still $6$ others to account for. These $6$ basis elements map to the $6$ narrow sectors in the B-model which have a conjugacy class of size $3$. Specifically, we map the elements on the A-side whose polynomial has the same permutation structure as the group elements on the B-side. One explicit example is given by mapping the A-model element
\[\big\lfloor x_1^2x_2^2 + x_{1}^2x_{3}^2 + x_{2}^2x_{3}^2, (0, 0, 0, 0) \big\rceil\]
to the B-model element
\[\big\lfloor 1, \big(\tfrac{3}{4}, \tfrac{3}{4}, \tfrac{1}{4}, \tfrac{1}{4}\big) \big\rceil + \big\lfloor 1, \big(\tfrac{3}{4}, \tfrac{1}{4}, \tfrac{3}{4}, \tfrac{1}{4}\big) \big\rceil + \big\lfloor 1, \big(\tfrac{1}{4}, \tfrac{3}{4}, \tfrac{3}{4}, \tfrac{1}{4}\big) \big\rceil.\]
Notice that the term $x_1^2x_2^2$ has a power of $2$ for $x_1$ and $x_2$, and this corresponds to the first two components of the group element of $\big\lfloor 1, \big(\tfrac{3}{4}, \tfrac{3}{4}, \tfrac{1}{4}, \tfrac{1}{4}\big) \big\rceil$ having a larger value by $\tfrac{2}{4}$. The same correspondence can be noticed between $x_{1}^2x_{3}^2$ and $\big\lfloor 1, \big(\tfrac{3}{4}, \tfrac{1}{4}, \tfrac{3}{4}, \tfrac{1}{4}\big) \big\rceil$, as well as $x_{2}^2x_{3}^2$ and $\big\lfloor 1, \big(\tfrac{1}{4}, \tfrac{3}{4}, \tfrac{3}{4}, \tfrac{1}{4}\big) \big\rceil$.

All $6$ elements of this type are given below in Table~\ref{tab:mirrormappermutations}. The bidegree is left out, but all of the following elements have a bidegree of $(1, 1)$.
\begin{table}[h!]\centering
\begin{longtable}{ll} \toprule
 \multicolumn{2}{c}{Mirror map: basis elements following permutations structure} \\
{A-model} &{B-model} \\
 \midrule
 \endfirsthead
 \toprule
{A-model} &{B-model} \\
 \midrule
 \endhead
$\big\lfloor x_1^2x_2^2 + x_{1}^2x_{3}^2 + x_{2}^2x_{3}^2, (0, 0, 0, 0) \big\rceil $& $\big\lfloor 1, \big(\tfrac{3}{4}, \tfrac{3}{4}, \tfrac{1}{4}, \tfrac{1}{4}\big) \big\rceil + \big\lfloor 1, \big(\tfrac{3}{4}, \tfrac{1}{4}, \tfrac{3}{4}, \tfrac{1}{4}\big) \big\rceil + \big\lfloor 1, \big(\tfrac{1}{4}, \tfrac{3}{4}, \tfrac{3}{4}, \tfrac{1}{4}\big) \big\rceil$\tsep{1pt}\\
$\big\lfloor x_1x_2x_3^2 + x_1^2x_2x_3 + x_1x_2^2x_3, (0, 0, 0, 0) \big\rceil$ & $\big\lfloor 1, \big(\tfrac{2}{4}, \tfrac{2}{4}, \tfrac{3}{4}, \tfrac{1}{4}\big) \big\rceil + \big\lfloor 1, \big(\tfrac{3}{4}, \tfrac{2}{4}, \tfrac{2}{4}, \tfrac{1}{4}\big) \big\rceil + \big\lfloor 1, \big(\tfrac{2}{4}, \tfrac{3}{4}, \tfrac{2}{4}, \tfrac{1}{4}\big) \big\rceil$\tsep{1pt}\\
$\big\lfloor x_1x_2x_4^2 + x_1x_3x_4^2 + x_2x_3x_4^2, (0, 0, 0, 0) \big\rceil$ & $\big\lfloor 1, \big(\tfrac{2}{4}, \tfrac{2}{4}, \tfrac{1}{4}, \tfrac{3}{4}\big) \big\rceil + \big\lfloor 1, \big(\tfrac{2}{4}, \tfrac{1}{4}, \tfrac{2}{4}, \tfrac{3}{4}\big) \big\rceil + \big\lfloor 1, \big(\tfrac{1}{4}, \tfrac{2}{4}, \tfrac{2}{4}, \tfrac{3}{4}\big) \big\rceil$\tsep{1pt}\\
$\big\lfloor x_1x_2^2x_4 + x_2x_3^2x_4 + x_1^2x_3x_4, (0, 0, 0, 0) \big\rceil$ & $\big\lfloor 1, \big(\tfrac{2}{4}, \tfrac{3}{4}, \tfrac{1}{4}, \tfrac{2}{4}\big) \big\rceil + \big\lfloor 1, \big(\tfrac{1}{4}, \tfrac{2}{4}, \tfrac{3}{4}, \tfrac{2}{4}\big) \big\rceil + \big\lfloor 1, \big(\tfrac{3}{4}, \tfrac{1}{4}, \tfrac{2}{4}, \tfrac{2}{4}\big)\big \rceil$\tsep{1pt}\\
$\big\lfloor x_1^2x_2x_4 + x_2^2x_3x_4 + x_1x_3^2x_4, (0, 0, 0, 0) \big\rceil$ & $\big\lfloor 1, \big(\tfrac{3}{4}, \tfrac{2}{4}, \tfrac{1}{4}, \tfrac{2}{4}\big) \big\rceil + \big\lfloor 1, \big(\tfrac{1}{4}, \tfrac{3}{4}, \tfrac{2}{4}, \tfrac{2}{4}\big) \big\rceil + \big\lfloor 1, \big(\tfrac{2}{4}, \tfrac{1}{4}, \tfrac{3}{4}, \tfrac{2}{4}\big) \big\rceil$\tsep{1pt}\\
$\big\lfloor x_1^2x_4^2 + x_2^2x_4^2 + x_3^2x_4^2, (0, 0, 0, 0) \big\rceil$ & $\big\lfloor 1, \big(\tfrac{3}{4}, \tfrac{1}{4}, \tfrac{1}{4}, \tfrac{3}{4}\big) \big\rceil + \big\lfloor 1, \big(\tfrac{1}{4}, \tfrac{3}{4}, \tfrac{1}{4}, \tfrac{3}{4}\big) \big\rceil + \big\lfloor 1, \big(\tfrac{1}{4}, \tfrac{1}{4}, \tfrac{3}{4}, \tfrac{3}{4}\big) \big\rceil$\tsep{1pt}\\
\bottomrule
\caption{The basis elements following permutation structure.}\label{tab:mirrormappermutations}
\end{longtable}
\end{table}

Again, notice that the polynomials of the A-model elements have the same permutation structure of the group elements in the B-model.

There are now $12$ basis elements left to be mapped in both models, with $6$ being twisted broad sectors from $g = (123)$ or $g = (132)$ and $6$ being narrow sectors, where $g$ is a product of a permutation and a power of $j_W$. While they all have the same bidegree, we expect that the mirror map will map broad sectors to narrow sectors and narrow sectors to broad sectors, so we will do the same here. Unlike the previous element mappings, it is not as clear exactly which A-model and B-model elements below should map to each other. We give the correspondence in Table~\ref{tab:mirrormap_remaining}.

\begin{table}[h!]\centering
\begin{longtable}{ll} \toprule
 \multicolumn{2}{c}{Mirror map, remaining basis elements} \\
{A-model} &{B-model} \\
 \midrule
$\big\lfloor (x_1 + x_2 + x_3)^2, (123) \big\rceil $ & $\lfloor 1, (123)j_W \rceil + (15 \text{ others})$\tsep{1pt}\\
$\lfloor (x_1 + x_2 + x_3)x_4, (123) \rceil$ & $\big\lfloor 1, (123)(j_W)^2 \big\rceil + (15 \text{ others})$\tsep{1pt}\\
$\big\lfloor (x_4)^2, (123) \big\rceil \hspace{1cm}$ & $\big\lfloor 1, (123)(j_W)^3 \big\rceil + (15 \text{ others})$\tsep{1pt}\\
$\big\lfloor (x_1 + x_2 + x_3)^2, (132) \big\rceil $ & $\lfloor 1, (132) j_W \rceil + (15 \text{ others})$\tsep{1pt}\\
$\lfloor (x_1 + x_2 + x_3)x_4, (132) \rceil $ & $\big\lfloor 1, (132)(j_W)^2 \big\rceil + (15 \text{ others})$\tsep{1pt}\\
$\big\lfloor (x_4)^2, (132) \big\rceil $ & $\big\lfloor 1, (132)(j_W)^3 \big\rceil + (15 \text{ others})$\tsep{1pt}\\
$\lfloor 1, (123)j_W \rceil $ & $\big\lfloor (x_1 + x_2 + x_3)^2, (123) \big\rceil+ (15 \text{ others})$\tsep{1pt}\\
$\big\lfloor 1, (123)(j_W)^2 \big\rceil$ & $\lfloor (x_1 + x_2 + x_3)x_4, (123) \rceil+ (15 \text{ others})$\tsep{1pt}\\
$\big\lfloor 1, (123)(j_W)^3 \big\rceil $ & $\big\lfloor (x_4)^2, (123) \big\rceil+ (15 \text{ others})$\tsep{1pt}\\
$\lfloor 1, (132) j_W \rceil$ & $\big\lfloor (x_1 + x_2 + x_3)^2, (132) \big\rceil+ (15 \text{ others})$\tsep{1pt}\\
$\big\lfloor 1, (132)(j_W)^2 \big\rceil$ & $ \lfloor (x_1 + x_2 + x_3)x_4, (132) \rceil+ (15 \text{ others})$\tsep{1pt}\\
$\big\lfloor 1, (132)(j_W)^3 \big\rceil$ & $\big\lfloor (x_4)^2, (132) \big\rceil+ (15 \text{ others})$\tsep{1pt}\\
\bottomrule
\caption{Mirror map, remaining basis elements.}\label{tab:mirrormap_remaining}
\end{longtable}
\end{table}

This completes the mirror map for this example, and we have explicitly shown in this example that as bigraded vector spaces,
$
\mathcal{A}_{W, G} \cong \mathcal{B}_{W^T, G^\star}.
$
\end{Example}

We will now generalize what we have seen in this example.

\subsection{Proof of the mirror map}

In this section we will generalize what we have observed in the previous example to prove a~general theorem, providing a canonical mirror map on certain natural subspaces.
Later we will find under certain conditions the A- and B-models are surprisingly not isomorphic as bigraded vector spaces. However, the restriction to the subspaces mentioned above holds nevertheless.

\newpage

\begin{Theorem}\label{thm:mirror}
Let $W$ be an invertible Fermat polynomial and $G \subseteq G_W^{\max}$ be an admissible group of the form $H\cdot K$, where $K \leq G$ is the subgroup of pure even permutations and $H \leq G$ is the subgroup of diagonal symmetries. Define $\mathcal{A}_0 \subseteq \mathcal{A}_{W, G}$ and $\mathcal{B}_0 \subseteq \mathcal{B}_{W^T, G^\star}$ to be the untwisted broad sectors for the A- and B-side, respectively. Let $nar' \leq H$ be the set of narrow diagonal symmetries. We will also denote $nar' \leq H^T$ to be the corresponding set on the B-side. Then there exist bigraded vector space isomorphisms
\[
\mathcal{A}_0 \xrightarrow{\sim} \mathcal{B}_{nar'} \qquad \text{and} \qquad \mathcal{A}_{nar'} \xrightarrow{\sim} \mathcal{B}_0.
\]
\end{Theorem}

\begin{Remark}We have restricted our attention to Fermat polynomials for simplicity. This seems to be the most natural setting for considering nonabelian symmetry groups, as the permutations allowed by the structure of $G^{\max}_W$ are somewhat restrictive. However, we expect a similar theorem to hold for more general invertible polynomials as well.
\end{Remark}

\begin{proof}[Proof of Theorem~\ref{thm:mirror}]
Let $W = x_1^{d_1} + \dots + x_N^{d_N}$, so $G_W^{\rm diag}$ is generated by the set
\[
\big\{\big(\tfrac{1}{d_1}, 0, \dots, 0\big), \big(0, \tfrac{1}{d_2}, \dots, 0\big), \dots, \big(0, 0, \dots, \tfrac{1}{d_N}\big)\big\}.
\]

Given $g=(a_1,\dots, a_N)\in G\cap G^{\rm diag}_W$, define $I_g = \{i \in \{1, \dots, N\} \,|\, a_i \neq 0 \}$ and consider the map
\[
\bigoplus_{g \in G \cap G_W^{\rm diag}} \cQ_{W_g} \cdot \omega_g \rightarrow \bigoplus_{g' \in G^\star \cap G_{W^T}^{\rm diag}} \cQ_{W^T_{g'}} \cdot \omega_{g'}
\]
given by
\begin{gather}\label{eq:mirrormap}
\Bigg \lfloor \bigwedge_{j \not\in I_g} x_i^{b_i} {\rm d}x_i, \big(\tfrac{a_1}{d_1}, \dots, \tfrac{a_N}{d_N}\big) \Bigg \rceil \mapsto \Bigg \lfloor \bigwedge_{j \in I_g} y_j^{a_j - 1} {\rm d}x_j, \big(\tfrac{b'_1}{d_1}, \dots, \tfrac{b'_N}{d_N}\big) \Bigg \rceil,
\end{gather}
for $g = \big(\tfrac{a_1}{d_1}, \dots, \tfrac{a_N}{d_N}\big),$ where $b'_i = b_i + 1$ if $i \not\in I_g$, and $b'_i = 0$ otherwise. Notice that although $W=W^T$ in this setting, we will continue to write $W^T$ to keep clear in our minds which side of the mirror map we are considering.

The map described in~\eqref{eq:mirrormap} is known as the \emph{map on the unprojected state spaces}, where invariance has not yet been considered. Notice that $G^\star \cap G_{W^T}^{\rm diag} = G^T$ in this case. This map was proven to be a bijection by Krawitz in~\cite{KrThesis}. For completeness, however, we will reprove the relevant part here, namely we will show that
\begin{gather}\label{eq:rmirrormap}
\cQ_W \rightarrow \bigoplus_{\substack{g' \in G_{W^T}^{\rm diag}\\ \Fix(g') = \{0\}}} \cQ_{W^T_{g'}}\cdot \omega_{g'}
\end{gather}
is a bijection. This is the case in \eqref{eq:mirrormap} for when $a_i = 0$ for all $i$ and $\big(\tfrac{b'_1}{d_1}, \dots, \tfrac{b'_N}{d_N}\big)$ is a diagonal symmetry with nonzero entries. We will also need that fact that
\[
\bigoplus_{\substack{g \in G_W^{\rm diag}\\ \Fix(g) = 0}} \cQ_{W_g}\cdot \omega_g \rightarrow\cQ_{W^T}
\]
is a bijection (this time as a map from the A-model to the B-model). However, the proof of this exactly mirrors the first one, so we will exclude it here. Before proceeding, note that since $W$ is Fermat, we know that the Milnor ring of~$W$ is{\samepage
\[
\cQ_W = \frac{\mathbb{C}[x_1, \dots, x_N]}{\big(x_1^{d_1 - 1}, \dots, x_N^{d_N - 1}\big)},
\]
which has a basis of elements of the form $\prod\limits_{i=1}^{N} x_i^{b_i}$, where $0 \leq b_i \leq d_i - 2$.}

First, to prove surjectivity of \eqref{eq:rmirrormap}, let
\[
\big\lfloor 1, \big(\tfrac{b'_1}{d_1}, \dots, \tfrac{b'_N}{d_N}\big) \big\rceil \in \bigoplus_{\substack{g' \in G_{W^T}^{\rm diag}\\ \Fix(g') = \{0\}}} \cQ_{W^T_{g'}}\cdot \omega_{g'},
\]
so $1 \leq b'_i \leq d_i - 1$. Since each $b'_i \neq 0$, let $b_i = b'_i - 1$, and notice $0 \leq b_i \leq d_i - 2$. The preimage of $\big\lfloor 1, \big(\tfrac{b'_1}{d_1}, \dots, \tfrac{b'_N}{d_N}\big) \big\rceil$ is
\[
\Bigg\lfloor \bigwedge_{i=1}^{N} x_i^{b_i} {\rm d}x_i, (0, \dots, 0) \Bigg\rceil.
\]
Notice $W|_{\Fix((0, \dots, 0))} = W$, so $\prod\limits_{i=1}^{N} x_i^{b_i} \in \cQ_W$, as desired.

To prove injectivity, let
\[
\big\lfloor 1, \big(\tfrac{b'_1}{d_1}, \dots, \tfrac{b'_N}{d_N}\big) \big\rceil = \big\lfloor 1, \big(\tfrac{c'_1}{d_1}, \dots, \tfrac{c'_N}{d_N}\big) \big\rceil \in \bigoplus_{\substack{g' \in G_{W^T}^{\rm diag}\\ \Fix(g') = \{0\}}} \cQ_{W^T_{g'}}\cdot \omega_{g'}
\]
be two elements in the image of the given map. This would imply that $b'_i = c'_i$ (mod $d_i$) for all $i$,
so $\prod\limits_{i=1}^{n} x_i^{b_i} = \prod\limits_{i=1}^{n} x_i^{c_i}$ in $\cQ_W$.
Thus their preimages are equal:
\[
\Bigg\lfloor \bigwedge_{i=1}^{n} x_i^{b_i} {\rm d}x_i, (0, \dots, 0) \Bigg\rceil = \Bigg\lfloor \bigwedge_{i=1}^{n} x_i^{c_i} {\rm d}x_i, (0, \dots, 0) \Bigg\rceil.
\]
So \eqref{eq:rmirrormap} is also injective.

Next, we look at the invariant subspaces of the preimage and image of the above map. Specifically, we aim to show that
\[
\mathcal{A}_0 \rightarrow \mathcal{B}_{nar'}
\]
is a bijection under \eqref{eq:rmirrormap}, where
\[
\mathcal{A}_0 = \big(\cQ_W\big)^G \qquad\text{and} \qquad \mathcal{B}_{nar'} =   \Bigg(\bigoplus_{\substack{g' \in G_{W^T}^{\rm diag}\\ \Fix(g') = \{0\}}} \cQ_{W^T_{g'}}\cdot \omega_{g'}\Bigg)^{G^\star}.
\]

We will need to show that $\Big\lfloor \sum\limits_{r=1}^m \Big(\bigwedge\limits_{i =1}^N x_i^{b_{i,r}} {\rm d}x_i\Big), (0, \dots, 0)\Big\rceil$ is invariant under $G$ if and only if its image is fixed by the elements of $G^\star$. Since $G$ is generated by~$K$ and~$H$, we
can consider all $\sigma \in K$ and $h \in H$ separately. In other words, we will prove this statement in two parts. First, we will show that $\Big\lfloor \sum\limits_{r=1}^m \Big(\bigwedge\limits_{i =1}^N x_i^{b_{i,r}} {\rm d}x_i\Big), (0, \dots, 0)\Big\rceil$ is fixed by $H$ if and only if its image is fixed by~$H^T$, and secondly the same element is invariant under $K$ if and only if its image is fixed by~$K$.

\textbf{Part 1:} Let $h \in H$, so $h$ is a diagonal symmetry of the form $\big(\tfrac{h_1}{d_1}, \dots, \tfrac{h_N}{d_N}\big)$. Note that, since~$h$ acts diagonally, it fixes terms of polynomials independently,
so we only need to consider~$P$ as a~monomial. That is, we need to consider
$P = \bigwedge\limits_{i=1}^N x_i^{b_{i}} {\rm d}x_i$ fixed by $h$.

We can see that
\[
h \cdot \lfloor P, (0, \dots, 0) \rceil = \lfloor h \cdot P, (0, \dots, 0) \rceil = \bigg\lfloor {\rm e}^{2\pi {\rm i} \sum\limits_{i = 1}^N \tfrac{h_i b'_i}{d_i}}P, (0, \dots, 0) \bigg\rceil.
\]

The monomial $P$ is fixed by $h$ if and only if $\sum\limits_{i = 1}^N \tfrac{h_i b'_i}{d_i} \in \mathbb{Z}$, which we can rewrite as
\[
\sum_{i = 1}^N \tfrac{h_i b'_i}{d_i} = h A_W\big(\tfrac{b'_1}{d_1}, \dots, \tfrac{b'_N}{d_N}\big) \in \mathbb{Z}.
\]
By definition, this happens if and only if $\big(\tfrac{b'_1}{d_1}, \dots, \tfrac{b'_N}{d_N}\big) \in H^T$, so
\[
\big\lfloor 1, \big(\tfrac{b'_1}{d_1}, \dots, \tfrac{b'_N}{d_N}\big) \big\rceil \in \Bigg(\bigoplus_{\substack{g \in G_{W^T}^{\rm diag}\\ \Fix(g) = \{0\}}} \cQ_{W^T|_{\Fix(g)}}\cdot \omega_{g}\Bigg)^{H^T}.
\]
This proves the first part.

\textbf{Part 2:} Let $\sigma \in K$, and assume that $\sigma$ fixes $\Big\lfloor \sum\limits_{r=1}^m \Big(\bigwedge\limits_{i=1}^N x_i^{b_{i,r}} {\rm d}x_i\Big), (0, \dots, 0)\Big\rceil \in \mathcal{A}_0$. That is, if $\bigwedge\limits_{i = 1}^N x_i^{b_i} {\rm d}x_i$ is a single term of the sum $\sum\limits_{r=1}^m \Big(\bigwedge\limits_{i=1}^N x_i^{b_{i,r}} {\rm d}x_i\Big)$, then
\[
\sigma\cdot \bigwedge_{i = 1}^N x_i^{b_i} {\rm d}x_i = \bigwedge_{i = 1}^N x_i^{b_\sigma(i)} {\rm d}x_i
\]
must be another term in the sum. Note that $\sigma$ fixes the volume form because $\sigma$ is an even permutation.

Now consider $\sum\limits_{r=1}^m \big\lfloor 1, \big(\tfrac{b'_{1,r}}{d_1}, \dots, \tfrac{b'_{N,r}}{d_N}\big) \big\rceil \in \mathcal{B}_{nar'}$, which is the image of
\[
\Bigg\lfloor \sum_{r=1}^m \Bigg(\bigwedge_{i=1}^N x_i^{b_{i,r}} {\rm d}x_i\Bigg), (0, \dots, 0)\Bigg\rceil \in \mathcal{A}_0.
\]
Notice
\[
\sigma\cdot \big\lfloor 1, \big(\tfrac{b'_1}{d_1}, \dots, \tfrac{b'_N}{d_N}\big) \big\rceil = \big\lfloor 1, \big(\tfrac{b'_{\sigma(1)}}{d_{\sigma(1)}}, \dots, \tfrac{b'_{\sigma(n)}}{d_{\sigma(N)}}\big) \big\rceil
\]
is the image of
\[
\sigma\cdot\bigwedge_{i = 1}^N x_i^{b_i} {\rm d}x_i = \bigwedge_{i = 1}^N x_i^{b_\sigma(i)} {\rm d}x_i.
\]

Thus $\sigma$ fixes $\sum\limits_{r=1}^m \big\lfloor 1, \big(\tfrac{b'_{1,r}}{d_1}, \dots, \tfrac{b'_{N,r}}{d_N}\big) \big\rceil \in \mathcal{B}_{nar'}$ if and only if it fixes $\bigwedge\limits_{i = 1}^N x_i^{b_i} {\rm d}x_i$ as well.

The work to show that $\mathcal{A}_{nar'} \rightarrow \mathcal{B}_0$ is a bijection follows similarly. Thus we have proved that the A- and B-models restricted to these important subspaces are isomorphic as vector spaces via a canonical isomorphism. It remains to show that the corresponding elements from either side have the same bidegree.

Recall that the A-model bigrading from Definition~\ref{def:Agrading} is
\[
(\deg P + \age g - \age j_W, N_g - \deg P + \age g - \age j_W).
\]
If we restrict to $\mathcal{A}_0$, then $\age g = 0$ and $N_g = N$, so the above definition reduces to
\[
(\deg P - \age j_W, N - \deg P - \age j_W).
\]
The B-model bigrading from Definition~\ref{def:Bgrading} was
\[
\big(\deg P' + \age g' - \age j_W, \deg P' + \age (g')^{-1} - \age j_W\big).
\]
When we consider elements of $\mathcal{B}_{nar'}$, we find that $\deg P' = 0$, so the bigrading becomes
\[
\big(\age g' - \age j_W,\age (g')^{-1} - \age j_W\big).
\]

Consider the corresponding elements
\[
\Bigg\lfloor \sum_{r=1}^m \Bigg(\bigwedge_{i=1}^N x_i^{b_{i,r}} {\rm d}x_i\Bigg), (0, \dots, 0)\Bigg\rceil \in \mathcal{A}_0\qquad\text{and} \qquad \sum_{r=1}^m \big\lfloor 1, \big(\tfrac{b'_{1,r}}{d_1}, \dots, \tfrac{b'_{N,r}}{d_N}\big) \big\rceil \in \mathcal{B}_{nar'}.
\]
By Lemmas~\ref{lem:Abidegree} and~\ref{lem:bidegree}, we only need to focus on one term in each sum. Thus to show that the mirror map preserves bidegree, we must prove that $\deg \Big(\bigwedge\limits_{i=1}^N x_i^{b_{i,r}} {\rm d}x_i\Big) = \age \big(\tfrac{b'_{1,r}}{d_1}, \dots, \tfrac{b'_{N,r}}{d_N}\big)$ and $ N - \deg \Big(\prod\limits_{i=1}^N x_i^{b_{i,r}} {\rm d}x_i\Big) = \age \big(\tfrac{b'_{1,r}}{d_1}, \dots, \tfrac{b'_{N,r}}{d_N}\big)^{-1}.$ Observe that
\[
\deg\Bigg(\bigwedge_{j=1}^N x_j^{b_i} {\rm d}x_j\Bigg) = \sum_{i =1}^N \frac{b_i + 1}{d_i}.
\]
and
\[
\age \big(\tfrac{b'_1}{d_1}, \dots, \tfrac{b'_N}{d_N}\big)
=\sum_{i=1}^N\frac{b'_i}{d_i}.
\]

It follows that $N - \deg \Big(\bigwedge\limits_{i=1}^N x_i^{b_{i,r}} {\rm d}x_i\Big) = \age \big(\tfrac{b'_{1,r}}{d_1}, \dots, \tfrac{b'_{N,r}}{d_N}\big)^{-1}$ since if $g$ is narrow, then it is known that $\age g = N - \age g^{-1}$ (Mukai~\cite{Mukai}). This establishes that the first isomorphism $\mathcal{A}_0 \rightarrow \mathcal{B}_{nar'}$ preserves bidegree.

For the other isomorphism, $\mathcal{A}_{nar'} \rightarrow \mathcal{B}_0,$ we need to show is that the corresponding elements from $\mathcal{A}_{nar'}$ and $\mathcal{B}_0$ also have the same bidegree. The bigrading of elements from these sectors is
\[
(\age g - \age j_W, \age g - \age j_W) \qquad \text{and} \qquad (\deg P' - \age j_W, \deg P' - \age j_W),
\]
respectively, where $g \in G$ and $P' \in (\cQ_W)^{G^\star}$. This means that all we need to show is that $\age g = \deg P'$, which follows the exact same work as above.

Thus we have shown that the maps
\[
\mathcal{A}_0 \rightarrow \mathcal{B}_{nar'} \qquad \text{and} \qquad \mathcal{B}_0 \rightarrow \mathcal{A}_{nar'}
\]
are bigraded vector space isomorphisms. This gives us the partial mirror map.
\end{proof}

\section{Two examples}

In this section, we give two more examples that are perhaps more illuminating than Example~\ref{eg:mirrormap}. The first is an example where A- and B-models are isomorphic as bigraded vector spaces, and the second is an example where the bigrading fails.

\subsection{Good example}

While the example we began in Section~\ref{s:statespace} was a great starting place, twenty of the $24$ basis elements had the same bidegree of $(1, 1)$. Moving up to a higher degree polynomial will create~A- and B-models with larger bases and more variety in their bigrading, illuminating a clearer picture of the mirror map. With Theorem~\ref{thm:mirror}, we know what most of the map will look like.

\begin{Example}\label{eg:good}
Let $W = x_1^5 + x_2^5 + x_3^5 + x_4^5 + x_5^5$ and $G = \langle j_W, (12)(34)\rangle$, where
\[
j_W = \big(\tfrac{1}{5}, \tfrac{1}{5}, \tfrac{1}{5}, \tfrac{1}{5}, \tfrac{1}{5}\big)
\qquad \text{and}\qquad (12)(34) =
\begin{pmatrix}
 0 & 1 & 0 & 0 & 0\\
 1 & 0 & 0 & 0 & 0\\
 0 & 0 & 0 & 1 & 0\\
 0 & 0 & 1 & 0 & 0\\
 0 & 0 & 0 & 0 & 1\\
\end{pmatrix}.
\]
Then $W^T = W$ and the non-abelian dual group of $G$ is
$
G^\star = \SL_{W^T}^{\rm diag}\cdot \langle (12)(34) \rangle,
$
where $\SL_W^{\rm diag} = \big\langle \big(\tfrac{1}{5}, \tfrac{1}{5}, \tfrac{1}{5}, \tfrac{1}{5}, \tfrac{1}{5}\big), \big(\tfrac{2}{5}, \tfrac{1}{5}, \tfrac{1}{5}, \tfrac{1}{5}, 0\big), \big(\tfrac{1}{5}, \tfrac{2}{5}, \tfrac{1}{5}, \tfrac{1}{5}, 0\big), \big(\tfrac{1}{5}, \tfrac{1}{5}, \tfrac{2}{5}, \tfrac{1}{5}, 0\big) \big\rangle.$ We will denote $K = \big(\tfrac{2}{5}, \tfrac{1}{5}, \tfrac{1}{5}, \tfrac{1}{5}, 0\big)$, $L= \big(\tfrac{1}{5}, \tfrac{2}{5}, \tfrac{1}{5}, \tfrac{1}{5}, 0\big)$, and $M = \big(\tfrac{1}{5}, \tfrac{1}{5}, \tfrac{2}{5}, \tfrac{1}{5}, 0\big)$.

As with the previous example, the goal is to show that
\[\mathcal{A}_{W, G} \cong \mathcal{B}_{W^T, G^\star}\]
as bigraded vector spaces. This example follows the same recipe as Example~\ref{eg:mirrormap}, so we leave the details to the reader and simply provide the mirror map.

The first eight elements listed in Table~\ref{tab:mirrormapcenter} correspond to the elements indexed by $j_W$ (on either the A-side or the B-side). Recall that~$j_W$ lies in the center of~$G$.

\begin{longtable}{lll} \toprule
 \multicolumn{3}{c}{Mirror map: first basis elements} \\
{bidegree}& {A-model} &{B-model} \\
 \midrule
 \endfirsthead
 \toprule
{bidegree}& {A-model} &{B-model} \\
 \midrule
 \endhead
$(0, 0)$ & $\lfloor 1, j_W \rceil $& $\lfloor 1, (0, 0, 0, 0, 0) \rceil$\tsep{1pt}\\
$(1, 1)$  & $\big\lfloor 1, (j_W)^2 \big\rceil$& $\lfloor x_1x_2x_3x_4, (0, 0, 0, 0, 0) \rceil$\tsep{1pt}\\
$(2, 2)$  & $\big\lfloor 1, (j_W)^3 \big\rceil $& $\big\lfloor x_1^2x_2^2x_3^2x_4^2, (0, 0, 0, 0, 0) \big\rceil$\tsep{1pt}\\
$(3, 3)$  & $\big\lfloor 1, (j_W)^4 \big\rceil$& $\big\lfloor x_1^3x_2^3x_3^3x_4^3, (0, 0, 0, 0, 0) \big\rceil$\tsep{1pt}\\
$(0, 3)$  & $\lfloor 1, (0, 0, 0, 0, 0) \rceil $& $\lfloor 1, j_W \rceil$\tsep{1pt}\\
$(1, 2)$  & $\lfloor x_1x_2x_3x_4, (0, 0, 0, 0, 0) \rceil $& $\big\lfloor 1, (j_W)^2 \big\rceil$\tsep{1pt}\\
$(2, 1)$   & $\big\lfloor x_1^2x_2^2x_3^2x_4^2, (0, 0, 0, 0, 0) \big\rceil$& $\big\lfloor 1, (j_W)^3 \big\rceil$\tsep{1pt}\\
$(3, 0)$  &$ \big\lfloor x_1^3x_2^3x_3^3x_4^3, (0, 0, 0, 0, 0) \big\rceil $& $\big\lfloor 1, (j_W)^4 \big\rceil$\tsep{1pt}\\
\bottomrule
\caption{Mirror map: first few basis elements.}\label{tab:mirrormapcenter}
\end{longtable}

The next two tables contain the basis elements in $\mathcal{A}_0$ and the corresponding narrow sectors in $\mathcal{B}_{nar'}$ as in Theorem~\ref{thm:mirror}. The elements in Table~\ref{tab:degree11} all have a bidegree of~$(1, 2),$ while those in Table~\ref{tab:degree21} have a bidegree of~$(2, 1)$.



{\small

\begin{longtable}{ll} \toprule
\multicolumn{2}{c}{Mirror map: basis elements of degree $(1,2)$} \\
{A-model} &{B-model}	\\
	\midrule
	\endfirsthead
	\toprule
	\multicolumn{2}{c}{Basis elements of degree $(1,2)$} \\
{A-model} &{B-model}	\\
	 \midrule
	\endhead
$\big\lfloor x_1^3x_2^2 + x_1^2x_2^3, (0, 0, 0, 0, 0) \big\rceil  $ & $\big\lfloor 1, \big(\tfrac{4}{5}, \tfrac{3}{5}, \tfrac{1}{5}, \tfrac{1}{5}, \tfrac{1}{5}\big) \big\rceil + \big\lfloor 1, \big(\tfrac{3}{5}, \tfrac{4}{5}, \tfrac{1}{5}, \tfrac{1}{5}, \tfrac{1}{5}\big) \big\rceil$\tsep{1pt}\\
$\big\lfloor x_3^3x_4^2 + x_3^2x_4^3, (0, 0, 0, 0, 0) \big\rceil $ & $\big\lfloor 1, \big(\tfrac{1}{5}, \tfrac{1}{5}, \tfrac{4}{5}, \tfrac{3}{5}, \tfrac{1}{5}\big) \big\rceil + \big\lfloor 1, \big(\tfrac{1}{5}, \tfrac{1}{5}, \tfrac{3}{5}, \tfrac{4}{5}, \tfrac{1}{5}\big) \big\rceil$\tsep{1pt}\\
$\big\lfloor x_1^2x_3^3 + x_2^2x_4^3, (0, 0, 0, 0, 0) \big\rceil  $ & $\big\lfloor 1, \big(\tfrac{3}{5}, \tfrac{1}{5}, \tfrac{4}{5}, \tfrac{1}{5}, \tfrac{1}{5}\big) \big\rceil + \big\lfloor 1, \big(\tfrac{1}{5}, \tfrac{3}{5}, \tfrac{1}{5}, \tfrac{4}{5}, \tfrac{1}{5}\big) \big\rceil$\tsep{1pt}\\
$\big\lfloor x_1^3x_3^2 + x_2^3x_4^2, (0, 0, 0, 0, 0) \big\rceil  $ &$ \big\lfloor 1, \big(\tfrac{4}{5}, \tfrac{1}{5}, \tfrac{3}{5}, \tfrac{1}{5}, \tfrac{1}{5}\big) \big\rceil + \big\lfloor 1, \big(\tfrac{1}{5}, \tfrac{4}{5}, \tfrac{1}{5}, \tfrac{3}{5}, \tfrac{1}{5}\big) \big\rceil$\tsep{1pt}\\
$\big\lfloor x_1^3x_4^2 + x_2^3x_3^2, (0, 0, 0, 0, 0) \big\rceil  $ & $\big\lfloor 1, \big(\tfrac{4}{5}, \tfrac{1}{5}, \tfrac{1}{5}, \tfrac{3}{5}, \tfrac{1}{5}\big) \big\rceil + \big\lfloor 1, \big(\tfrac{1}{5}, \tfrac{4}{5}, \tfrac{3}{5}, \tfrac{1}{5}, \tfrac{1}{5}\big) \big\rceil$\tsep{1pt}\\
$\big\lfloor x_1^2x_4^3 + x_2^2x_3^3, (0, 0, 0, 0, 0) \big\rceil  $ & $\big\lfloor 1, \big(\tfrac{3}{5}, \tfrac{1}{5}, \tfrac{1}{5}, \tfrac{4}{5}, \tfrac{1}{5}\big) \rceil + \big\lfloor 1, \big(\tfrac{1}{5}, \tfrac{3}{5}, \tfrac{4}{5}, \tfrac{1}{5}, \tfrac{1}{5}\big) \big\rceil$\tsep{1pt}\\
$\big\lfloor x_1^3x_5^2 + x_2^3x_5^2, (0, 0, 0, 0, 0) \big\rceil  $ & $\big\lfloor 1, \big(\tfrac{4}{5}, \tfrac{1}{5}, \tfrac{1}{5}, \tfrac{1}{5}, \tfrac{3}{5}\big) \big\rceil + \big\lfloor 1, \big(\tfrac{1}{5}, \tfrac{4}{5}, \tfrac{1}{5}, \tfrac{1}{5}, \tfrac{3}{5}\big) \big\rceil$\tsep{1pt}\\
$\big\lfloor x_3^3x_5^2 + x_4^3x_5^2, (0, 0, 0, 0, 0) \big\rceil  $ & $\big\lfloor 1, \big(\tfrac{1}{5}, \tfrac{1}{5}, \tfrac{4}{5}, \tfrac{1}{5}, \tfrac{3}{5}\big) \big\rceil + \big\lfloor 1, \big(\tfrac{1}{5}, \tfrac{1}{5}, \tfrac{1}{5}, \tfrac{4}{5}, \tfrac{3}{5}\big) \big\rceil$\tsep{1pt}\\
$\big\lfloor x_1^2x_5^3 + x_2^2x_5^3, (0, 0, 0, 0, 0) \big\rceil  $ &$ \big\lfloor 1, \big(\tfrac{3}{5}, \tfrac{1}{5}, \tfrac{1}{5}, \tfrac{1}{5}, \tfrac{4}{5}\big) \rceil + \big\lfloor 1, \big(\tfrac{1}{5}, \tfrac{3}{5}, \tfrac{1}{5}, \tfrac{1}{5}, \tfrac{4}{5}\big) \big\rceil$\tsep{1pt}\\
$\big\lfloor x_2^2x_5^3 + x_4^2x_5^3, (0, 0, 0, 0, 0) \big\rceil  $ & $\big\lfloor 1, \big(\tfrac{1}{5}, \tfrac{3}{5}, \tfrac{1}{5}, \tfrac{1}{5}, \tfrac{4}{5}\big) \big\rceil + \big\lfloor 1, \big(\tfrac{1}{5}, \tfrac{1}{5}, \tfrac{1}{5}, \tfrac{3}{5}, \tfrac{4}{5}\big) \big\rceil$\tsep{1pt}\\
$\big\lfloor x_1^3x_2x_3 + x_1x_2^3x_4, (0, 0, 0, 0, 0) \big\rceil  $ & $\big\lfloor 1, \big(\tfrac{4}{5}, \tfrac{2}{5}, \tfrac{2}{5}, \tfrac{1}{5}, \tfrac{1}{5}\big) \big\rceil + \big\lfloor 1, \big(\tfrac{2}{5}, \tfrac{4}{5}, \tfrac{1}{5}, \tfrac{2}{5}, \tfrac{1}{5}\big) \big\rceil$\tsep{1pt}\\
$\big\lfloor x_1x_3^3x_4 + x_2x_3x_4^3, (0, 0, 0, 0, 0) \big\rceil  $ & $\big\lfloor 1, \big(\tfrac{2}{5}, \tfrac{1}{5}, \tfrac{4}{5}, \tfrac{2}{5}, \tfrac{1}{5}\big) \big\rceil + \big\lfloor 1, \big(\tfrac{1}{5}, \tfrac{2}{5}, \tfrac{2}{5}, \tfrac{4}{5}, \tfrac{1}{5}\big) \big\rceil$\tsep{1pt}\\
$\big\lfloor x_1x_2^3x_3 + x_1^3x_2x_4, (0, 0, 0, 0, 0) \big\rceil  $ & $\big\lfloor 1, \big(\tfrac{2}{5}, \tfrac{4}{5}, \tfrac{2}{5}, \tfrac{1}{5}, \tfrac{1}{5}\big) \big\rceil + \big\lfloor 1, \big(\tfrac{4}{5}, \tfrac{2}{5}, \tfrac{1}{5}, \tfrac{2}{5}, \tfrac{1}{5}\big) \big\rceil$\tsep{1pt}\\
$\big\lfloor x_1x_3x_4^3 + x_2x_3^3x_4, (0, 0, 0, 0, 0) \big\rceil  $ & $\big\lfloor 1, \big(\tfrac{2}{5}, \tfrac{1}{5}, \tfrac{2}{5}, \tfrac{4}{5}, \tfrac{1}{5}\big) \big\rceil + \big\lfloor 1, \big(\tfrac{1}{5}, \tfrac{2}{5}, \tfrac{4}{5}, \tfrac{2}{5}, \tfrac{1}{5}\big) \big\rceil$\tsep{1pt}\\
$\big\lfloor x_1x_2x_3^3 + x_1x_2x_4^3, (0, 0, 0, 0, 0) \big\rceil  $ & $\big\lfloor 1, \big(\tfrac{2}{5}, \tfrac{2}{5}, \tfrac{4}{5}, \tfrac{1}{5}, \tfrac{1}{5}\big) \big\rceil + \big\lfloor 1, \big(\tfrac{2}{5}, \tfrac{2}{5}, \tfrac{1}{5}, \tfrac{4}{5}, \tfrac{1}{5}\big) \big\rceil$\tsep{1pt}\\
$\big\lfloor x_1^3x_3x_4 + x_2^3x_3x_4, (0, 0, 0, 0, 0) \big\rceil  $ & $\big\lfloor 1, \big(\tfrac{4}{5}, \tfrac{1}{5}, \tfrac{2}{5}, \tfrac{2}{5}, \tfrac{1}{5}\big) \big\rceil + \big\lfloor 1, \big(\tfrac{1}{5}, \tfrac{4}{5}, \tfrac{2}{5}, \tfrac{2}{5}, \tfrac{1}{5}\big) \big\rceil$\tsep{1pt}\\
$\big\lfloor x_1x_2^3x_5 + x_1^3x_2x_5, (0, 0, 0, 0, 0) \big\rceil  $ & $\big\lfloor 1, \big(\tfrac{2}{5}, \tfrac{4}{5}, \tfrac{1}{5}, \tfrac{1}{5}, \tfrac{2}{5}\big) \rceil + \big\lfloor 1, \big(\tfrac{4}{5}, \tfrac{2}{5}, \tfrac{1}{5}, \tfrac{1}{5}, \tfrac{2}{5}\big) \big\rceil$\tsep{1pt}\\
$\big\lfloor x_3x_4^3x_5 + x_3^3x_4x_5, (0, 0, 0, 0, 0) \big\rceil   $ & $\big\lfloor 1, \big(\tfrac{1}{5}, \tfrac{1}{5}, \tfrac{2}{5}, \tfrac{4}{5}, \tfrac{2}{5}\big) \big\rceil + \big\lfloor 1, \big(\tfrac{1}{5}, \tfrac{1}{5}, \tfrac{4}{5}, \tfrac{2}{5}, \tfrac{2}{5}\big) \big\rceil$\tsep{1pt}\\
$\big\lfloor x_1x_3^3x_5 + x_2x_4^3x_5, (0, 0, 0, 0, 0) \big\rceil   $& $\big\lfloor 1, \big(\tfrac{2}{5}, \tfrac{1}{5}, \tfrac{4}{5}, \tfrac{1}{5}, \tfrac{2}{5}\big) \big\rceil + \big\lfloor 1, \big(\tfrac{1}{5}, \tfrac{2}{5}, \tfrac{1}{5}, \tfrac{4}{5}, \tfrac{2}{5}\big) \big\rceil$\tsep{1pt}\\
$\big\lfloor x_1^3x_3x_5 + x_2^3x_4x_5, (0, 0, 0, 0, 0) \big\rceil   $& $\big\lfloor 1, \big(\tfrac{4}{5}, \tfrac{1}{5}, \tfrac{2}{5}, \tfrac{1}{5}, \tfrac{2}{5}\big) \big\rceil + \big\lfloor 1, \big(\tfrac{1}{5}, \tfrac{4}{5}, \tfrac{1}{5}, \tfrac{2}{5}, \tfrac{2}{5}\big) \big\rceil$\tsep{1pt}\\
$\big\lfloor x_1x_4^3x_5 + x_2x_3^3x_5, (0, 0, 0, 0, 0) \big\rceil  $ & $\big\lfloor 1, \big(\tfrac{2}{5}, \tfrac{1}{5}, \tfrac{1}{5}, \tfrac{4}{5}, \tfrac{2}{5}\big) \big\rceil + \big\lfloor 1, \big(\tfrac{1}{5}, \tfrac{2}{5}, \tfrac{4}{5}, \tfrac{1}{5}, \tfrac{2}{5}\big) \big\rceil$\tsep{1pt}\\
$\big\lfloor x_1^3x_4x_5 + x_2^3x_3x_5, (0, 0, 0, 0, 0) \big\rceil   $& $\big\lfloor 1, \big(\tfrac{4}{5}, \tfrac{1}{5}, \tfrac{1}{5}, \tfrac{2}{5}, \tfrac{2}{5}\big) \big\rceil + \big\lfloor 1, \big(\tfrac{1}{5}, \tfrac{4}{5}, \tfrac{2}{5}, \tfrac{1}{5}, \tfrac{2}{5}\big) \big\rceil$\tsep{1pt}\\
$\big\lfloor x_1x_2x_5^3, (0, 0, 0, 0, 0) \big\rceil   $& $\big\lfloor 1, \big(\tfrac{2}{5}, \tfrac{2}{5}, \tfrac{1}{5}, \tfrac{1}{5}, \tfrac{4}{5}\big) \big\rceil$\tsep{1pt}\\
$\big\lfloor x_3x_4x_5^3, (0, 0, 0, 0, 0) \big\rceil   $& $\big\lfloor 1, \big(\tfrac{1}{5}, \tfrac{1}{5}, \tfrac{2}{5}, \tfrac{2}{5}, \tfrac{4}{5}\big) \big\rceil$\tsep{1pt}\\
$\big\lfloor x_1x_3x_5^3 + x_2x_4x_5^3, (0, 0, 0, 0, 0) \big\rceil   $& $\big\lfloor 1, \big(\tfrac{2}{5}, \tfrac{1}{5}, \tfrac{2}{5}, \tfrac{1}{5}, \tfrac{4}{5}\big) \big\rceil + \big\lfloor 1, \big(\tfrac{1}{5}, \tfrac{2}{5}, \tfrac{1}{5}, \tfrac{2}{5}, \tfrac{4}{5}\big) \big\rceil$\tsep{1pt}\\
$\big\lfloor x_1x_4x_5^3 + x_2x_3x_5^3, (0, 0, 0, 0, 0) \big\rceil   $& $\big\lfloor 1, \big(\tfrac{2}{5}, \tfrac{1}{5}, \tfrac{1}{5}, \tfrac{2}{5}, \tfrac{4}{5}\big) \big\rceil + \big\lfloor 1, \big(\tfrac{1}{5}, \tfrac{2}{5}, \tfrac{2}{5}, \tfrac{1}{5}, \tfrac{4}{5}\big) \big\rceil$\tsep{1pt}\\
$\big\lfloor x_1^2x_2^2x_3 + x_1^2x_2^2x_4, (0, 0, 0, 0, 0) \big\rceil   $& $\big\lfloor 1, \big(\tfrac{3}{5}, \tfrac{3}{5}, \tfrac{2}{5}, \tfrac{1}{5}, \tfrac{1}{5}\big) \big\rceil + \big\lfloor 1, \big(\tfrac{3}{5}, \tfrac{3}{5}, \tfrac{1}{5}, \tfrac{2}{5}, \tfrac{1}{5}\big) \big\rceil$\tsep{1pt}\\
$\big\lfloor x_1x_2^2x_3^2 + x_1^2x_2x_4^2, (0, 0, 0, 0, 0) \big\rceil   $& $\big\lfloor 1, \big(\tfrac{2}{5}, \tfrac{3}{5}, \tfrac{3}{5}, \tfrac{1}{5}, \tfrac{1}{5}\big) \big\rceil + \big\lfloor 1, \big(\tfrac{3}{5}, \tfrac{2}{5}, \tfrac{1}{5}, \tfrac{3}{5}, \tfrac{1}{5}\big) \big\rceil$\tsep{1pt}\\
$\big\lfloor x_1^2x_2x_3^2 + x_1x_2^2x_4^2, (0, 0, 0, 0, 0) \big\rceil   $& $\big\lfloor 1, \big(\tfrac{3}{5}, \tfrac{2}{5}, \tfrac{3}{5}, \tfrac{1}{5}, \tfrac{1}{5}\big) \big\rceil + \big\lfloor 1, \big(\tfrac{2}{5}, \tfrac{3}{5}, \tfrac{1}{5}, \tfrac{3}{5}, \tfrac{1}{5}\big) \big\rceil$\tsep{1pt}\\
$\big\lfloor x_1^2x_3^2x_4 + x_2^2x_3x_4^2, (0, 0, 0, 0, 0) \big\rceil   $& $\big\lfloor 1, \big(\tfrac{3}{5}, \tfrac{1}{5}, \tfrac{3}{5}, \tfrac{2}{5}, \tfrac{1}{5}\big) \big\rceil + \big\lfloor 1, \big(\tfrac{1}{5}, \tfrac{3}{5}, \tfrac{2}{5}, \tfrac{3}{5}, \tfrac{1}{5}\big) \big\rceil$\tsep{1pt}\\
$\big\lfloor x_1x_3^2x_4^2 + x_2x_3^2x_4^2, (0, 0, 0, 0, 0) \big\rceil   $& $\big\lfloor 1, \big(\tfrac{2}{5}, \tfrac{1}{5}, \tfrac{3}{5}, \tfrac{3}{5}, \tfrac{1}{5}\big) \big\rceil + \big\lfloor 1, \big(\tfrac{1}{5}, \tfrac{2}{5}, \tfrac{3}{5}, \tfrac{3}{5}, \tfrac{1}{5}\big) \big\rceil$\tsep{1pt}\\
$\big\lfloor x_1^2x_3x_4^2 + x_2^2x_3^2x_4, (0, 0, 0, 0, 0) \big\rceil   $& $\big\lfloor 1, \big(\tfrac{3}{5}, \tfrac{1}{5}, \tfrac{2}{5}, \tfrac{3}{5}, \tfrac{1}{5}\big) \big\rceil + \big\lfloor 1, \big(\tfrac{1}{5}, \tfrac{3}{5}, \tfrac{3}{5}, \tfrac{2}{5}, \tfrac{1}{5}\big) \big\rceil$\tsep{1pt}\\
$\big\lfloor x_1^2x_2^2x_5, (0, 0, 0, 0, 0) \big\rceil   $& $\big\lfloor 1, \big(\tfrac{3}{5}, \tfrac{3}{5}, \tfrac{1}{5}, \tfrac{1}{5}, \tfrac{2}{5}\big) \big\rceil$\tsep{1pt}\\
$\big\lfloor x_3^2x_4^2x_5, (0, 0, 0, 0, 0) \big\rceil   $& $\big\lfloor 1, \big(\tfrac{1}{5}, \tfrac{1}{5}, \tfrac{3}{5}, \tfrac{3}{5}, \tfrac{2}{5}\big) \big\rceil$\tsep{1pt}\\
$\big\lfloor x_1^2x_3^2x_5 + x_2^2x_4^2x_5, (0, 0, 0, 0, 0) \big\rceil   $& $\big\lfloor 1, \big(\tfrac{3}{5}, \tfrac{1}{5}, \tfrac{3}{5}, \tfrac{1}{5}, \tfrac{2}{5}\big) \rceil + \big\lfloor 1, \big(\tfrac{1}{5}, \tfrac{3}{5}, \tfrac{1}{5}, \tfrac{3}{5}, \tfrac{2}{5}\big) \big\rceil$\tsep{1pt}\\
$\big\lfloor x_1^2x_4^2x_5 + x_2^2x_3^2x_5, (0, 0, 0, 0, 0) \big\rceil   $& $\big\lfloor 1, \big(\tfrac{3}{5}, \tfrac{1}{5}, \tfrac{1}{5}, \tfrac{3}{5}, \tfrac{2}{5}\big) \big\rceil + \big\lfloor 1, \big(\tfrac{1}{5}, \tfrac{3}{5}, \tfrac{3}{5}, \tfrac{1}{5}, \tfrac{2}{5}\big) \big\rceil$\tsep{1pt}\\
$\big\lfloor x_1x_2^2x_5^2 + x_1^2x_2x_5^2, (0, 0, 0, 0, 0) \big\rceil   $& $\big\lfloor 1, \big(\tfrac{2}{5}, \tfrac{3}{5}, \tfrac{1}{5}, \tfrac{1}{5}, \tfrac{3}{5}\big) \big\rceil + \big\lfloor 1, \big(\tfrac{3}{5}, \tfrac{2}{5}, \tfrac{1}{5}, \tfrac{1}{5}, \tfrac{3}{5}\big) \big\rceil$\tsep{1pt}\\
$\big\lfloor x_3x_4^2x_5^2 + x_3^2x_4x_5^2, (0, 0, 0, 0, 0) \big\rceil   $& $\big\lfloor 1, \big(\tfrac{1}{5}, \tfrac{1}{5}, \tfrac{2}{5}, \tfrac{3}{5}, \tfrac{3}{5}\big) \big\rceil + \big\lfloor 1, \big(\tfrac{1}{5}, \tfrac{1}{5}, \tfrac{3}{5}, \tfrac{2}{5}, \tfrac{3}{5}\big) \big\rceil$\tsep{1pt}\\
$\big\lfloor x_1x_3^2x_5^2 + x_2x_4^2x_5^2, (0, 0, 0, 0, 0) \big\rceil   $& $\big\lfloor 1, \big(\tfrac{2}{5}, \tfrac{1}{5}, \tfrac{3}{5}, \tfrac{1}{5}, \tfrac{3}{5}\big) \big\rceil + \big\lfloor 1, \big(\tfrac{1}{5}, \tfrac{2}{5}, \tfrac{1}{5}, \tfrac{3}{5}, \tfrac{3}{5}\big) \big\rceil$\tsep{1pt}\\
$\big\lfloor x_1^2x_3x_5^2 + x_2^2x_4x_5^2, (0, 0, 0, 0, 0) \big\rceil   $& $\big\lfloor 1, \big(\tfrac{3}{5}, \tfrac{1}{5}, \tfrac{2}{5}, \tfrac{1}{5}, \tfrac{3}{5}\big) \big\rceil + \big\lfloor 1, \big(\tfrac{1}{5}, \tfrac{3}{5}, \tfrac{1}{5}, \tfrac{2}{5}, \tfrac{3}{5}\big) \big\rceil$\tsep{1pt}\\
$\big\lfloor x_1x_4^2x_5^2 + x_2x_3^2x_5^2, (0, 0, 0, 0, 0) \big\rceil   $& $\big\lfloor 1, \big(\tfrac{2}{5}, \tfrac{1}{5}, \tfrac{1}{5}, \tfrac{3}{5}, \tfrac{3}{5}\big) \big\rceil + \big\lfloor 1, \big(\tfrac{1}{5}, \tfrac{2}{5}, \tfrac{3}{5}, \tfrac{1}{5}, \tfrac{3}{5}\big) \big\rceil$\tsep{1pt}\\
$\big\lfloor x_1^2x_4x_5^2 + x_2^2x_3x_5^2, (0, 0, 0, 0, 0) \big\rceil   $& $\big\lfloor 1, \big(\tfrac{3}{5}, \tfrac{1}{5}, \tfrac{1}{5}, \tfrac{2}{5}, \tfrac{3}{5}\big) \big\rceil + \big\lfloor 1, \big(\tfrac{1}{5}, \tfrac{3}{5}, \tfrac{2}{5}, \tfrac{1}{5}, \tfrac{3}{5}\big) \big\rceil$\tsep{1pt}\\
$\big\lfloor x_1^2x_2x_3x_4 + x_1x_2^2x_3x_4, (0, 0, 0, 0, 0) \big\rceil   $& $\big\lfloor 1, \big(\tfrac{3}{5}, \tfrac{2}{5}, \tfrac{2}{5}, \tfrac{2}{5}, \tfrac{1}{5}\big) \big\rceil + \big\lfloor 1, \big(\tfrac{2}{5}, \tfrac{3}{5}, \tfrac{2}{5}, \tfrac{2}{5}, \tfrac{1}{5}\big) \big\rceil$\tsep{1pt}\\
$\big\lfloor x_1x_2x_3^2x_4 + x_1x_2x_3x_4^2, (0, 0, 0, 0, 0) \big\rceil   $& $\big\lfloor 1, \big(\tfrac{2}{5}, \tfrac{2}{5}, \tfrac{3}{5}, \tfrac{2}{5}, \tfrac{1}{5}\big) \big\rceil + \big\lfloor 1, \big(\tfrac{2}{5}, \tfrac{2}{5}, \tfrac{2}{5}, \tfrac{3}{5}, \tfrac{1}{5}\big) \big\rceil$\tsep{1pt}\\
$\big\lfloor x_1^2x_2x_3x_5 + x_1x_2^2x_4x_5, (0, 0, 0, 0, 0) \big\rceil   $& $\big\lfloor 1, \big(\tfrac{3}{5}, \tfrac{2}{5}, \tfrac{2}{5}, \tfrac{1}{5}, \tfrac{2}{5}\big) \big\rceil + \big\lfloor 1, \big(\tfrac{2}{5}, \tfrac{3}{5}, \tfrac{1}{5}, \tfrac{2}{5}, \tfrac{2}{5}\big) \big\rceil$\tsep{1pt}\\
$\big\lfloor x_1x_3^2x_4x_5 + x_2x_3x_4^2x_5, (0, 0, 0, 0, 0) \big\rceil   $& $\big\lfloor 1, \big(\tfrac{2}{5}, \tfrac{1}{5}, \tfrac{3}{5}, \tfrac{2}{5}, \tfrac{2}{5}\big) \big\rceil + \big\lfloor 1, \big(\tfrac{1}{5}, \tfrac{2}{5}, \tfrac{2}{5}, \tfrac{3}{5}, \tfrac{2}{5}\big) \big\rceil$\tsep{1pt}\\
$\big\lfloor x_1x_2^2x_3x_5 + x_1^2x_2x_4x_5, (0, 0, 0, 0, 0) \big\rceil   $& $\big\lfloor 1, \big(\tfrac{2}{5}, \tfrac{3}{5}, \tfrac{2}{5}, \tfrac{1}{5}, \tfrac{2}{5}\big) \big\rceil + \big\lfloor 1, \big(\tfrac{3}{5}, \tfrac{2}{5}, \tfrac{1}{5}, \tfrac{2}{5}, \tfrac{2}{5}\big) \big\rceil$\tsep{1pt}\\
$\big\lfloor x_1x_3x_4^2x_5 + x_2x_3^2x_4x_5, (0, 0, 0, 0, 0) \big\rceil   $& $\big\lfloor 1, \big(\tfrac{2}{5}, \tfrac{1}{5}, \tfrac{2}{5}, \tfrac{3}{5}, \tfrac{2}{5}\big) \big\rceil + \big\lfloor 1, \big(\tfrac{1}{5}, \tfrac{2}{5}, \tfrac{3}{5}, \tfrac{2}{5}, \tfrac{2}{5}\big) \big\rceil$\tsep{1pt}\\
$\big\lfloor x_1x_2x_3^2x_5 + x_1x_2x_4^2x_5, (0, 0, 0, 0, 0) \big\rceil   $& $\big\lfloor 1, \big(\tfrac{2}{5}, \tfrac{2}{5}, \tfrac{3}{5}, \tfrac{1}{5}, \tfrac{2}{5}\big) \big\rceil + \big\lfloor 1, \big(\tfrac{2}{5}, \tfrac{2}{5}, \tfrac{1}{5}, \tfrac{3}{5}, \tfrac{2}{5}\big) \big\rceil$\tsep{1pt}\\
$\big\lfloor x_1^2x_3x_4x_5 + x_2^2x_3x_4x_5, (0, 0, 0, 0, 0) \big\rceil   $& $\big\lfloor 1, \big(\tfrac{3}{5}, \tfrac{1}{5}, \tfrac{2}{5}, \tfrac{2}{5}, \tfrac{2}{5}\big) \big\rceil + \big\lfloor 1, \big(\tfrac{1}{5}, \tfrac{3}{5}, \tfrac{2}{5}, \tfrac{2}{5}, \tfrac{2}{5}\big) \big\rceil$\tsep{1pt}\\
$\big\lfloor x_1x_2x_3x_5^2 + x_1x_2x_4x_5^2, (0, 0, 0, 0, 0) \big\rceil   $& $\big\lfloor 1, \big(\tfrac{2}{5}, \tfrac{2}{5}, \tfrac{2}{5}, \tfrac{1}{5}, \tfrac{3}{5}\big) \big\rceil + \big\lfloor 1, \big(\tfrac{2}{5}, \tfrac{2}{5}, \tfrac{1}{5}, \tfrac{2}{5}, \tfrac{3}{5}\big) \big\rceil$\tsep{1pt}\\
$\big\lfloor x_1x_3x_4x_5^2 + x_2x_3x_4x_5^2, (0, 0, 0, 0, 0) \big\rceil   $& $\big\lfloor 1, \big(\tfrac{2}{5}, \tfrac{1}{5}, \tfrac{2}{5}, \tfrac{2}{5}, \tfrac{3}{5}\big) \big\rceil + \big\lfloor 1, \big(\tfrac{1}{5}, \tfrac{2}{5}, \tfrac{2}{5}, \tfrac{2}{5}, \tfrac{3}{5}\big) \big\rceil$\tsep{1pt}\\
\bottomrule
\caption{Elements of degree $(1,2)$.}\label{tab:degree11}
\end{longtable}
}
{ \small
\begin{longtable}{ll} \toprule
\multicolumn{2}{c}{Mirror map: basis elements of degree $(2,1)$} \\
{A-model} &{B-model}	\\
	\midrule
	\endfirsthead
	\toprule
	\multicolumn{2}{c}{Basis elements of degree $(2,1)$} \\
{A-model} &{B-model}	\\
	 \midrule
	\endhead
$\big\lfloor x_1^3x_2^3x_3^3x_4 + x_1^3x_2^3x_3x_4^3, (0, 0, 0, 0, 0) \big\rceil   $& $\big\lfloor 1, \big(\tfrac{4}{5}, \tfrac{4}{5}, \tfrac{4}{5}, \tfrac{2}{5}, \tfrac{1}{5}\big) \big\rceil + \big\lfloor 1, \big(\tfrac{4}{5}, \tfrac{4}{5}, \tfrac{2}{5}, \tfrac{4}{5}, \tfrac{1}{5}\big) \big\rceil$\tsep{1.2pt}\\
$\big\lfloor x_1^3x_2x_3^3x_4^3 + x_1x_2^3x_3^3x_4^3, (0, 0, 0, 0, 0) \big\rceil   $& $\big\lfloor 1, \big(\tfrac{4}{5}, \tfrac{2}{5}, \tfrac{4}{5}, \tfrac{4}{5}, \tfrac{1}{5}\big) \big\rceil + \big\lfloor 1, \big(\tfrac{2}{5}, \tfrac{4}{5}, \tfrac{4}{5}, \tfrac{4}{5}, \tfrac{1}{5}\big) \big\rceil$\tsep{1.2pt}\\
$\big\lfloor x_1^3x_2^3x_3^3x_5 + x_1^3x_2^3x_4^3x_5, (0, 0, 0, 0, 0) \big\rceil   $& $\big\lfloor 1, \big(\tfrac{4}{5}, \tfrac{4}{5}, \tfrac{4}{5}, \tfrac{1}{5}, \tfrac{2}{5}\big) \big\rceil + \big\lfloor 1, \big(\tfrac{4}{5}, \tfrac{4}{5}, \tfrac{1}{5}, \tfrac{4}{5}, \tfrac{2}{5}\big) \big\rceil$\tsep{1.2pt}\\
$\big\lfloor x_1^3x_2^3x_4^3x_5 + x_2^3x_3^3x_4^3x_5, (0, 0, 0, 0, 0) \big\rceil   $& $\big\lfloor 1, \big(\tfrac{4}{5}, \tfrac{4}{5}, \tfrac{1}{5}, \tfrac{4}{5}, \tfrac{2}{5}\big) \big\rceil + \big\lfloor 1, \big(\tfrac{1}{5}, \tfrac{4}{5}, \tfrac{4}{5}, \tfrac{4}{5}, \tfrac{2}{5}\big) \big\rceil$\tsep{1.2pt}\\
$\big\lfloor x_1^3x_2^3x_3x_5^3 + x_1^3x_2^3x_4x_5^3, (0, 0, 0, 0, 0) \big\rceil   $& $\big\lfloor 1, \big(\tfrac{4}{5}, \tfrac{4}{5}, \tfrac{2}{5}, \tfrac{1}{5}, \tfrac{4}{5}\big) \big\rceil + \big\lfloor 1, \big(\tfrac{4}{5}, \tfrac{4}{5}, \tfrac{1}{5}, \tfrac{2}{5}, \tfrac{4}{5}\big) \big\rceil$\tsep{1.2pt}\\
$\big\lfloor x_1x_3^3x_4^3x_5^3 + x_2x_3^3x_4^3x_5^3, (0, 0, 0, 0, 0) \big\rceil   $& $\big\lfloor 1, \big(\tfrac{2}{5}, \tfrac{1}{5}, \tfrac{4}{5}, \tfrac{4}{5}, \tfrac{4}{5}\big) \big\rceil + \big\lfloor 1, \big(\tfrac{1}{5}, \tfrac{2}{5}, \tfrac{4}{5}, \tfrac{4}{5}, \tfrac{4}{5}\big) \big\rceil$\tsep{1.2pt}\\
$\big\lfloor x_1^3x_2x_3^3x_5^3 + x_1x_2^3x_4^3x_5^3, (0, 0, 0, 0, 0) \big\rceil   $& $\big\lfloor 1, \big(\tfrac{4}{5}, \tfrac{2}{5}, \tfrac{4}{5}, \tfrac{1}{5}, \tfrac{4}{5}\big) \big\rceil + \big\lfloor 1, \big(\tfrac{2}{5}, \tfrac{4}{5}, \tfrac{1}{5}, \tfrac{4}{5}, \tfrac{4}{5}\big) \big\rceil$\tsep{1.2pt}\\
$\big\lfloor x_1^3x_3^3x_4x_5^3 + x_2^3x_3x_4^3x_5^3, (0, 0, 0, 0, 0) \big\rceil   $& $\big\lfloor 1, \big(\tfrac{4}{5}, \tfrac{1}{5}, \tfrac{4}{5}, \tfrac{2}{5}, \tfrac{4}{5}\big) \big\rceil + \big\lfloor 1, \big(\tfrac{1}{5}, \tfrac{4}{5}, \tfrac{2}{5}, \tfrac{4}{5}, \tfrac{4}{5}\big) \big\rceil$\tsep{1.2pt}\\
$\big\lfloor x_1x_2^3x_3^3x_5^3 + x_1^3x_2x_4^3x_5^3, (0, 0, 0, 0, 0) \big\rceil   $& $\big\lfloor 1, \big(\tfrac{2}{5}, \tfrac{4}{5}, \tfrac{4}{5}, \tfrac{1}{5}, \tfrac{4}{5}\big) \big\rceil + \big\lfloor 1, \big(\tfrac{4}{5}, \tfrac{2}{5}, \tfrac{1}{5}, \tfrac{4}{5}, \tfrac{4}{5}\big) \big\rceil$\tsep{1.2pt}\\
$\big\lfloor x_1^3x_3x_4^3x_5^3 + x_2^3x_3^3x_4x_5^3, (0, 0, 0, 0, 0) \big\rceil   $& $\big\lfloor 1, \big(\tfrac{4}{5}, \tfrac{1}{5}, \tfrac{2}{5}, \tfrac{4}{5}, \tfrac{4}{5}\big) \big\rceil + \big\lfloor 1, \big(\tfrac{1}{5}, \tfrac{4}{5}, \tfrac{4}{5}, \tfrac{2}{5}, \tfrac{4}{5}\big) \big\rceil$\tsep{1.2pt}\\
$\big\lfloor x_1^3x_2^3x_3^2x_4^2, (0, 0, 0, 0, 0) \big\rceil   $& $\big\lfloor 1, \big(\tfrac{4}{5}, \tfrac{4}{5}, \tfrac{3}{5}, \tfrac{3}{5}, \tfrac{1}{5}\big) \big\rceil$\tsep{1.2pt}\\
$\big\lfloor x_1^2x_2^2x_3^3x_4^3, (0, 0, 0, 0, 0) \big\rceil   $& $\big\lfloor 1, \big(\tfrac{3}{5}, \tfrac{3}{5}, \tfrac{4}{5}, \tfrac{4}{5}, \tfrac{1}{5}\big) \big\rceil$\tsep{1.2pt}\\
$\big\lfloor x_1^3x_2^2x_3^3x_4^2 + x_1^2x_2^3x_3^2x_4^3, (0, 0, 0, 0, 0) \big\rceil   $& $\big\lfloor 1, \big(\tfrac{4}{5}, \tfrac{3}{5}, \tfrac{4}{5}, \tfrac{3}{5}, \tfrac{1}{5}\big) \big\rceil + \big\lfloor 1, \big(\tfrac{3}{5}, \tfrac{4}{5}, \tfrac{3}{5}, \tfrac{4}{5}, \tfrac{1}{5}\big) \big\rceil$\tsep{1.2pt}\\
$\big\lfloor x_1^3x_2^2x_3^2x_4^3 + x_1^2x_2^3x_3^3x_4^2, (0, 0, 0, 0, 0) \big\rceil   $& $\big\lfloor 1, \big(\tfrac{4}{5}, \tfrac{3}{5}, \tfrac{3}{5}, \tfrac{4}{5}, \tfrac{1}{5}\big) \big\rceil + \big\lfloor 1, \big(\tfrac{3}{5}, \tfrac{4}{5}, \tfrac{4}{5}, \tfrac{3}{5}, \tfrac{1}{5}\big) \big\rceil$\tsep{1.2pt}\\
$\big\lfloor x_1^3x_2^3x_3^2x_5^2 + x_1^3x_2^3x_4^2x_5^2, (0, 0, 0, 0, 0) \big\rceil   $& $\big\lfloor 1, \big(\tfrac{4}{5}, \tfrac{3}{5}, \tfrac{3}{5}, \tfrac{4}{5}, \tfrac{1}{5}\big) \big\rceil + \big\lfloor 1, \big(\tfrac{4}{5}, \tfrac{4}{5}, \tfrac{1}{5}, \tfrac{3}{5}, \tfrac{3}{5}\big) \big\rceil$\tsep{1.2pt}\\
$\big\lfloor x_1^2x_3^3x_4^3x_5^2 + x_2^2x_3^3x_4^3x_5^2, (0, 0, 0, 0, 0) \big\rceil   $& $\big\lfloor 1, \big(\tfrac{3}{5}, \tfrac{1}{5}, \tfrac{4}{5}, \tfrac{4}{5}, \tfrac{3}{5}\big) \big\rceil + \big\lfloor 1, \big(\tfrac{1}{5}, \tfrac{3}{5}, \tfrac{4}{5}, \tfrac{4}{5}, \tfrac{3}{5}\big) \big\rceil$\tsep{1.2pt}\\
$\big\lfloor x_1^3x_2^2x_3^3x_5^2 + x_1^2x_2^3x_4^3x_5^2, (0, 0, 0, 0, 0) \big\rceil   $& $\big\lfloor 1, \big(\tfrac{4}{5}, \tfrac{3}{5}, \tfrac{4}{5}, \tfrac{1}{5}, \tfrac{3}{5}\big) \big\rceil + \big\lfloor 1, \big(\tfrac{3}{5}, \tfrac{4}{5}, \tfrac{1}{5}, \tfrac{4}{5}, \tfrac{3}{5}\big) \big\rceil$\tsep{1.2pt}\\
$\big\lfloor x_1^3x_3^3x_4^2x_5^2 + x_2^3x_3^2x_4^3x_5^2, (0, 0, 0, 0, 0) \big\rceil   $& $\big\lfloor 1, \big(\tfrac{4}{5}, \tfrac{1}{5}, \tfrac{4}{5}, \tfrac{3}{5}, \tfrac{3}{5}\big) \big\rceil + \big\lfloor 1, \big(\tfrac{1}{5}, \tfrac{4}{5}, \tfrac{3}{5}, \tfrac{4}{5}, \tfrac{3}{5}\big) \big\rceil$\tsep{1.2pt}\\
$\big\lfloor x_1^2x_2^3x_3^3x_5^2 + x_1^3x_2^2x_4^3x_5^2, (0, 0, 0, 0, 0) \big\rceil   $& $\big\lfloor 1, \big(\tfrac{3}{5}, \tfrac{4}{5}, \tfrac{4}{5}, \tfrac{1}{5}, \tfrac{3}{5}\big) \big\rceil + \big\lfloor 1, \big(\tfrac{4}{5}, \tfrac{3}{5}, \tfrac{1}{5}, \tfrac{4}{5}, \tfrac{3}{5}\big) \big\rceil$\tsep{1.2pt}\\
$\big\lfloor x_1^3x_3^2x_4^3x_5^2 + x_2^3x_3^3x_4^2x_5^2, (0, 0, 0, 0, 0) \big\rceil   $& $\big\lfloor 1, \big(\tfrac{4}{5}, \tfrac{1}{5}, \tfrac{3}{5}, \tfrac{4}{5}, \tfrac{3}{5}\big) \big\rceil + \big\lfloor 1, \big(\tfrac{1}{5}, \tfrac{4}{5}, \tfrac{4}{5}, \tfrac{3}{5}, \tfrac{3}{5}\big) \big\rceil$\tsep{1.2pt}\\
$\big\lfloor x_1^3x_2^2x_3^2x_5^3 + x_1^2x_2^3x_4^2x_5^3, (0, 0, 0, 0, 0) \big\rceil   $& $\big\lfloor 1, \big(\tfrac{4}{5}, \tfrac{3}{5}, \tfrac{3}{5}, \tfrac{1}{5}, \tfrac{4}{5}\big) \big\rceil + \big\lfloor 1, \big(\tfrac{3}{5}, \tfrac{4}{5}, \tfrac{1}{5}, \tfrac{3}{5}, \tfrac{4}{5}\big) \big\rceil$\tsep{1.2pt}\\
$\big\lfloor x_1^2x_3^3x_4^2x_5^3 + x_2^2x_3^2x_4^3x_5^3, (0, 0, 0, 0, 0) \big\rceil   $& $\big\lfloor 1, \big(\tfrac{3}{5}, \tfrac{1}{5}, \tfrac{4}{5}, \tfrac{3}{5}, \tfrac{4}{5}\big) \big\rceil + \big\lfloor 1, \big(\tfrac{1}{5}, \tfrac{3}{5}, \tfrac{3}{5}, \tfrac{4}{5}, \tfrac{4}{5}\big) \big\rceil$\tsep{1.2pt}\\
$\big\lfloor x_1^2x_2^3x_3^2x_5^3 + x_1^3x_2^2x_4^2x_5^3, (0, 0, 0, 0, 0) \big\rceil   $& $\big\lfloor 1, \big(\tfrac{3}{5}, \tfrac{4}{5}, \tfrac{3}{5}, \tfrac{1}{5}, \tfrac{4}{5}\big) \big\rceil + \big\lfloor 1, \big(\tfrac{4}{5}, \tfrac{3}{5}, \tfrac{1}{5}, \tfrac{3}{5}, \tfrac{4}{5}\big) \big\rceil$\tsep{1.2pt}\\
$\big\lfloor x_1^2x_3^2x_4^3x_5^3 + x_2^2x_3^3x_4^2x_5^3, (0, 0, 0, 0, 0) \big\rceil   $& $\big\lfloor 1, \big(\tfrac{3}{5}, \tfrac{1}{5}, \tfrac{3}{5}, \tfrac{4}{5}, \tfrac{4}{5}\big) \big\rceil + \big\lfloor 1, \big(\tfrac{1}{5}, \tfrac{3}{5}, \tfrac{4}{5}, \tfrac{3}{5}, \tfrac{4}{5}\big) \big\rceil$\tsep{1.2pt}\\
$\big\lfloor x_1^2x_2^2x_3^3x_5^3 + x_1^2x_2^2x_4^3x_5^3, (0, 0, 0, 0, 0) \big\rceil   $& $\big\lfloor 1, \big(\tfrac{3}{5}, \tfrac{3}{5}, \tfrac{4}{5}, \tfrac{1}{5}, \tfrac{4}{5}\big) \big\rceil + \big\lfloor 1, \big(\tfrac{3}{5}, \tfrac{3}{5}, \tfrac{1}{5}, \tfrac{4}{5}, \tfrac{4}{5}\big) \big\rceil$\tsep{1.2pt}\\
$\big\lfloor x_1^3x_3^2x_4^2x_5^3 + x_2^3x_3^2x_4^2x_5^3, (0, 0, 0, 0, 0) \big\rceil   $& $\big\lfloor 1, \big(\tfrac{4}{5}, \tfrac{1}{5}, \tfrac{3}{5}, \tfrac{3}{5}, \tfrac{4}{5}\big) \big\rceil + \big\lfloor 1, \big(\tfrac{1}{5}, \tfrac{4}{5}, \tfrac{3}{5}, \tfrac{3}{5}, \tfrac{4}{5}\big) \big\rceil$\tsep{1.2pt}\\
$\big\lfloor x_1^3x_2^2x_3^2x_4^2x_5 + x_1^2x_2^3x_3^2x_4^2x_5, (0, 0, 0, 0, 0) \big\rceil   $& $\big\lfloor 1, \big(\tfrac{4}{5}, \tfrac{3}{5}, \tfrac{3}{5}, \tfrac{3}{5}, \tfrac{2}{5}\big) \big\rceil + \big\lfloor 1, \big(\tfrac{3}{5}, \tfrac{4}{5}, \tfrac{3}{5}, \tfrac{3}{5}, \tfrac{2}{5}\big) \big\rceil$\tsep{1.2pt}\\
$\big\lfloor x_1^2x_2^2x_3^3x_4^2x_5 + x_1^2x_2^2x_3^2x_4^3x_5, (0, 0, 0, 0, 0) \big\rceil   $& $\big\lfloor 1, \big(\tfrac{3}{5}, \tfrac{3}{5}, \tfrac{4}{5}, \tfrac{3}{5}, \tfrac{2}{5}\big) \big\rceil + \big\lfloor 1, \big(\tfrac{3}{5}, \tfrac{3}{5}, \tfrac{3}{5}, \tfrac{4}{5}, \tfrac{2}{5}\big) \big\rceil$\tsep{1.2pt}\\
$\big\lfloor x_1^3x_2^2x_3^2x_4x_5^2 + x_1^2x_2^3x_3x_4^2x_5^2, (0, 0, 0, 0, 0) \big\rceil   $& $\big\lfloor 1, \big(\tfrac{4}{5}, \tfrac{3}{5}, \tfrac{3}{5}, \tfrac{2}{5}, \tfrac{3}{5}\big) \big\rceil + \big\lfloor 1, \big(\tfrac{3}{5}, \tfrac{4}{5}, \tfrac{2}{5}, \tfrac{3}{5}, \tfrac{3}{5}\big) \big\rceil$\tsep{1.2pt}\\
$\big\lfloor x_1^2x_2x_3^3x_4^2x_5^2 + x_1x_2^2x_3^2x_4^3x_5^2, (0, 0, 0, 0, 0) \big\rceil   $& $\big\lfloor 1, \big(\tfrac{3}{5}, \tfrac{2}{5}, \tfrac{4}{5}, \tfrac{3}{5}, \tfrac{3}{5}\big) \big\rceil + \big\lfloor 1, \big(\tfrac{2}{5}, \tfrac{3}{5}, \tfrac{3}{5}, \tfrac{4}{5}, \tfrac{3}{5}\big) \big\rceil$\tsep{1.2pt}\\
$\big\lfloor x_1^3x_2^2x_3x_4^2x_5^2 + x_1^2x_2^3x_3^2x_4x_5^2, (0, 0, 0, 0, 0) \big\rceil   $& $\big\lfloor 1, \big(\tfrac{4}{5}, \tfrac{3}{5}, \tfrac{2}{5}, \tfrac{3}{5}, \tfrac{3}{5}\big) \big\rceil + \big\lfloor 1, \big(\tfrac{3}{5}, \tfrac{4}{5}, \tfrac{3}{5}, \tfrac{2}{5}, \tfrac{3}{5}\big) \big\rceil$\tsep{1.2pt}\\
$\big\lfloor x_1x_2^2x_3^3x_4^2x_5^2 + x_1^2x_2x_3^2x_4^3x_5^2, (0, 0, 0, 0, 0) \big\rceil   $& $\big\lfloor 1, \big(\tfrac{2}{5}, \tfrac{3}{5}, \tfrac{4}{5}, \tfrac{3}{5}, \tfrac{3}{5}\big) \rceil + \big\lfloor 1, \big(\tfrac{3}{5}, \tfrac{2}{5}, \tfrac{3}{5}, \tfrac{4}{5}, \tfrac{3}{5}\big) \big\rceil$\tsep{1.2pt}\\
$\big\lfloor x_1^3x_2x_3^2x_4^2x_5^2 + x_1x_2^3x_3^2x_4^2x_5^2, (0, 0, 0, 0, 0) \big\rceil   $& $\big\lfloor 1, \big(\tfrac{4}{5}, \tfrac{2}{5}, \tfrac{3}{5}, \tfrac{3}{5}, \tfrac{3}{5}\big) \big\rceil + \big\lfloor 1, \big(\tfrac{2}{5}, \tfrac{4}{5}, \tfrac{3}{5}, \tfrac{3}{5}, \tfrac{3}{5}\big) \big\rceil$\tsep{1.2pt}\\
$\big\lfloor x_1^2x_2^2x_3^3x_4x_5^2 + x_1^2x_2^2x_3x_4^3x_5^2, (0, 0, 0, 0, 0) \big\rceil   $& $\big\lfloor 1, \big(\tfrac{3}{5}, \tfrac{3}{5}, \tfrac{4}{5}, \tfrac{2}{5}, \tfrac{3}{5}\big) \big\rceil + \big\lfloor 1, \big(\tfrac{3}{5}, \tfrac{3}{5}, \tfrac{2}{5}, \tfrac{4}{5}, \tfrac{3}{5}\big) \big\rceil$\tsep{1.2pt}\\
$\big\lfloor x_1^2x_2^2x_3^2x_4x_5^3 + x_1^2x_2^2x_3x_4^2x_5^3, (0, 0, 0, 0, 0) \big\rceil   $& $\big\lfloor 1, \big(\tfrac{3}{5}, \tfrac{3}{5}, \tfrac{3}{5}, \tfrac{2}{5}, \tfrac{4}{5}\big) \big\rceil + \big\lfloor 1, \big(\tfrac{3}{5}, \tfrac{3}{5}, \tfrac{2}{5}, \tfrac{3}{5}, \tfrac{4}{5}\big) \big\rceil$\tsep{1.2pt}\\
$\big\lfloor x_1^2x_2x_3^2x_4^2x_5^3 + x_1x_2^2x_3^2x_4^2x_5^3, (0, 0, 0, 0, 0) \big\rceil   $& $\big\lfloor 1, \big(\tfrac{3}{5}, \tfrac{2}{5}, \tfrac{3}{5}, \tfrac{3}{5}, \tfrac{4}{5}\big) \big\rceil + \big\lfloor 1, \big(\tfrac{2}{5}, \tfrac{3}{5}, \tfrac{3}{5}, \tfrac{3}{5}, \tfrac{4}{5}\big) \big\rceil$\tsep{1.2pt}\\
$\big\lfloor x_1^3x_2^3x_3^2x_4x_5 + x_1^3x_2^3x_3x_4^2x_5, (0, 0, 0, 0, 0) \big\rceil   $& $\big\lfloor 1, \big(\tfrac{4}{5}, \tfrac{4}{5}, \tfrac{3}{5}, \tfrac{2}{5}, \tfrac{2}{5}\big) \big\rceil + \big\lfloor 1, \big(\tfrac{4}{5}, \tfrac{4}{5}, \tfrac{2}{5}, \tfrac{3}{5}, \tfrac{2}{5}\big) \big\rceil$\tsep{1.2pt}\\
$\big\lfloor x_1^2x_2x_3^3x_4^3x_5 + x_1x_2^2x_3^3x_4^3x_5, (0, 0, 0, 0, 0) \big\rceil   $& $\big\lfloor 1, \big(\tfrac{3}{5}, \tfrac{2}{5}, \tfrac{4}{5}, \tfrac{4}{5}, \tfrac{2}{5}\big) \big\rceil + \big\lfloor 1, \big(\tfrac{2}{5}, \tfrac{3}{5}, \tfrac{4}{5}, \tfrac{4}{5}, \tfrac{2}{5}\big) \big\rceil$\tsep{1.2pt}\\
$\big\lfloor x_1^3x_2^2x_3^3x_4x_5 + x_1^2x_2^3x_3x_4^3x_5, (0, 0, 0, 0, 0) \big\rceil   $& $\big\lfloor 1, \big(\tfrac{4}{5}, \tfrac{3}{5}, \tfrac{4}{5}, \tfrac{2}{5}, \tfrac{2}{5}\big) \big\rceil + \big\lfloor 1, \big(\tfrac{3}{5}, \tfrac{4}{5}, \tfrac{2}{5}, \tfrac{4}{5}, \tfrac{2}{5}\big) \big\rceil$\tsep{1.2pt}\\
$\big\lfloor x_1^3x_2x_3^3x_4^2x_5 + x_1x_2^3x_3^2x_4^3x_5, (0, 0, 0, 0, 0) \big\rceil   $& $\big\lfloor 1, \big(\tfrac{4}{5}, \tfrac{2}{5}, \tfrac{4}{5}, \tfrac{3}{5}, \tfrac{2}{5}\big) \big\rceil + \big\lfloor 1, \big(\tfrac{2}{5}, \tfrac{4}{5}, \tfrac{3}{5}, \tfrac{4}{5}, \tfrac{2}{5}\big) \big\rceil$\tsep{1.2pt}\\
$\big\lfloor x_1^3x_2^2x_3x_4^3x_5 + x_1^2x_2^3x_3^3x_4x_5, (0, 0, 0, 0, 0) \big\rceil   $& $\big\lfloor 1, \big(\tfrac{4}{5}, \tfrac{3}{5}, \tfrac{2}{5}, \tfrac{4}{5}, \tfrac{2}{5}\big) \big\rceil + \big\lfloor 1, \big(\tfrac{3}{5}, \tfrac{4}{5}, \tfrac{4}{5}, \tfrac{2}{5}, \tfrac{2}{5}\big) \big\rceil$\tsep{1.2pt}\\
$\big\lfloor x_1x_2^3x_3^3x_4^2x_5 + x_1^3x_2x_3^2x_4^3x_5, (0, 0, 0, 0, 0) \big\rceil   $& $\big\lfloor 1, \big(\tfrac{2}{5}, \tfrac{4}{5}, \tfrac{4}{5}, \tfrac{3}{5}, \tfrac{2}{5}\big) \big\rceil + \big\lfloor 1, \big(\tfrac{4}{5}, \tfrac{2}{5}, \tfrac{3}{5}, \tfrac{4}{5}, \tfrac{2}{5}\big) \big\rceil$\tsep{1.2pt}\\
$\big\lfloor x_1^3x_2^3x_3x_4x_5^2, (0, 0, 0, 0, 0) \big\rceil   $& $\big\lfloor 1, \big(\tfrac{4}{5}, \tfrac{4}{5}, \tfrac{2}{5}, \tfrac{2}{5}, \tfrac{3}{5}\big) \rceil$\tsep{1.2pt}\\
$\big\lfloor x_1x_2x_3^3x_4^3x_5^2, (0, 0, 0, 0, 0) \big\rceil   $& $\big\lfloor 1, \big(\tfrac{2}{5}, \tfrac{2}{5}, \tfrac{4}{5}, \tfrac{4}{5}, \tfrac{3}{5}\big) \big\rceil$\tsep{1.2pt}\\
$\big\lfloor x_1^3x_2x_3^3x_4x_5^2 + x_1x_2^3x_3x_4^3x_5^2, (0, 0, 0, 0, 0) \big\rceil   $& $\big\lfloor 1, \big(\tfrac{4}{5}, \tfrac{2}{5}, \tfrac{4}{5}, \tfrac{2}{5}, \tfrac{3}{5}\big) \big\rceil + \big\lfloor 1, \big(\tfrac{2}{5}, \tfrac{4}{5}, \tfrac{2}{5}, \tfrac{4}{5}, \tfrac{3}{5}\big) \big\rceil$\tsep{1.2pt}\\
$\big\lfloor x_1^3x_2x_3x_4^3x_5^2 + x_1x_2^3x_3^3x_4x_5^2, (0, 0, 0, 0, 0) \big\rceil   $& $\big\lfloor 1, \big(\tfrac{4}{5}, \tfrac{2}{5}, \tfrac{2}{5}, \tfrac{4}{5}, \tfrac{3}{5}\big) \big\rceil + \big\lfloor 1, \big(\tfrac{2}{5}, \tfrac{4}{5}, \tfrac{4}{5}, \tfrac{2}{5}, \tfrac{3}{5}\big) \big\rceil$\tsep{1.2pt}\\
$\big\lfloor x_1^3x_2^2x_3x_4x_5^3 + x_1^2x_2^3x_3x_4x_5^3, (0, 0, 0, 0, 0) \big\rceil   $& $\big\lfloor 1, \big(\tfrac{4}{5}, \tfrac{3}{5}, \tfrac{2}{5}, \tfrac{2}{5}, \tfrac{4}{5}\big) \big\rceil + \big\lfloor 1, \big(\tfrac{3}{5}, \tfrac{4}{5}, \tfrac{2}{5}, \tfrac{2}{5}, \tfrac{4}{5}\big) \big\rceil$\tsep{1.2pt}\\
$\big\lfloor x_1x_2x_3^3x_4^2x_5^3 + x_1x_2x_3^2x_4^3x_5^3, (0, 0, 0, 0, 0) \big\rceil   $& $\big\lfloor 1, \big(\tfrac{2}{5}, \tfrac{2}{5}, \tfrac{4}{5}, \tfrac{3}{5}, \tfrac{4}{5}\big) \big\rceil + \big\lfloor 1, \big(\tfrac{2}{5}, \tfrac{2}{5}, \tfrac{3}{5}, \tfrac{4}{5}, \tfrac{4}{5}\big) \big\rceil$\tsep{1.2pt}\\
$\big\lfloor x_1^3x_2x_3^2x_4x_5^3 + x_1x_2^3x_3x_4^2x_5^3, (0, 0, 0, 0, 0) \big\rceil   $& $\big\lfloor 1, \big(\tfrac{4}{5}, \tfrac{2}{5}, \tfrac{3}{5}, \tfrac{2}{5}, \tfrac{4}{5}\big) \big\rceil + \big\lfloor 1, \big(\tfrac{2}{5}, \tfrac{4}{5}, \tfrac{2}{5}, \tfrac{3}{5}, \tfrac{4}{5}\big) \big\rceil$\tsep{1.2pt}\\
$\big\lfloor x_1^2x_2x_3^3x_4x_5^3 + x_1x_2^2x_3x_4^3x_5^3, (0, 0, 0, 0, 0) \big\rceil   $& $\big\lfloor 1, \big(\tfrac{3}{5}, \tfrac{2}{5}, \tfrac{4}{5}, \tfrac{2}{5}, \tfrac{4}{5}\big) \big\rceil + \big\lfloor 1, \big(\tfrac{2}{5}, \tfrac{3}{5}, \tfrac{2}{5}, \tfrac{4}{5}, \tfrac{4}{5}\big) \big\rceil$\tsep{1.2pt}\\
$\big\lfloor x_1^3x_2x_3x_4^2x_5^3 + x_1x_2^3x_3^2x_4x_5^3, (0, 0, 0, 0, 0) \big\rceil   $& $\big\lfloor 1, \big(\tfrac{4}{5}, \tfrac{2}{5}, \tfrac{2}{5}, \tfrac{3}{5}, \tfrac{4}{5}\big) \big\rceil + \big\lfloor 1, \big(\tfrac{2}{5}, \tfrac{4}{5}, \tfrac{3}{5}, \tfrac{2}{5}, \tfrac{4}{5}\big) \big\rceil$\tsep{1.2pt}\\
$\big\lfloor x_1x_2^2x_3^3x_4x_5^3 + x_1^2x_2x_3x_4^3x_5^3, (0, 0, 0, 0, 0) \big\rceil   $& $\big\lfloor 1, \big(\tfrac{2}{5}, \tfrac{3}{5}, \tfrac{4}{5}, \tfrac{2}{5}, \tfrac{4}{5}\big) \big\rceil + \big\lfloor 1, \big(\tfrac{3}{5}, \tfrac{2}{5}, \tfrac{2}{5}, \tfrac{4}{5}, \tfrac{4}{5}\big) \big\rceil$ \tsep{1.2pt}\\
\bottomrule
\caption{Basis elements of degree $(2,1)$.}\label{tab:degree21}
\end{longtable}} 

The following basis elements in Table~\ref{tab:example2Rest} are not described by Theorem~\ref{thm:mirror}, yet we are still able to find the same number of each bidegree on either side. Notice the narrow sectors on the A-side correspond with broad sectors on the B-side and vice versa.

{\small	
\begin{longtable}{llll} \toprule
	\multicolumn{4}{c}{Mirror map: remaining basis elements} \\
{bidegree} &{A-model} &{B-model}	\\
	\midrule
	\endfirsthead
	\toprule
	\multicolumn{4}{c}{Basis elements of degree $(1,2)$} \\
{bidegree} &{A-model} &{B-model}	\\
	\midrule
	\endhead
$(1, 1) $& $\lfloor 1, ((12)(34))j_W \rceil $&$ \lfloor (x_1+x_2)(x_3+x_4), (12)(34) \rceil+ (24 \text{ others})$\tsep{1.2pt}\\
$(2, 2) $& $\big\lfloor 1, ((12)(34))(j_W)^2 \big\rceil $&$ \big\lfloor (x_1+x_2)^2(x_3+x_4)^2x_5^3, (12)(34) \big\rceil+ (24 \text{ others})$\tsep{1.2pt}\\
$(1, 1) $& $\big\lfloor 1, ((12)(34))(j_W)^3 \big\rceil $&$ \big\lfloor x_5^2, (12)(34) \big\rceil+ (24 \text{ others})$\tsep{1.2pt}\\
$(2, 2) $& $\big\lfloor 1, ((12)(34))(j_W)^4 \big\rceil $&$ \big\lfloor (x_1+x_2)^3(x_3+x_4)^3x_5, (12)(34) \big\rceil+ (24 \text{ others})$\tsep{1.2pt}\\
$(1, 2) $& $\big\lfloor (x_1+x_2)^2, (12)(34) \big\rceil $&$ \lfloor 1, ((12)(34))j_W \rceil + (24 \text{ others})$\tsep{1.2pt}\\
$(1, 2) $& $\big\lfloor (x_3+x_4)^2, (12)(34) \big\rceil $&$ \big\lfloor 1, ((12)(34))j_W^3 \big\rceil + (24 \text{ others})$\tsep{1.2pt}\\
$(1, 2) $& $\lfloor (x_1 + x_2)(x_3 + x_4), (12)(34) \rceil $& $\big\lfloor 1, ((12)(34))j_W^2K \big\rceil + (24 \text{ others})$\tsep{1.2pt}\\
$(1, 2) $& $\big\lfloor x_5^2, (12)(34) \big\rceil $& $\big\lfloor 1, ((12)(34))j_W^2L \big\rceil + (24 \text{ others})$\tsep{1.2pt}\\
$(1, 2) $& $\lfloor (x_1+x_2)x_5, (12)(34) \rceil $&$ \big\lfloor 1, ((12)(34))j_W^2K^4 \big\rceil + (24 \text{ others})$\tsep{1.2pt}\\
$(1, 2) $& $\lfloor (x_3+x_4)x_5, (12)(34) \rceil $&$ \big\lfloor 1, ((12)(34))j_W^2L^4 \big\rceil + (24 \text{ others})$\tsep{1.2pt}\\
$(2, 1) $& $\big\lfloor (x_1 + x_2)^3(x_3 + x_4)^3x_5, (12)(34) \big\rceil $& $\big\lfloor 1, ((12)(34))j_W^2 \big\rceil + (24 \text{ others})$\\
$(2, 1) $& $\big\lfloor (x_1 + x_2)^3(x_3+x_4)^2x_5^2, (12)(34) \big\rceil $& $\big\lfloor 1, ((12)(34))j_W^4 \big\rceil + (24 \text{ others})$\tsep{1.2pt}\\
$(2, 1) $& $\big\lfloor (x_1 + x_2)^2(x_3 + x_4)^3x_5^2, (12)(34) \big\rceil  $& $\big\lfloor 1, ((12)(34))j_W^3K \big\rceil + (24 \text{ others})$\tsep{1.2pt}\\
$(2, 1) $& $\big\lfloor (x_1 + x_2)^3(x_3 + x_4)x_5^3, (12)(34) \big\rceil $& $\big\lfloor 1, ((12)(34))j_W^3L \big\rceil + (24 \text{ others})$\tsep{1.2pt}\\
$(2, 1) $& $\big\lfloor (x_1 + x_2)(x_3 + x_4)^3x_5^3, (12)(34) \big\rceil $& $\big\lfloor 1, ((12)(34))j_W^3K^4 \big\rceil + (24 \text{ others})$\tsep{1.2pt}\\
$(2, 1) $& $\big\lfloor (x_1 + x_2)^2(x_3 + x_4)^2x_5^3, (12)(34) \big\rceil $& $\big\lfloor 1, ((12)(34))j_W^3L^4 \big\rceil + (24 \text{ others})$\tsep{1.2pt}\\
\bottomrule
\caption{Remaining basis elements.}\label{tab:example2Rest}
\end{longtable}} 
\end{Example}

\subsection{Bad example}

Unfortunately, the correspondence previously shown does not hold for all pairs $(W,G)$. In Example~\ref{eg:bad} we will see an example where the A- and B-model state spaces are isomorphic as vector spaces, however, the given bidegrees may not match. Ebeling and Gusein-Zade describe
a condition in \cite{Ebel2} on subgroups $K \leq G$ of pure permutations, which we will now describe. They conjecture that this condition is necessary for the Milnor fibers associated to $(W,G)$ and $\big(W^T,G^\star\big)$ to have the same orbifold Euler characteristic. It appears that this condition must also hold for $\mathcal{A}_{W, G}$ and $\mathcal{B}_{W^T, G^\star}$ to be isomorphic as bigraded vector spaces.

\begin{Definition}[Ebeling--Gusein-Zade \cite{Ebel2}]
Let $K$ be the subgroup of pure permutations in a~group $G \leq G_W^{\max}$. We say that $K$ satisfies the \emph{parity condition} (PC) if for each subgroup $T \leq K$ one has
\[
\dim \big(\mathbb{C}^N\big)^T \equiv N \pmod 2,
\]
where $\big(\mathbb{C}^N\big)^T = \big\{x \in \mathbb{C}^N\colon \sigma x = x \text{ for all } \sigma \in T\big\}.$
\end{Definition}

\begin{Example}
Consider $K = \langle (12)(34) \rangle$ from Example~\ref{eg:good}. Then $\big(\mathbb{C}^5\big)^K$ has dimension~$3$ and~$\big(\mathbb{C}^5\big)^{\{(1)\}}$ has dimension~$5$, which are both equal to~$5$ (mod~$2$). Thus $K$ satisfies the parity condition, and we have seen in Example~\ref{eg:good} that the bigraded state spaces of $\mathcal{A}_{W, G}$ and $\mathcal{B}_{W^T, G^\star}$ are isomorphic.
\end{Example}

\begin{Example}\label{eg:bad}
In this example, we will consider the Klein~4 group as our group of permutations. As we will see, this group does not satisfy PC. We will examine exactly where the mirror map fails to preserve the bigrading.

Let $W = x_1^5 + x_2^5 + x_3^5 + x_4^5 + x_5^5$ and $G = \langle j_W, (12)(34), (13)(24) \rangle$, where
\[
j_W = \big(\tfrac{1}{5}, \tfrac{1}{5}, \tfrac{1}{5}, \tfrac{1}{5}, \tfrac{1}{5}\big),
(12)(34) =
\begin{pmatrix}
 0 & 1 & 0 & 0 & 0\\
 1 & 0 & 0 & 0 & 0\\
 0 & 0 & 0 & 1 & 0\\
 0 & 0 & 1 & 0 & 0\\
 0 & 0 & 0 & 0 & 1
\end{pmatrix}, \qquad 
(13)(24) =
\begin{pmatrix}
 0 & 0 & 1 & 0 & 0\\
 0 & 0 & 0 & 1 & 0\\
 1 & 0 & 0 & 0 & 0\\
 0 & 1 & 0 & 0 & 0\\
 0 & 0 & 0 & 0 & 1
\end{pmatrix}.
\]
Then $W^T = W$ and the non-abelian dual group of $G$ is
\[
G^\star = \SL_W^{\rm diag}\cdot \langle (12)(34), (13)(24) \rangle,
\]
where
\[
\SL_W^{\rm diag} = \big\langle j_W, \big(0, \tfrac{2}{5}, \tfrac{1}{5}, \tfrac{1}{5}, \tfrac{1}{5}\big), \big(0, \tfrac{1}{5}, \tfrac{2}{5}, \tfrac{1}{5}, \tfrac{1}{5}\big), \big(0, \tfrac{1}{5}, \tfrac{1}{5}, \tfrac{2}{5}, \tfrac{1}{5}\big) \big\rangle.\]
This choice of group does not satisfy the PC above. Indeed, if $K = \langle (12)(34), (13)(24) \rangle$, then $(\mathbb{C}^5)^K = 2 \neq 5 \pmod 2$.

While the mirror map works for the sectors described in Theorem~\ref{thm:mirror}, it is not an isomorphism when considering the entire A- and B-models. We will show this explicitly by computing the basis elements of $\mathcal{A}_{W, G}$ and $\mathcal{B}_{W^T, G^\star}$, and then computing their bidegree. However, as proved in Theorem~\ref{thm:mirror}, the restricted mirror map is still an isomorphism, which we will list first.

The eight elements listed below in Table~\ref{tab:badexamplefirstelements} are the exact same as those from Example~\ref{eg:good}.
\begin{longtable}{lll} \toprule
\multicolumn{3}{c}{Mirror map: first basis elements} \\
{bidegree}& {A-model} &{B-model} \\
\midrule
\endfirsthead
\toprule
{bidegree}& {A-model} &{B-model} \\
\midrule
\endhead
$(0, 0)$  & $\lfloor 1, j_W \rceil $ & $\lfloor 1, (0, 0, 0, 0, 0) \rceil$\tsep{1pt}\\
$(1, 1)$  & $\big\lfloor 1, (j_W)^2 \big\rceil$ & $\lfloor x_1x_2x_3x_4, (0, 0, 0, 0, 0) \rceil$\tsep{1pt}\\
$(2, 2)$ & $\big\lfloor 1, (j_W)^3 \big\rceil$ & $\big\lfloor x_1^2x_2^2x_3^2x_4^2, (0, 0, 0, 0, 0) \big\rceil$\tsep{1pt}\\
$(3, 3)$ & $\big\lfloor 1, (j_W)^4 \big\rceil$ & $\big\lfloor x_1^3x_2^3x_3^3x_4^3, (0, 0, 0, 0, 0) \big\rceil$\tsep{1pt}\\
$(0, 3)$  & $\lfloor 1, (0, 0, 0, 0, 0) \rceil$ & $\lfloor 1, j_W \rceil$\tsep{1pt}\\
$(1, 2)$ & $\lfloor x_1x_2x_3x_4, (0, 0, 0, 0, 0) \rceil$ & $\big\lfloor 1, (j_W)^2 \big\rceil$\tsep{1pt}\\
$(2, 1)$& $\big\lfloor x_1^2x_2^2x_3^2x_4^2, (0, 0, 0, 0, 0) \big\rceil$ & $\big\lfloor 1, (j_W)^3 \big\rceil$\tsep{1pt}\\
$(3, 0)$& $\big\lfloor x_1^3x_2^3x_3^3x_4^3, (0, 0, 0, 0, 0) \big\rceil$ & $\big\lfloor 1, (j_W)^4 \big\rceil$\tsep{1pt}\\
\bottomrule
\caption{In our bad example, these are the first elements that match up in the mirror map.}\label{tab:badexamplefirstelements}
\end{longtable}

The following $28$ corresponding elements in Table~\ref{tab:baddegree12} have a bidegree of $(1, 2)$. On the A side these come from untwisted sectors, and on the B side these are from the narrow sectors, again following the recipe from Theorem~\ref{thm:mirror}.

{\small	
\begin{longtable}{@{\,}ll@{\,}} \toprule
\multicolumn{2}{c}{Basis elements of degree $(1,2)$} \\
{A-model} &{B-model}	\\
	\midrule
	\endfirsthead
	\toprule
{A-model} &{B-model}	\\
	 \midrule
	\endhead
$\big\lfloor x_1^3x_2^2 + x_1^2x_2^3 + x_3^3x_4^2 + x_3^2x_4^3, (0, 0, 0, 0, 0) \big\rceil   $& $\big\lfloor 1, \big(\tfrac{4}{5}, \tfrac{3}{5}, \tfrac{1}{5}, \tfrac{1}{5}, \tfrac{1}{5}\big) \big\rceil + (3 \text{ others})$\tsep{1pt}\\
$\big\lfloor x_1^3x_3^2 + x_1^2x_3^3 + x_2^3x_4^2 + x_2^2x_4^3, (0, 0, 0, 0, 0) \big\rceil   $& $\big\lfloor 1, \big(\tfrac{4}{5}, \tfrac{1}{5}, \tfrac{3}{5}, \tfrac{1}{5}, \tfrac{1}{5}\big) \big\rceil + (3 \text{ others})$\tsep{1pt}\\
$\big\lfloor x_1^3x_4^2 + x_1^2x_4^3 + x_2^3x_3^2 + x_2^2x_3^3, (0, 0, 0, 0, 0) \big\rceil   $& $\big\lfloor 1, \big(\tfrac{4}{5}, \tfrac{1}{5}, \tfrac{1}{5}, \tfrac{3}{5}, \tfrac{1}{5}\big) \big\rceil + (3 \text{ others})$\tsep{1pt}\\
$\big\lfloor x_1^3x_5^2 + x_2^3x_5^2 + x_3^3x_5^2 + x_4^3x_5^2, (0, 0, 0, 0, 0) \big\rceil   $& $\big\lfloor 1, \big(\tfrac{4}{5}, \tfrac{1}{5}, \tfrac{1}{5}, \tfrac{1}{5}, \tfrac{3}{5}\big) \big\rceil + (3 \text{ others})$\tsep{1pt}\\
$\big\lfloor x_1^2x_5^3 + x_2^2x_5^3 + x_2^3x_5^3 + x_4^2x_5^3, (0, 0, 0, 0, 0) \big\rceil   $& $\big\lfloor 1, \big(\tfrac{3}{5}, \tfrac{1}{5}, \tfrac{1}{5}, \tfrac{1}{5}, \tfrac{4}{5}\big) \big\rceil + (3 \text{ others})$\tsep{1pt}\\
$\big\lfloor x_1^3x_2x_3 + x_1x_2^3x_4 + x_1x_3^3x_4 + x_2x_3x_4^3, (0, 0, 0, 0, 0) \big\rceil   $& $\big\lfloor 1, \big(\tfrac{4}{5}, \tfrac{2}{5}, \tfrac{2}{5}, \tfrac{1}{5}, \tfrac{1}{5}\big) \big\rceil + (3 \text{ others})$\tsep{1pt}\\
$\big\lfloor x_1x_2^3x_3 + x_1^3x_2x_4 + x_1x_3x_4^3 + x_2x_3^3x_4, (0, 0, 0, 0, 0) \big\rceil   $& $\big\lfloor 1, \big(\tfrac{2}{5}, \tfrac{4}{5}, \tfrac{2}{5}, \tfrac{1}{5}, \tfrac{1}{5}\big) \big\rceil + (3 \text{ others})$\tsep{1pt}\\
$\big\lfloor x_1x_2x_3^3 + x_1x_2x_4^3 + x_1^3x_3x_4 + x_2^3x_3x_4, (0, 0, 0, 0, 0) \big\rceil   $& $\big\lfloor 1, \big(\tfrac{2}{5}, \tfrac{2}{5}, \tfrac{4}{5}, \tfrac{1}{5}, \tfrac{1}{5}\big) \big\rceil + (3 \text{ others})$\tsep{1pt}\\
$\big\lfloor x_1x_2^3x_5 + x_1^3x_2x_5 + x_3x_4^3x_5 + x_3^3x_4x_5, (0, 0, 0, 0, 0) \big\rceil   $& $\big\lfloor 1, \big(\tfrac{2}{5}, \tfrac{4}{5}, \tfrac{1}{5}, \tfrac{1}{5}, \tfrac{2}{5}\big) \big\rceil + (3 \text{ others})$\tsep{1pt}\\
$\big\lfloor x_1x_3^3x_5 + x_1^3x_3x_5 + x_2x_4^3x_5 + x_2^3x_4x_5, (0, 0, 0, 0, 0) \big\rceil   $& $\big\lfloor 1, \big(\tfrac{2}{5}, \tfrac{1}{5}, \tfrac{4}{5}, \tfrac{1}{5}, \tfrac{2}{5}\big) \big\rceil + (3 \text{ others})$\tsep{1pt}\\
$\big\lfloor x_1x_4^3x_5 + x_1^3x_4x_5 + x_2x_3^3x_5 + x_2^3x_3x_5, (0, 0, 0, 0, 0) \big\rceil   $& $\big\lfloor 1, \big(\tfrac{2}{5}, \tfrac{1}{5}, \tfrac{1}{5}, \tfrac{4}{5}, \tfrac{2}{5}\big) \big\rceil + (3 \text{ others})$\tsep{1pt}\\
$\big\lfloor x_1x_2x_5^3 + x_3x_4x_5^3, (0, 0, 0, 0, 0) \big\rceil   $& $\big\lfloor 1, \big(\tfrac{2}{5}, \tfrac{2}{5}, \tfrac{1}{5}, \tfrac{1}{5}, \tfrac{4}{5}\big) \big\rceil + \big\lfloor 1, \big(\tfrac{1}{5}, \tfrac{1}{5}, \tfrac{2}{5}, \tfrac{2}{5}, \tfrac{4}{5}\big) \big\rceil$\tsep{1pt}\\
$\big\lfloor x_1x_3x_5^3 + x_2x_4x_5^3, (0, 0, 0, 0, 0) \big\rceil   $& $\big\lfloor 1, \big(\tfrac{2}{5}, \tfrac{1}{5}, \tfrac{2}{5}, \tfrac{1}{5}, \tfrac{4}{5}\big) \big\rceil + \big\lfloor 1, \big(\tfrac{1}{5}, \tfrac{2}{5}, \tfrac{1}{5}, \tfrac{2}{5}, \tfrac{4}{5}\big) \big\rceil$\tsep{1pt}\\
$\big\lfloor x_1x_4x_5^3 + x_2x_3x_5^3, (0, 0, 0, 0, 0) \big\rceil   $& $\big\lfloor 1, \big(\tfrac{2}{5}, \tfrac{1}{5}, \tfrac{1}{5}, \tfrac{2}{5}, \tfrac{4}{5}\big) \big\rceil + \big\lfloor 1, \big(\tfrac{1}{5}, \tfrac{2}{5}, \tfrac{2}{5}, \tfrac{1}{5}, \tfrac{4}{5}\big) \big\rceil$\tsep{1pt}\\
$\big\lfloor x_1x_2^2x_3^2 + x_1^2x_2x_4^2 + x_1^2x_3x_4^2 + x_2^2x_3^2x_4, (0, 0, 0, 0, 0) \big\rceil   $& $\big\lfloor 1, \big(\tfrac{2}{5}, \tfrac{3}{5}, \tfrac{3}{5}, \tfrac{1}{5}, \tfrac{1}{5}\big) \big\rceil + (3 \text{ others})$\tsep{1pt}\\
$\big\lfloor x_1^2x_2x_3^2 + x_1x_2^2x_4^2 + x_1^2x_3^2x_4 + x_2^2x_3x_4^2, (0, 0, 0, 0, 0) \big\rceil   $& $\big\lfloor 1, \big(\tfrac{3}{5}, \tfrac{2}{5}, \tfrac{3}{5}, \tfrac{1}{5}, \tfrac{1}{5}\big) \big\rceil + (3 \text{ others})$\tsep{1pt}\\
$\big\lfloor x_1^2x_2^2x_3 + x_1^2x_2^2x_4 + x_1x_3^2x_4^2 + x_2x_3^2x_4^2, (0, 0, 0, 0, 0) \big\rceil   $& $\big\lfloor 1, \big(\tfrac{3}{5}, \tfrac{3}{5}, \tfrac{2}{5}, \tfrac{1}{5}, \tfrac{1}{5}\big) \big\rceil + (3 \text{ others})$\tsep{1pt}\\
$\big\lfloor x_1^2x_2^2x_5 + x_3^2x_4^2x_5, (0, 0, 0, 0, 0) \big\rceil   $& $\big\lfloor 1, \big(\tfrac{3}{5}, \tfrac{3}{5}, \tfrac{1}{5}, \tfrac{1}{5}, \tfrac{2}{5}\big) \big\rceil + \big\lfloor 1, \big(\tfrac{1}{5}, \tfrac{1}{5}, \tfrac{3}{5}, \tfrac{3}{5}, \tfrac{2}{5}\big) \big\rceil$\tsep{1pt}\\
$\big\lfloor x_1^2x_3^2x_5 + x_2^2x_4^2x_5, (0, 0, 0, 0, 0) \big\rceil   $& $\big\lfloor 1, \big(\tfrac{3}{5}, \tfrac{1}{5}, \tfrac{3}{5}, \tfrac{1}{5}, \tfrac{2}{5}\big) \big\rceil + \big\lfloor 1, \big(\tfrac{1}{5}, \tfrac{3}{5}, \tfrac{1}{5}, \tfrac{3}{5}, \tfrac{2}{5}\big) \big\rceil$\tsep{1pt}\\
$\big\lfloor x_1^2x_4^2x_5 + x_2^2x_3^2x_5, (0, 0, 0, 0, 0) \big\rceil   $& $\big\lfloor 1, \big(\tfrac{3}{5}, \tfrac{1}{5}, \tfrac{1}{5}, \tfrac{3}{5}, \tfrac{2}{5}\big) \big\rceil + \big\lfloor 1, \big(\tfrac{1}{5}, \tfrac{3}{5}, \tfrac{3}{5}, \tfrac{1}{5}, \tfrac{2}{5}\big) \big\rceil$\tsep{1pt}\\
$\big\lfloor x_1x_2^2x_5^2 + x_1^2x_2x_5^2 + x_3x_4^2x_5^2 + x_3^2x_4x_5^2, (0, 0, 0, 0, 0) \big\rceil   $& $\big\lfloor 1, \big(\tfrac{2}{5}, \tfrac{3}{5}, \tfrac{1}{5}, \tfrac{1}{5}, \tfrac{3}{5}\big) \big\rceil + (3 \text{ others})$\tsep{1pt}\\
$\big\lfloor x_1x_3^2x_5^2 + x_1^2x_3x_5^2 + x_2x_4^2x_5^2 + x_2^2x_4x_5^2, (0, 0, 0, 0, 0) \big\rceil   $& $\big\lfloor 1, \big(\tfrac{2}{5}, \tfrac{1}{5}, \tfrac{3}{5}, \tfrac{1}{5}, \tfrac{3}{5}\big) \big\rceil + (3 \text{ others})$\tsep{1pt}\\
$\big\lfloor x_1x_4^2x_5^2 + x_1^2x_4x_5^2 + x_2x_3^2x_5^2 + x_2^2x_3x_5^2, (0, 0, 0, 0, 0) \big\rceil   $& $\big\lfloor 1, \big(\tfrac{2}{5}, \tfrac{1}{5}, \tfrac{1}{5}, \tfrac{3}{5}, \tfrac{3}{5}\big) \big\rceil + (3 \text{ others})$\tsep{1pt}\\
$\big\lfloor x_1^2x_2x_3x_4 + x_1x_2^2x_3x_4 + x_1x_2x_3^2x_4 + x_1x_2x_3x_4^2, (0, 0, 0, 0, 0) \big\rceil $& $\big\lfloor 1, \big(\tfrac{3}{5}, \tfrac{2}{5}, \tfrac{2}{5}, \tfrac{2}{5}, \tfrac{1}{5}\big) \big\rceil + (3 \text{ others})$\tsep{1pt}\\
$\big\lfloor x_1^2x_2x_3x_5 + x_1x_2^2x_4x_5 + x_1x_3^2x_4x_5 + x_2x_3x_4^2x_5, (0, 0, 0, 0, 0) \big\rceil $& $\big\lfloor 1, \big(\tfrac{3}{5}, \tfrac{2}{5}, \tfrac{2}{5}, \tfrac{1}{5}, \tfrac{2}{5}\big) \big\rceil + (3 \text{ others})$\tsep{1pt}\\
$\big\lfloor x_1x_2^2x_3x_5 + x_1^2x_2x_4x_5 + x_1x_3x_4^2x_5 + x_2x_3^2x_4x_5, (0, 0, 0, 0, 0) \big\rceil $& $\big\lfloor 1, \big(\tfrac{2}{5}, \tfrac{3}{5}, \tfrac{1}{5}, \tfrac{2}{5}, \tfrac{2}{5}\big) \big\rceil + (3 \text{ others})$\tsep{1pt}\\
$\big\lfloor x_1x_2x_3^2x_5 + x_1x_2x_4^2x_5 + x_1^2x_3x_4x_5 + x_2^2x_3x_4x_5, (0, 0, 0, 0, 0) \big\rceil $& $\big\lfloor 1, \big(\tfrac{2}{5}, \tfrac{2}{5}, \tfrac{3}{5}, \tfrac{1}{5}, \tfrac{2}{5}\big) \big\rceil + (3 \text{ others})$\tsep{1pt}\\
$\big\lfloor x_1x_2x_3x_5^2 + x_1x_2x_4x_5^2 + x_1x_3x_4x_5^2 + x_2x_3x_4x_5^2, (0, 0, 0, 0, 0) \big\rceil $& $\big\lfloor 1, \big(\tfrac{2}{5}, \tfrac{2}{5}, \tfrac{2}{5}, \tfrac{1}{5}, \tfrac{3}{5}\big) \big\rceil + (3 \text{ others})$\tsep{1pt}\\
\bottomrule 
\caption{In our bad example, the mirror map works on these basis elements of degree $(1,2)$ by Theorem~\ref{thm:mirror}.}\label{tab:baddegree12}
\end{longtable}}

The following $28$ corresponding elements in Table~\ref{tab:baddegree21} have a bidegree of $(2, 1)$. As with the previous page, the basis elements on the A-side come from the untwisted broad sector and the elements on the B side are narrow, also following the recipe from Theorem~\ref{thm:mirror}.

{\footnotesize	
\begin{longtable}{@{\,}ll@{\,}} \toprule
\multicolumn{2}{c}{Basis elements of degree $(2,1)$} \\
{A-model} &{B-model}	\\
	\midrule
	\endfirsthead
	\toprule
{A-model} &{B-model}	\\
	 \midrule
	\endhead
$\big\lfloor x_1^3x_2^3x_3^3x_4 + x_1^3x_2^3x_3x_4^3 + x_1^3x_2x_3^3x_4^3 + x_1x_2^3x_3^3x_4^3, (0, 0, 0, 0, 0) \big\rceil  $ & $\big\lfloor 1, \big(\tfrac{4}{5}, \tfrac{4}{5}, \tfrac{4}{5}, \tfrac{2}{5}, \tfrac{1}{5}\big) \big\rceil + (3 \text{ others})$\tsep{1pt}\\
$\big\lfloor x_1^3x_2^3x_3^3x_5 + x_1^3x_2^3x_4^3x_5 + x_1^3x_2^3x_4^3x_5 + x_2^3x_3^3x_4^3x_5, (0, 0, 0, 0, 0) \big\rceil  $ & $\big\lfloor 1, \big(\tfrac{4}{5}, \tfrac{4}{5}, \tfrac{4}{5}, \tfrac{1}{5}, \tfrac{2}{5}\big) \big\rceil + (3 \text{ others})$\tsep{1pt}\\
$\big\lfloor x_1^3x_2^3x_3x_5^3 + x_1^3x_2^3x_4x_5^3 + x_1x_3^3x_4^3x_5^3 + x_2x_3^3x_4^3x_5^3, (0, 0, 0, 0, 0) \big\rceil  $ & $\big\lfloor 1, \big(\tfrac{4}{5}, \tfrac{4}{5}, \tfrac{2}{5}, \tfrac{1}{5}, \tfrac{4}{5}\big) \big\rceil + (3 \text{ others})$\tsep{1pt}\\
$\big\lfloor x_1^3x_2x_3^3x_5^3 + x_1x_2^3x_4^3x_5^3 + x_1^3x_3^3x_4x_5^3 + x_2^3x_3x_4^3x_5^3, (0, 0, 0, 0, 0) \big\rceil  $ & $\big\lfloor 1, \big(\tfrac{4}{5}, \tfrac{2}{5}, \tfrac{4}{5}, \tfrac{1}{5}, \tfrac{4}{5}\big) \big\rceil + (3 \text{ others})$\tsep{1pt}\\
$\big\lfloor x_1x_2^3x_3^3x_5^3 + x_1^3x_2x_4^3x_5^3 + x_1^3x_3x_4^3x_5^3 + x_2^3x_3^3x_4x_5^3, (0, 0, 0, 0, 0) \big\rceil  $ & $\big\lfloor 1, \big(\tfrac{2}{5}, \tfrac{4}{5}, \tfrac{4}{5}, \tfrac{1}{5}, \tfrac{4}{5}\big) \big\rceil + (3 \text{ others})$\tsep{1pt}\\
$\big\lfloor x_1^3x_2^3x_3^2x_4^2 + x_1^2x_2^2x_3^3x_4^3, (0, 0, 0, 0, 0) \big\rceil  $ & $\big\lfloor 1, \big(\tfrac{4}{5}, \tfrac{4}{5}, \tfrac{3}{5}, \tfrac{3}{5}, \tfrac{1}{5}\big) \big\rceil + \big\lfloor 1, \big(\tfrac{3}{5}, \tfrac{3}{5}, \tfrac{4}{5}, \tfrac{4}{5}, \tfrac{1}{5}\big) \big\rceil$\tsep{1pt}\\
$\big\lfloor x_1^3x_2^2x_3^3x_4^2 + x_1^2x_2^3x_3^2x_4^3, (0, 0, 0, 0, 0) \big\rceil  $ & $\big\lfloor 1, \big(\tfrac{4}{5}, \tfrac{3}{5}, \tfrac{4}{5}, \tfrac{3}{5}, \tfrac{1}{5}\big) \big\rceil + \big\lfloor 1, \big(\tfrac{3}{5}, \tfrac{4}{5}, \tfrac{3}{5}, \tfrac{4}{5}, \tfrac{1}{5}\big) \big\rceil$\tsep{1pt}\\
$\big\lfloor x_1^3x_2^2x_3^2x_4^3 + x_1^2x_2^3x_3^3x_4^2, (0, 0, 0, 0, 0) \big\rceil  $ & $\big\lfloor 1, \big(\tfrac{4}{5}, \tfrac{3}{5}, \tfrac{3}{5}, \tfrac{4}{5}, \tfrac{1}{5}\big) \big\rceil + \big\lfloor 1, \big(\tfrac{3}{5}, \tfrac{4}{5}, \tfrac{4}{5}, \tfrac{3}{5}, \tfrac{1}{5}\big) \big\rceil$\tsep{1pt}\\
$\big\lfloor x_1^3x_2^3x_3^2x_5^2 + x_1^3x_2^3x_4^2x_5^2 + x_1^2x_3^3x_4^3x_5^2 + x_2^2x_3^3x_4^3x_5^2, (0, 0, 0, 0, 0) \big\rceil  $ & $\big\lfloor 1, \big(\tfrac{4}{5}, \tfrac{4}{5}, \tfrac{3}{5}, \tfrac{1}{5}, \tfrac{3}{5}\big) \big\rceil + (3 \text{ others})$\tsep{1pt}\\
$\big\lfloor x_1^3x_2^2x_3^3x_5^2 + x_1^2x_2^3x_4^3x_5^2 + x_1^3x_3^3x_4^2x_5^2 + x_2^3x_3^2x_4^3x_5^2, (0, 0, 0, 0, 0) \big\rceil  $ & $\big\lfloor 1, \big(\tfrac{4}{5}, \tfrac{2}{5}, \tfrac{4}{5}, \tfrac{1}{5}, \tfrac{3}{5}\big) \big\rceil + (3 \text{ others})$\tsep{1pt}\\
$\big\lfloor x_1^2x_2^3x_3^3x_5^2 + x_1^3x_2^2x_4^3x_5^2 + x_1^3x_3^2x_4^3x_5^2 + x_2^3x_3^3x_4^2x_5^2, (0, 0, 0, 0, 0) \big\rceil  $ & $\big\lfloor 1, \big(\tfrac{3}{5}, \tfrac{4}{5}, \tfrac{4}{5}, \tfrac{1}{5}, \tfrac{3}{5}\big) \big\rceil + (3 \text{ others})$\tsep{1pt}\\
$\big\lfloor x_1^3x_2^2x_3^2x_5^3 + x_1^2x_2^3x_4^2x_5^3 + x_1^2x_3^3x_4^2x_5^3 + x_2^2x_3^2x_4^3x_5^3, (0, 0, 0, 0, 0) \big\rceil  $ & $\big\lfloor 1, \big(\tfrac{4}{5}, \tfrac{4}{5}, \tfrac{3}{5}, \tfrac{1}{5}, \tfrac{4}{5}\big) \big\rceil + (3 \text{ others})$\tsep{1pt}\\
$\big\lfloor x_1^2x_2^3x_3^2x_5^3 + x_1^3x_2^2x_4^2x_5^3 + x_1^2x_3^2x_4^3x_5^3 + x_2^2x_3^3x_4^2x_5^3, (0, 0, 0, 0, 0) \big\rceil  $ & $\big\lfloor 1, \big(\tfrac{3}{5}, \tfrac{4}{5}, \tfrac{3}{5}, \tfrac{1}{5}, \tfrac{4}{5}\big) \big\rceil + (3 \text{ others})$\tsep{1pt}\\
$\big\lfloor x_1^2x_2^2x_3^3x_5^3 + x_1^2x_2^2x_4^3x_5^3 + x_1^3x_3^2x_4^2x_5^3 + x_2^3x_3^2x_4^2x_5^3, (0, 0, 0, 0, 0) \big\rceil  $ & $\big\lfloor 1, \big(\tfrac{3}{5}, \tfrac{3}{5}, \tfrac{4}{5}, \tfrac{1}{5}, \tfrac{4}{5}\big) \big\rceil + (3 \text{ others})$\tsep{1pt}\\
$\big\lfloor x_1^3x_2^2x_3^2x_4^2x_5 + x_1^2x_2^3x_3^2x_4^2x_5 + x_1^2x_2^2x_3^3x_4^2x_5 + x_1^2x_2^2x_3^2x_4^3x_5, (0, 0, 0, 0, 0) \big\rceil  $ & $\big\lfloor 1, \big(\tfrac{4}{5}, \tfrac{3}{5}, \tfrac{3}{5}, \tfrac{3}{5}, \tfrac{2}{5}\big) \big\rceil + (3 \text{ others})$\tsep{1pt}\\
$\big\lfloor x_1^3x_2^2x_3^2x_4x_5^2 + x_1^2x_2^3x_3x_4^2x_5^2 + x_1^2x_2x_3^3x_4^2x_5^2 + x_1x_2^2x_3^2x_4^3x_5^2, (0, 0, 0, 0, 0) \big\rceil  $ & $\big\lfloor 1, \big(\tfrac{4}{5}, \tfrac{3}{5}, \tfrac{3}{5}, \tfrac{2}{5}, \tfrac{3}{5}\big) \big\rceil + (3 \text{ others})$\tsep{1pt}\\
$\big\lfloor x_1^3x_2^2x_3x_4^2x_5^2 + x_1^2x_2^3x_3^2x_4x_5^2 + x_1x_2^2x_3^3x_4^2x_5^2 + x_1^2x_2x_3^2x_4^3x_5^2, (0, 0, 0, 0, 0) \big\rceil  $ & $\big\lfloor 1, \big(\tfrac{4}{5}, \tfrac{3}{5}, \tfrac{2}{5}, \tfrac{3}{5}, \tfrac{3}{5}\big) \big\rceil + (3 \text{ others})$\tsep{1pt}\\
$\big\lfloor x_1^3x_2x_3^2x_4^2x_5^2 + x_1x_2^3x_3^2x_4^2x_5^2 + x_1^2x_2^2x_3^3x_4x_5^2 + x_1^2x_2^2x_3x_4^3x_5^2, (0, 0, 0, 0, 0) \big\rceil  $ & $\big\lfloor 1, \big(\tfrac{4}{5}, \tfrac{2}{5}, \tfrac{3}{5}, \tfrac{3}{5}, \tfrac{3}{5}\big) \big\rceil + (3 \text{ others})$\tsep{1pt}\\
$\big\lfloor x_1^2x_2^2x_3^2x_4x_5^3 + x_1^2x_2^2x_3^1x_4^2x_5^3 + x_1^2x_2x_3^2x_4^2x_5^3 + x_1x_2^2x_3^2x_4^2x_5^3, (0, 0, 0, 0, 0) \big\rceil  $ & $\big\lfloor 1, \big(\tfrac{3}{5}, \tfrac{3}{5}, \tfrac{3}{5}, \tfrac{2}{5}, \tfrac{4}{5}\big) \big\rceil + (3 \text{ others})$\tsep{1pt}\\
$\big\lfloor x_1^3x_2^3x_3^2x_4x_5 + x_1^3x_2^3x_3x_4^2x_5 + x_1^2x_2x_3^3x_4^3x_5 + x_1x_2^2x_3^3x_4^3x_5, (0, 0, 0, 0, 0) \big\rceil  $ & $\big\lfloor 1, \big(\tfrac{4}{5}, \tfrac{4}{5}, \tfrac{3}{5}, \tfrac{2}{5}, \tfrac{2}{5}\big) \big\rceil + (3 \text{ others})$\tsep{1pt}\\
$\big\lfloor x_1^3x_2^2x_3^3x_4x_5 + x_1^2x_2^3x_3x_4^3x_5 + x_1^3x_2x_3^3x_4^2x_5 + x_1x_2^3x_3^2x_4^3x_5, (0, 0, 0, 0, 0) \big\rceil  $ & $\big\lfloor 1, \big(\tfrac{4}{5}, \tfrac{3}{5}, \tfrac{4}{5}, \tfrac{2}{5}, \tfrac{2}{5}\big) \big\rceil + (3 \text{ others})$\tsep{1pt}\\
$\big\lfloor x_1^3x_2^2x_3x_4^3x_5 + x_1^2x_2^3x_3^3x_4x_5 + x_1x_2^3x_3^3x_4^2x_5 + x_1^3x_2x_3^2x_4^3x_5, (0, 0, 0, 0, 0) \big\rceil  $ & $\big\lfloor 1, \big(\tfrac{4}{5}, \tfrac{3}{5}, \tfrac{2}{5}, \tfrac{4}{5}, \tfrac{2}{5}\big) \big\rceil + (3 \text{ others})$\tsep{1pt}\\
$\big\lfloor x_1^3x_2^3x_3x_4x_5^2 + x_1x_2x_3^3x_4^3x_5^2, (0, 0, 0, 0, 0) \big\rceil  $ & $\big\lfloor 1, \big(\tfrac{4}{5}, \tfrac{4}{5}, \tfrac{2}{5}, \tfrac{2}{5}, \tfrac{3}{5}\big) \big\rceil + \big\lfloor 1, \big(\tfrac{2}{5}, \tfrac{2}{5}, \tfrac{4}{5}, \tfrac{4}{5}, \tfrac{3}{5}\big) \big\rceil$\tsep{1pt}\\
$\big\lfloor x_1^3x_2x_3^3x_4x_5^2 + x_1x_2^3x_3x_4^3x_5^2, (0, 0, 0, 0, 0) \big\rceil  $ & $\big\lfloor 1, \big(\tfrac{4}{5}, \tfrac{2}{5}, \tfrac{4}{5}, \tfrac{2}{5}, \tfrac{3}{5}\big) \big\rceil + \big\lfloor 1, \big(\tfrac{2}{5}, \tfrac{4}{5}, \tfrac{2}{5}, \tfrac{4}{5}, \tfrac{3}{5}\big) \big\rceil$\tsep{1pt}\\
$\big\lfloor x_1^3x_2x_3x_4^3x_5^2 + x_1x_2^3x_3^3x_4x_5^2, (0, 0, 0, 0, 0) \big\rceil  $ & $\big\lfloor 1, \big(\tfrac{4}{5}, \tfrac{2}{5}, \tfrac{2}{5}, \tfrac{4}{5}, \tfrac{3}{5}\big) \big\rceil + \big\lfloor 1, \big(\tfrac{2}{5}, \tfrac{4}{5}, \tfrac{4}{5}, \tfrac{2}{5}, \tfrac{3}{5}\big) \big\rceil$\tsep{1pt}\\
$\big\lfloor x_1^3x_2^2x_3x_4x_5^3 + x_1^2x_2^3x_3x_4x_5^3 + x_1x_2x_3^3x_4^2x_5^3 + x_1x_2x_3^2x_4^3x_5^3, (0, 0, 0, 0, 0) \big\rceil  $ & $\big\lfloor 1, \big(\tfrac{4}{5}, \tfrac{3}{5}, \tfrac{2}{5}, \tfrac{2}{5}, \tfrac{4}{5}\big) \big\rceil + (3 \text{ others})$\tsep{1pt}\\
$\big\lfloor x_1^3x_2x_3^2x_4x_5^3 + x_1x_2^3x_3x_4^2x_5^3 + x_1^2x_2x_3^3x_4x_5^3 + x_1x_2^2x_3x_4^3x_5^3, (0, 0, 0, 0, 0) \big\rceil  $ & $\big\lfloor 1, \big(\tfrac{4}{5}, \tfrac{2}{5}, \tfrac{3}{5}, \tfrac{2}{5}, \tfrac{4}{5}\big) \big\rceil + (3 \text{ others})$\tsep{1pt}\\
$\big\lfloor x_1^3x_2x_3x_4^2x_5^3 + x_1x_2^3x_3^2x_4x_5^3 + x_1x_2^2x_3^3x_4x_5^3 + x_1^2x_2x_3x_4^3x_5^3, (0, 0, 0, 0, 0) \big\rceil  $ & $\big\lfloor 1, \big(\tfrac{4}{5}, \tfrac{2}{5}, \tfrac{2}{5}, \tfrac{3}{5}, \tfrac{4}{5}\big) \big\rceil + (3 \text{ others})$\tsep{1pt}\\
\bottomrule 
\caption{In our bad example, the mirror map works on these basis elements of degree $(2,1)$ by Theorem~\ref{thm:mirror}.}\label{tab:baddegree21}
\end{longtable}}	

While everything has followed as expected in the previous tables, the mirror map doesn't hold when we look at sectors generated by permutation matrices, as we can see in the following table. The total number of basis elements is the same on both the A- and B-side, but there is no way to make the bidegrees match up. In the bottom half of the table, we have made an appropriate change of variables, in order to express the polynomials in a more succinct manner. For example, in the sector indexed by $(12)(34)$, we make the change of variables $y_1=x_1+x_2$, $y_2=x_3+x_4$ and $y_3=x_5$, and in the sector indexed by $(13)(24)$, we use the (slightly abusive notation) change of variables $y_1=x_1+x_3$, $y_2=x_2+x_4$ and $y_3=x_5$, etc.

{\small	
\begin{longtable}{llll} \toprule
{bidegree} &{A-model} & {bidegree} &{B-model}	\\
	\midrule
	\endfirsthead
	\toprule
{bidegree} &{A-model} & {bidegree} &{B-model}	\\
	 \midrule
	\endhead
$(1, 1)$ &$\lfloor 1, ((12)(34))j_W \rceil $& $(1, 2) $& $\lfloor 1, ((12)(34))j_W \rceil + (24 \text{ others})$\tsep{1pt}\\
$(2, 2) $& $\big\lfloor 1, ((12)(34))(j_W)^2 \big\rceil $& $(2, 1) $& $\big\lfloor 1, ((12)(34))(j_W)^2 \big\rceil + (24 \text{ others})$\tsep{1pt}\\
$(1, 1) $& $\big\lfloor 1, ((12)(34))(j_W)^3 \big\rceil $& $(1, 2) $& $\big\lfloor 1, ((12)(34))(j_W)^3 \big\rceil + (24 \text{ others})$\tsep{1pt}\\
$(2, 2) $& $\big\lfloor 1, ((12)(34))(j_W)^4 \big\rceil $& $(2, 1) $& $\big\lfloor 1, ((12)(34))(j_W)^4 \big\rceil + (24 \text{ others})$\tsep{1pt}\\
$(1, 1) $& $\lfloor 1, ((13)(24))j_W \rceil $& $(1, 2) $& $\lfloor 1, ((13)(24))j_W \rceil + (24 \text{ others})$\tsep{1pt}\\
$(2, 2) $& $\big\lfloor 1, ((13)(24))(j_W)^2 \big\rceil $& $(2, 1) $& $\big\lfloor 1, ((13)(24))(j_W)^2 \big\rceil + (24 \text{ others})$\tsep{1pt}\\
$(1, 1) $& $\big\lfloor 1, ((13)(24))(j_W)^3 \big\rceil $& $(1, 2) $& $\big\lfloor 1, ((13)(24))(j_W)^3 \big\rceil + (24 \text{ others})$\tsep{1pt}\\
$(2, 2) $& $\big\lfloor 1, ((13)(24))(j_W)^4 \big\rceil $& $(2, 1) $& $\big\lfloor 1, ((13)(24))(j_W)^4 \big\rceil + (24 \text{ others})$\tsep{1pt}\\
$(1, 1) $& $\lfloor 1, ((14)(23))j_W \rceil $& $(1, 2) $& $\lfloor 1, ((14)(23))j_W \rceil + (24 \text{ others})$\tsep{1pt}\\
$(2, 2) $& $\big\lfloor 1, ((14)(23))(j_W)^2 \big\rceil $& $(2, 1) $& $\big\lfloor 1, ((14)(23))(j_W)^2 \big\rceil + (24 \text{ others})$\tsep{1pt}\\
$(1, 1) $& $\big\lfloor 1, ((14)(23))(j_W)^3 \big\rceil $& $(1, 2) $& $\big\lfloor 1, ((14)(23))(j_W)^3 \big\rceil + (24 \text{ others})$\tsep{1pt}\\
$(2, 2) $& $\big\lfloor 1, ((14)(23))(j_W)^4 \big\rceil $& $(2, 1) $& $\big\lfloor 1, ((14)(23))(j_W)^4 \big\rceil + (24 \text{ others})$\tsep{1pt}\\
$(1, 2) $& $\big\lfloor y_1^2-y_2^2, (12)(34) \big\rceil $& $(1, 2) $& $\big\lfloor 1, ((12)(34))j_W^2K \big\rceil + (49 \text{ others})$\tsep{1pt}\\
$(1, 2) $& $\lfloor (y_1-y_2)y_3, (12)(34) \rceil $& $(1, 2) $& $\big\lfloor 1, ((12)(34))j_W^2L \big\rceil + (49 \text{ others})$\tsep{1pt}\\
$(2, 1) $& $\big\lfloor \big(y_1^3y_2-y_1y_2^3\big)y_3^3, (12)(34) \big\rceil $& $(2, 1)$ & $\big\lfloor 1, ((12)(34))j_W^3K \big\rceil + (49 \text{ others})$\tsep{1pt}\\
$(2, 1) $& $\big\lfloor \big(y_1^3y_2^2-y_1^2y_2^3\big)y_3^2, (12)(34) \big\rceil $& $(2, 1) $& $\big\lfloor 1, ((12)(34))j_W^3L \big\rceil + (49 \text{ others})$\tsep{1pt}\\
$(1, 2) $& $\big\lfloor y_1^2-y_2^2, (13)(24) \big\rceil $& $(1, 2) $& $\big\lfloor 1, ((13)(24))j_W^2K \big\rceil + (49 \text{ others})$\tsep{1pt}\\
$(1, 2) $& $\lfloor (y_1-y_2)y_3, (13)(24) \rceil $& $(1, 2) $& $\big\lfloor 1, ((13)(24))j_W^2L \big\rceil + (49 \text{ others})$\tsep{1pt}\\
$(2, 1) $& $\big\lfloor \big(y_1^3y_2-y_1y_2^3\big)y_3^3, (13)(24) \big\rceil $& $(2, 1) $& $\big\lfloor 1, ((13)(24))j_W^3K \big\rceil + (49 \text{ others})$\tsep{1pt}\\
$(2, 1) $& $\big\lfloor \big(y_1^3y_2^2-y_1^2y_2^3\big)y_3^2, (13)(24) \big\rceil $& $(2, 1) $& $\big\lfloor 1, ((13)(24))j_W^3L \big\rceil + (49 \text{ others})$\tsep{1pt}\\
$(1, 2) $& $\big\lfloor y_1^2-y_2^2, (14)(23) \big\rceil\ $& $(1, 2) $& $\big\lfloor 1, ((14)(23))j_W^2K \big\rceil + (49 \text{ others})$\tsep{1pt}\\
$(1, 2) $& $\lfloor (y_1-y_2)y_3, (14)(23) \rceil $& $(1, 2) $& $\big\lfloor 1, ((14)(23))j_W^2L \big\rceil + (49 \text{ others})$\tsep{1pt}\\
$(2, 1) $& $\big\lfloor \big(y_1^3y_2-y_1y_2^3\big)y_3^3, (14)(23) \big\rceil $& $(2, 1) $& $\big\lfloor 1, ((14)(23))j_W^3K \big\rceil + (49 \text{ others})$\tsep{1pt}\\
$(2, 1) $& $\big\lfloor \big(y_1^3y_2^2-y_1^2y_2^3\big)y_3^2, (14)(23) \big\rceil $& $(2, 1) $& $\big\lfloor 1, ((14)(23))j_W^3L \big\rceil + (49 \text{ others})$\tsep{1pt}\\
\bottomrule 
\caption{In our bad example, the mirror map works on these basis elements of degree $(2,1)$ by Theorem~\ref{thm:mirror}. In the sector indexed by $(12)(34)$, we make the change of variables $y_1=x_1+x_2$, $y_2=x_3+x_4$ and $y_3=x_5$; in the sector indexed by $(13)(24)$, we use the (slightly abusive notation) change of variables $y_1=x_1+x_3$, $y_2=x_2+x_4$ and $y_3=x_5$; and in the sector indexed by $(13)(24)$, we use the (slightly abusive notation) change of variables $y_1=x_1+x_4$, $y_2=x_2+x_3$ and $y_3=x_5$.}\label{tab:baddeg21}
\end{longtable}} 

Thus we see exactly in these places, the mirror map does not hold in full generality. So although the restricted mirror map is always an isomorphism, the entire state spaces are not always isomorphic as bigraded vector spaces.
\end{Example}

\subsection*{Acknowledgments}
The authors would like to thank Tyler Jarvis, Yongbin Ruan, and Yefeng Shen for their helpful remarks on early results. We would also like to thank the anonymous referees whose suggestions have greatly improved this article.

\pdfbookmark[1]{References}{ref}
\LastPageEnding

\end{document}